\documentclass[11pt]{amsart}
\usepackage{amssymb,mathrsfs}
\title{Harmonic Analysis on Homogeneous Spaces}

\author{Jae-Hyun Yang}
\address{Department of Mathematics, Inha University,
Incheon 402-751, Korea}
\email{jhyang@inha.ac.kr }

\begin{document}

\newtheorem{theorem}{Theorem}[section]
\newtheorem{lemma}{Lemma}[section]
\newtheorem{proposition}{Proposition}[section]
\newtheorem{remark}{Remark}[section]
\newtheorem{definition}{Definition}[section]

\renewcommand{\theequation}{\thesection.\arabic{equation}}
\renewcommand{\thetheorem}{\thesection.\arabic{theorem}}
\renewcommand{\thelemma}{\thesection.\arabic{lemma}}
\newcommand{\BR}{\mathbb R}
\newcommand{\BQ}{\mathbb Q}
\newcommand{\bn}{\bf n}
\def\charf {\mbox{{\text 1}\kern-.24em {\text l}}}
\newcommand{\BC}{\mathbb C}
\newcommand{\BZ}{\mathbb Z}



\begin{abstract}
This article is an expository paper. We first survey developments
over the past three decades in the theory of harmonic analysis on
reductive symmetric spaces. Next we deal with the particular
homogeneous space of non-reductive type, the so called
Siegel-Jacobi space that is important arithmetically and
geometrically. We present some new results on the Siegel-Jacobi
space.
\end{abstract}
\maketitle


\newcommand\tr{\triangleright}
\newcommand\al{\alpha}
\newcommand\g{\gamma}
\newcommand\gh{\Cal G^J}
\newcommand\G{\Gamma}
\newcommand\de{\delta}
\newcommand\e{\epsilon}
\newcommand\z{\zeta}
\newcommand\vth{\vartheta}
\newcommand\vp{\varphi}
\newcommand\om{\omega}
\newcommand\p{\pi}
\newcommand\la{\lambda}
\newcommand\lb{\lbrace}
\newcommand\lk{\lbrack}
\newcommand\rb{\rbrace}
\newcommand\rk{\rbrack}
\newcommand\s{\sigma}
\newcommand\w{\wedge}
\newcommand\fgj{{\frak g}^J}
\newcommand\lrt{\longrightarrow}
\newcommand\lmt{\longmapsto}
\newcommand\lmk{(\lambda,\mu,\kappa)}
\newcommand\Om{\Omega}
\newcommand\ka{\kappa}
\newcommand\ba{\backslash}
\newcommand\ph{\phi}
\newcommand\M{{\Cal M}}
\newcommand\bA{\bold A}
\newcommand\bH{\bold H}

\newcommand\Hom{\text{Hom}}
\newcommand\cP{\Cal P}
\newcommand\cH{\Cal H}

\newcommand\pa{\partial}

\newcommand\pis{\pi i \sigma}
\newcommand\sd{\,\,{\vartriangleright}\kern -1.0ex{<}\,}
\newcommand\wt{\widetilde}
\newcommand\fg{\frak g}
\newcommand\fk{\frak k}
\newcommand\fp{\frak p}
\newcommand\fs{\frak s}
\newcommand\fh{\frak h}
\newcommand\fa{\frak a}
\newcommand\fb{\frak b}
\newcommand\fq{\frak q}
\newcommand\fm{\frak m}
\newcommand\fn{\frak n}
\newcommand\Cal{\mathcal}
\newcommand\CH{\mathcal H}
\newcommand\wf{\widehat f}

\newcommand\CP{{\mathcal P}_n}
\newcommand\Hnm{{\mathbb H}_n \times {\mathbb C}^{(m,n)}}
\newcommand\Dnm{{\mathbb D}_n \times {\mathbb C}^{(m,n)}}
\newcommand\BD{\mathbb D}
\newcommand\BH{\mathbb H}
\newcommand\BHn{\BH_n}
\newcommand\BDn{\BD_n}

\newcommand\CCF{{\mathcal F}_n}
\newcommand\CM{{\mathcal M}}
\newcommand\CHX{{\mathcal H}_{\xi}}
\newcommand\Ggh{\Gamma_{g,h}}
\newcommand\Cmn{{\mathbb C}^{(m,n)}}
\newcommand\Yd{{{\partial}\over {\partial Y}}}
\newcommand\Vd{{{\partial}\over {\partial V}}}

\newcommand\Ys{Y^{\ast}}
\newcommand\Vs{V^{\ast}}
\newcommand\LO{L_{\Omega}}
\newcommand\fac{{\frak a}_{\mathbb C}^{\ast}}
\newcommand\OW{\overline{W}}
\newcommand\OP{\overline{P}}
\newcommand\OQ{\overline{Q}}
\newcommand\Dg{{\mathbb D}_g}
\newcommand\Hg{{\mathbb H}_g}
\newcommand\wG{\widehat G}

\newcommand\td{\bigtriangledown}
\newcommand\pdx{ {{\partial}\over{\partial x}} }
\newcommand\pdy{ {{\partial}\over{\partial y}} }
\newcommand\pdu{ {{\partial}\over{\partial u}} }
\newcommand\pdv{ {{\partial}\over{\partial v}} }
\newcommand\PZ{ {{\partial}\over {\partial Z}} }
\newcommand\PW{ {{\partial}\over {\partial W}} }
\newcommand\PZB{ {{\partial}\over {\partial{\overline Z}}} }
\newcommand\PWB{ {{\partial}\over {\partial{\overline W}}} }
\newcommand\PX{ {{\partial}\over{\partial X}} }
\newcommand\PY{ {{\partial}\over {\partial Y}} }
\newcommand\PU{ {{\partial}\over{\partial U}} }
\newcommand\PV{ {{\partial}\over{\partial V}} }
\renewcommand\th{\theta}
\renewcommand\l{\lambda}
\renewcommand\k{\kappa}
\newcommand\PWE{ \frac{\partial}{\partial \overline W}}

\def\charf {\mbox{{\text 1}\kern-.24em {\text l}}}

\vskip 1cm
%
%
\begin{section}{{\bf Introduction}}
\setcounter{equation}{0} \vskip 0.2cm  Let $G$ be a connected real
Lie group and let $H$ be a closed subgroup of $G$. A homogeneous
space is understood as a manifold with a transitive action of a
Lie group. We see that the coset space $G/H$ is a homogeneous
space. If $G/H$ carries a $G$-invariant measure and $H$ has at
most finitely many connected components, then we have the $
\textit{natural}$ unitary representation $\pi_H$ of $G$ on the
Hilbert space $L^2(G/H)$ of all square integrable functions on
$G/H$ defined by
\begin{equation}
\left( \pi_H(g)f\right)(x)=f(g^{-1}x), \end{equation} where $g\in
G,\ x\in G/H$ and $f\in L^2(G/H).$ $\pi_H$ is called the (left) $
\textit{regular representation}$ of $G$ for the homogeneous space
$G/H$. The fundamental problem in harmonic analysis on $G/H$ is
given as follows.

\vskip 0.3cm  {\bf Problem 1.} Decompose the regular
representation $\pi_H$ into irreducible unitary representations of
$G$.

\vskip 0.3cm We now assume that $\G$ is an arithmetic subgroup of
$G$. The right regular representation $\pi_\G$ of $G$ on the
Hilbert space $L^2(\G\ba G)$ of square integrable functions on
$\G\ba G$ is given by
\begin{equation}
\left( \pi_{\G}(g)f\right)(z)=f(zg), \end{equation} where $g\in
G,\ z\in \G\ba G$ and $f\in L^2(\G\ba G).$ Obviously $\pi_\G$ is a
unitary representation of $G$. Then we have a basic problem in
harmonic analysis on $\G\ba G.$

\vskip 0.3cm  {\bf Problem 2.} Find an explicit formula of the
irreducible decomposition of $\pi_\G.$ \vskip 0.3cm Problem 2 is
closely related to the theory of automorphic forms on $G$ for $\G$
and also to the geometry of the associated arithmetic variety. In
order to answer the above problems, we need to classify
irreducible unitary representations of $G.$ Let $\widehat G$ be
the unitary dual of $G$, that is, the set of all the equivalence
classes of irreducible unitary representations of $G$. We denote
by $\textrm{Disc}(G/H)$\,(resp. $ \textrm{Disc}(\G\ba G)$) the set
of discrete series representations of $G$ for $G/H$\,(resp. $\G\ba
G$). The first step for solving the above problems might be the
determination of the subsets $\textrm{Disc}(G/H)$ and
$\textrm{Disc}(\G\ba G)$ of ${\widehat G}$. For instance, it may
be meaningful to find a criterion under which
$\textrm{Disc}(G/H)$\,(resp. $ \textrm{Disc}(\G\ba G)$) is
non-empty. In fact, the study of the unitary dual $\widehat G$ has
served as a basic tool for the solution of the above problems. We
can say that Problem 1 is essentially equivalent to finding an
explicit direct integral decomposition
\begin{equation}
\pi_H\cong \int_{\wG}m_{\pi}\,\pi\, d\mu(\pi)
\end{equation}
of $\pi_H$ into irreducible unitary representations of $G$, where
$d\mu$ is a suitable measure (Plancherel measure) on $\wG$. If $G$
is a connected reductive real Lie group, it can be regarded as a
symmetric space $(G\times G)/ G$ for the left times right action
of $G\times G$. As is well known, in this case (called the {\it
group case}) an explicit decomposition of the form (1.3) was
determined by Harish-Chandra [26,\,27,\,28]. The components of
$L^2(G)$ fall in a finite number of series, each of which
corresponds to a particular cuspidal parabolic subgroup $P$ and is
a direct integral of $\pi_{P,\xi,\la}\otimes \pi_{P,\xi,\la}^*,$
where $\pi_{P,\xi,\la}$ belongs to the principal series induced
from $P=MAN,\ \xi$ is a discrete series representation of $M$ and
$\la$ is a continuous parameter ranging in the imaginary linear
dual $i\,\fa$ of the Lie algebra $\fa$ of $A.$ The unitary dual
$\wG$ has been investigated by Harish-Chandra, Gelfand school and
others. So far the classification of $\wG$ has not been solved
completely except for some special groups. Recently van den Ban
and Schlichkrull [13,\,14] obtained an explicit Plancherel formula
for a reductive symmetric space, which will be explained in
Section 4 in some detail. On the other hand, P. Delorme [21] gave
the Plancherel formula for a reductive symmetric space by a method
different from that of van den Ban and Schlichtkrull. In his paper
[19,\,20], the Maass-Selberg relations for the Eisenstein
integrals are established using the truncation method which was
inspired by J. Arthur [2]. We refer to [22] for a nice survey of
his works.

\vskip 0.2cm Though it is well known that Problem 1 is rather
different from Problem 2, there are some analogies between them.
In Problem 1, one has Harish-Chandra's theory of discrete series
and his theory of the generalized principal series, the associated
sphericalized matrix coefficients of $K$-vectors, the so called
Eisenstein integrals. In Problem 2, one has Langlands' general
theory of Eisenstein series and spectral analysis of
$L^2(\G\backslash G)$ applies. Harish-Chandra used the name
Eisenstein integrals to emphasize the analogy with Eisenstein
series. The reason why I deal with the above two unrelated
problems here is that I want to discuss the above two problems for
the Siegel-Jacobi space (see below) at the same time.

\vskip 0.3cm The last half part of this paper will deal with the
Siegel-Jacobi space that is meaningful arithmetically and
geometrically. I will explain this roughly.

\vskip 0.2cm We consider the Jacobi group $G^J=Sp(n,\BR)\ltimes
H_{\BR}^{(n,m)}$ which is the semidirect product of the symplectic
group $Sp(n,\BR)$ of degree $n$ and the Heisenberg group
$H_{\BR}^{(n,m)}$ (see (5.14)). $G^J$ is not a reductive Lie group
and its associated homogeneous space $\Hnm$ is not a homogeneous
space of reductive type, where
\begin{equation*}
\BHn=\left\{ \Omega\in \BC^{(n,n)}\,|\ \Omega=\,{}^t\Omega,\
\textrm{Im}\,\Omega > 0\ \right\}
\end{equation*}
is the Siegel upper half plane of degree $n$ and $\BC^{(m,n)}$ is
the space consisting of $m\times n$ complex matrices. In fact,
$G^J$ acts on the Siegel-Jacobi space $\Hnm$ transitively
\begin{equation*}\left(\begin{pmatrix} A&B\\
C&D\end{pmatrix},(\lambda,\mu;\kappa)\right)\cdot
(\Om,Z)=((A\Om+B)(C\Om+D)^{-1},(Z+\lambda \Om+\mu) (C\Om+D)^{-1}),
\end{equation*}
where $M=\begin{pmatrix} A&B\\
C&D\end{pmatrix} \in Sp(n,\BR),\ (\lambda,\mu; \kappa)\in
H_{\BR}^{(n,m)}$ and $(\Om,Z)\in \BH_n\times \BC^{(m,n)}.$ We
discuss harmonic analysis on $\Hnm$ which generalizes the theory
of harmonic analysis on $\BHn$. In order to develop the theory of
the harmonic analysis of the Siegel-Jacobi space, we need to
review the theory of harmonic analysis on reductive symmetric
spaces.

\vskip 0.2cm This article is organized as follows. In Section 2,
we review the harmonic analysis on $\BR^n$ and compact homogeneous
spaces. In Section 3, we survey some results in harmonic analysis
on semisimple symmetric spaces. In Section 4, we review harmonic
analysis on reductive symmetric spaces. We note that semisimple
symmetric spaces are reductive symmetric spaces. We briefly
describe some recent results on the Plancherel formula for a
reductive symmetric space obtained by van den Ban and
Schlichtkrull [5,\,13,\,14,\,15] and Delorme [19,\,20,\,21,\,22].
In Section 5, we discuss the Siegel-Jacobi space $\Hnm$ and
present some new results (Theorem 5.1, 5.2, 5.3, 5.4, 5.6, 5.7 and
5.8). We also give some open problems to be investigated in the
future.

\vskip 0.2cm  \noindent {\bf Notations\,:} \ We denote by $\BR$
and $\BC$ the field of real numbers, and the field of complex
numbers respectively. We denote by $\BZ$ and $\BZ^+$ be the set of
all integers and the set of all positive integers respectively.
The symbol ``:='' means that the expression on the right is the
definition of that on the left. We denote by $C(X)$\,(resp.
$C^{\infty}(X))$ the space of continuous (resp. smooth) functions
on a manifold $X$. We denote by $C_c^{\infty}(X)$ the space of all
smooth functions on $X$ with compact support. For two positive
integers $k$ and $l$, $F^{(k,l)}$ denotes the set of all $k\times
l$ matrices with entries in a commutative ring $F$. For a square
matrix $A\in F^{(k,k)}$ of degree $k$, $\textrm{tr}(A)$ denotes
the trace of $A$. For an $m\times n$ complex matrix $\Omega,\
\text{Re}\,\Omega$ ({\it resp.}\ $\textrm{Im}\,\Omega)$ denotes
the real ({\it resp.}\ imaginary) part of $\Omega.$ For any $M\in
F^{(k,l)},\ ^t\!M$ denotes the transpose matrix of $M$. For a
matrix $A\in F^{(k,k)}$ and $B\in F^{(k,l)},$ we write
$A[B]=\,^tBAB$. $I_n$ denotes the identity matrix of degree $n$.
For a square matrix $A$ of degree $n$, the inequality $A>0$ means
that $A=\,^t{\bar A}$ is a positive definite hermitian matrix,
i.e., $^t{\bar x}Ax>0$ for any nonzero $x\in\BC^{(n,1)}.$
Obviously $A>0$ means in the case of a real $A$ that $A=\,^tA$ is
a positive definite symmetric matrix, i.e., $^txAx>0$ for any
nonzero $x\in\BR^{(n,1)}.$ For two $n\times n$ square matrices $A$
and $B$, the inequality $A>B$ means that $A-B>0$, i.e., $A-B$ is a
positive definite hermitian matrix. In the real case, the
inequality $A>B$ means that $A-B$ is a positive definite symmetric
matrix

\end{section}
%
%
\begin{section}{{\bf Harmonic Analysis on $\BR^n$ and Compact Homogeneous Spaces }}
\setcounter{equation}{0} In this section, we review the harmonic
analysis on the Euclidean space $\BR^n$ and some examples of
compact homogeneous spaces.

\vskip 0.2cm\noindent {\bf 2.1. Harmonic Analysis on} $\BR^n$

\vskip 0.2cm Let $C_c^{\infty}(\BR^n)$ be the space of smooth
functions on $\BR^n$ with compact support. Let $(\,\,,\,)$ be the
standard inner product on $\BC^n$ and let $|x|$ be the norm of
$x=(x_i)\in\BC^n$ defined by $|x|^2=(x,x).$ The restriction of
$(\,\,,\,)$ to $\BR^n$ is the standard inner product on $\BR^n$
and is also denoted by the same notation $(\,\,,\,)$. We denote by
$\textrm{PW}(\BR^n)$ the space of rapidly decreasing entire
functions of exponential type. More precisely, $\varphi$ is an
element of $\textrm{PW}(\BR^n)$ if and only if it extends to an
entire function on $\BC^n$ for which there exists a positive
number $R>0$ such that the following condition holds for all
$N\in\BZ^+$\,:
\begin{equation*}
\sup_{\la\in\BC^n} \left(
1+|\la|\right)^N\,e^{-R\,|\textrm{Im}\,\la|}\,|\varphi(\lambda)|<\infty.
\end{equation*}

If $f$ is an element of $C_c^{\infty}(\BR^n)$, the Fourier
tranform of $f$ is defined by
\begin{equation*}
{\widehat f}(y)=\left( {1\over {2\pi}}\right)^{\frac n2}
\int_{\BR^n}f(x)\,e^{-i(x,y)}dx.
\end{equation*}

Then we have the well-known classical results\,: \vskip 0.2cm
$\bigstar$ {\bf Inversion Formula\,:} If $f\in
C_c^{\infty}(\BR^n)$, then
\begin{equation*}
f(x)=\left( {1\over {2\pi}}\right)^{\frac n2}
\int_{\BR^n}{\widehat f}(y)\,e^{i(y,x)}dy.
\end{equation*}

\vskip 0.2cm $\bigstar$ {\bf Plancherel Theorem\,:} $f\mapsto
{\widehat f}$ extends to an isometry of the Hilbert space
$L^2(\BR^n)$ onto $L^2(\BR^n)$.

\vskip 0.2cm $\bigstar$ {\bf Paley-Wiener Theorem\,:} $f\mapsto
{\widehat f}$ is a bijection of $C_c^{\infty}(\BR^n)$ onto the
Paley-Wiener space $\textrm{PW}(\BR^n)$.

\vskip 0.2cm The above classical Plancherel theorem can be
expressed using the unitary representations of the noncompact
abelian group $\BR^n$. First we view $\BR^n$ as a homogeneous
space of $\BR^n$ acting on itself by translation. The invariant
differential operators on $\BR^n$ are just the differential
operators with constant coefficients, and the joint eigenfunctions
are the constant multiples of the exponential functions. For
brevity, we set $G=\BR^n$. Then the unitary dual of $G$ is
isomorphic to $\BR^n$. For $\la=(\la_1,\cdots,\la_n)\in\BR^n$, we
write
\begin{equation*}
f_{\la}(x):=e^{i\,(\la,x)},\quad x\in\BR^n\ \ \ \textrm{and}\ \ \
V_{\la}=\BC\cdot f_{\la}.
\end{equation*}

A parameter $\la=(\la_i)\in\BR^n$ corresponds to the
one-dimensional representation $\pi_{\la}$ of $G$ on $V_{\la}$
defined by
\begin{equation*}
\pi_{\la}(g)f_{\la}(x):=f_{\la}(x-g),\quad g\in G.
\end{equation*}

\noindent Then we have the following Plancherel formula
\begin{equation*}
L^2(\BR^n)=\int_{\widehat G}\,V_{\la}\,{{d\la}\over
{(2\pi)^{n/2}}}=\int_{\BR^n}\,V_{\la}\,{{d\la}\over {(2\pi)^{n/
2}}},
\end{equation*}

\noindent where $d\la$ is the Lebesque measure on $\BR^n$. In
other words, $L^2(\BR^n)$ is the direct integral of
one-dimensional Hilbert spaces $V_{\la}$ over the  Lebesque
measure on $\BR^n$. The Fourier transform ${\widehat f}$ of $f$
can be considered as a function on $\wG$.

\vskip 0.2cm We can also consider $\BR^n$ as the homogeneous space
$M(n)/O(n)$, where $M(n)=O(n)\ltimes \BR^n$ is the motion group of
$\BR^n$ equipped with the following multiplication
\begin{equation*}
(k_1,a_1)\cdot (k_2,a_2)=(k_1k_2,\,k_2^{-1}a_1+a_2),\quad
k_1,k_2\in O(n),\ a_1,a_2\in \BR^n.
\end{equation*}

\noindent We observe that $O(n)$ is the stabilizer of the action
of $M(n)$ at $o=(0,\cdots,0)$. In this case, the Laplacian
\begin{equation*}
L_{\BR^n}={ {\partial^2}\over {\partial x_1^2}}+\cdots+{
{\partial^2}\over {\partial x_n^2}}
\end{equation*}

\noindent generates the algebra of all $M(n)$-invariant
differential operators on $\BR^n$. We may view a function on
$\BR^n$ as an $L^2(S^{n-1})$-valued function on $\BR^+$ by means
of the geodesic polar coordinates $x=r\,\omega\,(r>0,\ \omega\in
S^{n-1})$, where $S^{n-1}$ is the $(n-1)$-dimensional unit sphere
in $\BR^n$. Then we can write the Fourier transform ${\widehat f}$
of a function $f$ on $\BR^n$ in the form
\begin{equation*}
{\widehat f}(r,\omega)=\int_{\BR^n}f(y)\,e^{-i\,(r\,
\omega,\,y)}\,dy.
\end{equation*}

\noindent Let ${\mathcal H}:=L^2_B(\BR^+,\,r^{n-1}dr)$ be the
space of $L^2(S^{n-1})$-valued functions on $\BR^+$ which are
square integrable with respect to the measure $r^{n-1}dr$. Then
${\widehat f}\in {\mathcal H}$. The Fourier transform $f\mapsto
{\widehat f}$ maps $L^2(\BR^n)$ isometrically onto $\CH$. The
decomposition of the left regular representation $\pi_L$ of $M(n)$
on $L^2(\BR^n)$ can be interpreted as follows. First of all, for
each $t\in\BR$, we define the representation $\rho_t$ of $M(n)$ on
$L^2(S^{n-1})$ by
\begin{equation*}
\left(
\rho_t(k,a)\varphi\right)(\omega)=e^{\,i\,(t\,a,\,\omega)}\,\varphi
(k^{-1}\omega),\quad k\in O(n),\ a\in\BR^n,\ \omega\in S^{n-1}.
\end{equation*}

\noindent We can prove that $\rho_t$ is irreducible for $t\neq 0$,
and that $\rho_t\cong \rho_{-t}$. Next we define a unitary
representation $\rho$ of $M(n)$ on $\CH$ by
\begin{equation*}
\left( \rho(g)\phi\right)(t):=\rho_{-t} (g)(\phi(t)),\quad g\in
M(n),\ \phi\in \CH,\ t\in\BR^+.
\end{equation*}

\noindent Then $\rho$ is isomorphic to the direct integral of the
$\rho_{-t}\,(t\in\BR^+).$ Let $ \textbf{1}\in L^2(S^{n-1})$ be the
distinguished vector given by $ \textbf{1}(\omega)=1\ (\omega\in
S^{n-1}).$ For brevity, we set $H:=O(n)$. If $f\in
C_c^{\infty}(\BR^n)$, then
\begin{equation*}
{\widehat
f}(t,\omega)=\int_{M(n)}f(gH)\,\rho_{-t}(g)\textbf{1}\,dg=\rho_{-t}(f)\textbf{1}.
\end{equation*}

\noindent Therefore the map $f\mapsto \wf$ is a $M(n)$-equivariant
map from $C_c^{\infty}(\BR^n)$ into $\CH$. Hence the map $f\mapsto
\wf$ extends to an isometry of $L^2(\BR^n)$ onto $\CH$, and we
have the Plancherel decomposition of $\pi_L$ for the homogeneous
space $M(n)/H=\BR^n$. The inversion formula can be reformulated as
follows\,: for $f\in C_c^{\infty}(\BR^n)$, we get
\begin{equation*}
f(gH)=c\int_{\BR^+} \,\langle\, \wf(t),\,\rho_{-t}(g)\textbf{1}\,
\rangle \,r^{n-1}\,dr,
\end{equation*}

\noindent where $c$ is a constant and $\langle\,\,,\,\rangle$ is
an inner product on $L^2(S^{n-1}).$

\vskip 0.5cm\noindent {\bf 2.2. Harmonic Analysis on Compact
Homogeneous Spaces}

\vskip 0.2cm Let $G$ be a compact Lie group. It is known that the
unitary dual of $G$ is discrete. For $\tau\in \wG,$ we let
$\tau^*$ be the contragredient (or dual) representation of $G$. As
mentioned before, $G$ may be regarded as a homogeneous space
$(G\times G)/G$. Indeed the group $G\times G$ acts on $G$
transitively by
\begin{equation}
(g,h)\cdot x=gxh^{-1},\quad g,h,x\in G.
\end{equation}

\noindent The stabilizer $H$ at the identity element $e$ is given
by
\begin{equation*}
H=\left\{\ (g,g)\in G\times G\,|\ g\in G\,\right\}.
\end{equation*}

\noindent Thus $G$ is diffeomorphic to $(G\times G)/G.$ Peter and
Weyl [40] proved that the Hilbert space $L^2(G)$ is decomposed
discretely into a $(G\times G)$-invariant orthogonal direct sum
\begin{equation}
L^2(G)\cong L^2((G\times G)/G)\cong
{\widehat{\oplus}}_{\tau\in\wG}\,\,\tau\otimes \tau^*
\end{equation}

\noindent and each $\tau\in \wG$ is finite dimensional. This gives
the decomposition of the regular representation of $G\times G$ on
$L^2(G)$ into irreducible subrepresentations of $G\times G$. The
formula (2.2) is understood through the matrix coefficient map
$\Phi_{\tau}:V_{\tau}\otimes V_{\tau}^*\lrt L^2(G)$ defined by
\begin{equation}
\Phi_{\tau}(v\otimes
v^*)(g):=(\dim\,V_{\tau})^{1/2}\,\langle\,v,\,\tau^*(g)v^*\,\rangle,\quad
g\in G,\ v\in V_{\tau},\ v^*\in V^*_{\tau},
\end{equation}

\noindent where $V_{\tau}$ is the representation of $\tau\in \wG$
and $V_{\tau}^*$ is its dual space. It is easily seen that the
matrix coefficient map $\Phi_{\tau}\,(\tau\in \wG)$ is a $G\times
G$-homomorphism of $V_{\tau}\otimes V_{\tau}^*$ into $L^2(G)$ with
the left times right action. According to the Schur orthogonality
relations, $\Phi_{\tau}$ is an injective map preserving the inner
products. We identify $V_{\tau}\otimes V_{\tau}^*$ with the image
$\Phi_{\tau}(V_{\tau}\otimes V_{\tau}^*)$.

\vskip 0.2cm \noindent {\bf Example 2.1.} Let
$T=S^1\times\cdots\times S^1$ be the $n$-dimensional torus, where
$S^1$ is the unit circle in $\BC$. Then the unitary dual
${\widehat T}$ is isomorphic to $\BZ^n$. For each
$\al=(\al_1,\cdots,\al_n)\in\BZ^n$, we define $h_{\al}:T\lrt \BC$
by
\begin{equation*}
h_{\al}(t)=e^{i\,(\al_1 t_1+\cdots+\al_n t_n)},\quad
t=(t_1,\cdots,t_n)\in T.
\end{equation*}

\noindent We set $W_{\al}:=\BC\cdot h_{\al}.$ A parameter
$\al=(\al_i)\in\BZ^n$ corresponds to the one-dimensional
representation $\tau_{\al}$ of $T$ on $W_{\al}$ defined by
\begin{equation*}
\tau_{\al}(g)h_{\al}(t):=h_{\al}(t-g),\quad g,t\in T.
\end{equation*}

\noindent If $\psi\in L^2(T),$ according to Peter-Weyl Theorem,
$\psi$ can be written as the direct sum
\begin{equation}
\psi (t)=\left({1\over {2\pi}}\right)^{\frac
n2}\,\sum_{\al\in\BZ^n}\,{\widehat{\psi}}(\al)\,h_{\al}(t),
\end{equation}

\noindent where
\begin{equation*}
{\widehat \psi}(\al)=\left({1\over {2\pi}}\right)^{\frac n2}
\int_0^{2\pi}\cdots\int_0^{2\pi} \psi(t)\,e^{-i\,(\al_1
t_1+\cdots+\al_n t_n)}\,dt_1\cdots dt_n.
\end{equation*}

\noindent The formula (2.4) is nothing but the Fourier series
expansion of a periodic function of period $2\pi$.

\vskip 0.2cm \noindent {\bf Example 2.2.} Let
\begin{equation*}
G=SU(2)=\left\{ \begin{pmatrix} p&q\\
-{\overline q}&{\overline p}\end{pmatrix}\in \BC^{(2,2)}\,\Big|\
|p|^2+|q|^2=1\,\right\}
\end{equation*}

\noindent be the special unitary group of degree $2$. For a
nonnegative integer $n\geq 0$, we denote by $V_n$ the vector space
of complex homogeneous polynomials of degree $n$ in two complex
variables $z_1$ and $z_2.$ That is, $V_n=\sum_{k=0}^n\BC\cdot
z_1^kz_2^{n-k}.$ Let $\Phi_n$ be the representation of $G$ on
$V_n$ defined by
\begin{equation*}
\Phi_n(g)P \begin{pmatrix} z_1\\
z_2\end{pmatrix}:=P\left(g^{-1}\begin{pmatrix} z_1\\
z_2\end{pmatrix}\right),\quad g\in G,\ P\in V_n.
\end{equation*}

\noindent We see that $\Phi_n$ is an irreducible representation of
$G$ and that the unitary dual of $G$ is given by
\begin{equation*}
\wG=\left\{ \Phi_n\,|\ n\in \BZ,\ n\geq 0\ \right\}.
\end{equation*}

\noindent We recall that for $F\in L^2(G),$
\begin{equation*}
\Phi_n(F)=\int_G F(g)\Phi_n(g)\,dg.
\end{equation*}

\noindent Let $u_0,u_1,\cdots,u_n$ be an orthonormal basis of
$V_n$. The Hilbert-Schmidt norm $||\Phi_n(F)||_{\textrm{HS}}$ of
$\Phi_n(F)$ is defined by
\begin{equation*}
||\Phi_n(F)||^2_{\textrm{HS}}:=\sum_{i,j=0}^n|(\Phi_n(F)u_i,u_j)|^2.
\end{equation*}

\noindent We observe that it is well defined, i.e., is independent
of the choice of an orthonormal basis of $V_n$. According to
Peter-Weyl theorem, we get the following Plancherel formula
\begin{equation*}
\int_G
|F(x)|^2\,dx=\sum_{k=0}^{\infty}\,(k+1)\,||\Phi_k(F)||^2_{\textrm{HS}},\quad
F\in L^2(G).
\end{equation*}

\noindent Let $\Theta_n$ be the character of $\Phi_n$. Then if
$f\in C^{\infty}(G),$ we get the Fourier inversion formula
\begin{equation*}
f(I_2)=\sum_{n=0}^{\infty}\,(n+1)\,\Theta_n(f).
\end{equation*}

\vskip 0.2cm \noindent {\bf Example 2.3.} Let $H$ be a closed
subgroup of a compact Lie group $G$. Then the homogeneous space
$G/H$ inherits an invariant measure from the Haar measure on $G$.
Since $L^2(G/H)$ may be identified with the space of right
$H$-invariant functions in $L^2(G)$, according to the Peter-Weyl
theorem, we have the following orthogonal direct sum decomposition
\begin{equation}
L^2(G/H)\cong {\widehat\oplus}_{\tau\in\wG}\,\,V_{\tau}\otimes
\left( V_{\tau}^*\right)^H,
\end{equation}

\noindent where $\left(V_{\tau}^*\right)^H$ is the space of
$\tau^*(H)$-fixed vectors in $V_{\tau}^*$. The decomposition (2.5)
is orthogonal and equivariant with respect to the $G$-action ($G$
acts on the tensor product by its action on the first factors),
and hence the formula (2.5) gives the decomposition of the regular
representation $\pi_H$ for $G/H.$ We note that the decomposition
(2.5) contains only the representation $\tau\in \wG$ such that
$\left( V_{\tau}^*\right)^H\neq 0$, or equivalently
$V_{\tau}^H\neq 0.$

\vskip 0.2cm \noindent {\bf Example 2.4.} Let $S^n$ be the
$n$-dimensional sphere in $\BR^{n+1}$. Then $S^n$ may be regarded
as the homogeneous space $O(n+1)/O(n)\cong SO(n+1)/SO(n).$ If we
identify $\BR^{n+1}$ with $\BR\times S^n$ using the geodesic polar
coordinates, the Laplacian $L_n$ on $S^n$ is characterized as
\begin{equation*}
L_{\BR^{n+1}}={ {\partial^2}\over {\partial r^2}}+ {\frac nr}{
{\partial}\over {\partial r}}+{ 1\over {r^2}}\,L_n,
\end{equation*}

\noindent where $L_{\BR^{n+1}}$ is the Laplacian on
$\BR^{n+1}$\,(cf.\,[29]). If $n=2$, using the spherical polar
coordinates
\begin{equation*}
x_1=r\,\cos \psi\,\sin\theta,\quad x_2=r\,\sin
\psi\,\sin\theta,\quad x_3=r\,\cos\theta
\end{equation*}

\noindent on $\BR^3$, then
\begin{equation*}
L_{\BR^{3}}={ {\partial^2}\over {\partial r^2}}+ {\frac 2r}{
{\partial}\over {\partial r}}+{ 1\over {r^2}}\,L_2,
\end{equation*}

\noindent where
\begin{equation*}
L_2={ {\partial^2}\over {\partial \theta^2}}+ \cot\theta\, {
{\partial}\over {\partial \theta}}+(\sin\theta)^{-2}\, {
{\partial^2}\over {\partial\psi^2}}.
\end{equation*}

It is known that $L_n$ generates the algebra of $O(n+1)$-invariant
differential operators on $S^n$\,(cf.\,[29],\,p.\,17). Let
$s=(s_1,\cdots,s_{n+1})$ denote the Cartesian coordinates of
points on $S^n$. The eigenspaces of the Laplacian $L_n$ on $S^n$
are of the form
\begin{equation*}
E_k= \textrm{span\ of}\
\left\{\,f_{a,k}(s)=(a_1s_1+\cdots+a_{n+1}s_{n+1})^k,\ s\in
S^n\,\right\},
\end{equation*}

\noindent where $k\in\BZ^+\cup\left\{ 0\right\}$ and
$a=(a_1,\cdots,a_{n+1})$ is an isotropic vector in $\BC^{n+1}$,
i.e., $a_1^2+\cdots+a_{n+1}^2=0.$ The eigenvalue of a function in
$E_k$ is given by $-k(n-1+k).$ We observe that $E_k$ is invariant
under the action of $O(n+1)$. We can show that the eigenspace
representation $T_k\,(k\in\BZ^+\cup\left\{ 0\right\})$ of $O(n+1)$
on $E_k$ defined by
\begin{equation*}
\left( T_k(g)f\right)(s)=f(g^{-1}\cdot s),\quad g\in O(n+1), \
s\in S^n,\ f\in E_k
\end{equation*}

\noindent is irreducible.

\vskip 0.05cm Let $P_k$ be the space of homogeneous polynomials of
degree $k$ on $\BR^{n+1}$ and let $H_k\subset P_k$ the subspace of
harmonic polynomials.
It is easily checked that
\begin{equation*}
P_k=q\,P_{k-2}+H_k,\quad q(x):=x_1^2+\cdots+x_{n+1}^2.
\end{equation*}
We can see that the restriction map $p\mapsto p|_{S^n}$ of $H_k$
onto $E_k$ is one-to-one. The eigenfunctions of $L_n$ on $S^n$ are
called $\textit{spherical harmonic}$. The name is derived from the
fact that they coincide with the set of restrictions $\cup_{k\geq
0}\left( H_k|_{S^n}\right).$ \vskip 0.05cm The Hilbert space
$L^2(S^n)$ is decomposed into
\begin{equation}
L^2(S^n)=\sum_{k=0}^{\infty}E_k\quad ( \textrm{orthogonal sum }\
\textrm{for}\ \langle\,\,\,,\,\,\rangle\,),
\end{equation}

\noindent where $\langle\,\,\,,\,\,\rangle$ is the inner product
on $L^2(S^{n-1})$ corresponding to a certain normalized measure
(cf.\,[29],\,p.\,18). The regular representation $\rho$ of
$O(n+1)$ on $S^n$ is mutiplicity-free, that is, the decomposition
(2.6) is a multiplicity-free decomposition of $\rho$. Let $d(k)$
be the dimension of $E_k$. Let $\left\{ S_{k,m}\,|\ 1\leq m\leq
d(k)\,\right\}$ be an orthonormal basis of $E_k$. For $f\in
C^{\infty}(S^n)$, the series
\begin{equation*}
\sum_{k,m}a_{k,m}S_{k,m},\quad
a_{k,m}:=\langle\,f,S_{k,m}\,\rangle
\end{equation*}

\noindent converges absolutely and uniformly to $f$. The mapping
$f\mapsto \left\{ a_{k,m}\right\}$ sends $C^{\infty}(S^n)$ onto
the set of all sequences $\left\{ a_{k,m}\,|\ k\in\BZ^+ \cup
\left\{ 0\right\},\ 1\leq m\leq d(k)\,\right\}$ satisfying the
condition
\begin{equation*}
\sup_{k,m}|a_{k,m}|\,k^q < \infty
\end{equation*}

\noindent for each $q\in k\in\BZ^+ \cup \left\{ 0\right\}.$ For
more detail, we refer to [29].

\end{section}
%
%
\begin{section}{{\bf Semisimple Symmetric Spaces}}
\setcounter{equation}{0}

\vskip 0.2cm Let $G$ be a connected semisimple Lie group with
finite center. Let $\theta$ be a Cartan involution of $G$, that
is, an involution whose fixed subgroup
$K:=G^{\theta}=\left\{\,g\in G\,|\ \theta(g)=g\ \right\}$ is a
maximal compact subgroup of $G$. Let $\s$ be an involution of $G$
commuting with $\theta$. Let $H$ be an open subgroup of the group
$G^{\s}$ of fixed points for $\s$. Then the homogeneous space
$X=G/H$ is called a $\textsf{semisimple symmetric space}$.

\vskip 0.2cm Let $\fg$ be the Lie algebra of $G$, and let $\fk$
and $\fh$ be the subalgebras of $\fg$ corresponding to $K$ and $H$
respectively. The differential of $\theta$ (resp. $\s$) is also
denoted by the same symbol $\theta$ (resp. $\s$). Then we have the
following decomposition
\begin{equation}
\fg=\fk\oplus \fp=\fh\oplus \fq
\end{equation}

\noindent of $\fg$ into the $\pm 1$ eigenspaces for $\theta$ and
$\s$ respectively. We note that $\fp$ and $\fq$ are the orthogonal
complements of $\fk$ and $\fh$ with respect to the Killing form
respectively. Since $\theta$ and $\s$ commute with each other, we
get the joint decomposition
\begin{equation}
\fg=\fk\cap \fh \oplus \fh\cap \fp \oplus \fk\cap \fq \oplus
\fp\cap\fq.
\end{equation}

\vskip 0.1cm A $\textit{Cartan subspace}\ \fb$ for $G/H$ is
defined to a maximal abelian subspace of $\fq$ that consists of
semisimple elements. All Cartan subspaces have the same dimension,
which is called the $ \textsf{rank}$ of $G/H$. We say that a
Cartan subspace $\fb$ is $ \textit{fundamental}$ if $\fk\cap\fb$
is maximal abelian in $\fk\cap\fq$, and that it is $
\textit{maximal split}$ if $\fb\cap \fp$ is maximal abelian in
$\fq\cap\fp$. There exist, up to conjugation by $K\cap H,$ a
unique fundamental Cartan subspace and a unique maximal split
Cartan subspace for $G/H$. If the fundamental Cartan subspace for
$G/H$ is contained in $\fk$, it is called a $ \textit{compact
Cartan subspace}$. The dimension of the $\fp$-part of a maximal
split Cartan subspace for $G/H$ is called the $ \textsf{split
rank}$ of $G/H$.

\vskip 0.2cm Here are examples of semisimple symmetric spaces\,:
\begin{eqnarray*}
  SO(n+1)/SO(n)&=& \textrm{the}\ n\!-\!
\textrm{dimensional
sphere},\\
SL(n,\BR)/SO(m,n-m)&=& \textmd{the space of quadratic forms of
signature}\ (m,n-m)\ \textmd{in}\ \BR^n,\\
& & \ \textmd{where}\ 0\leq m\leq n,\\
Sp(n,\BR)/U(n) &=& \textmd{the Siegel upper half plane of degree}\
n,\\
SO_0(p,q+1)/SO_0(p,q)&=& \textmd{the hyperboloid}\ \BH^{p,q}\
\textmd{in}\ \BR^{p+q+1}.
\end{eqnarray*}

\noindent Here
\begin{equation*}
\BH^{p,q}=\left\{ x=(x_k)\in \BR^{p+q+1}\,\Big|\
\sum_{i=1}^px_i^2-\sum_{j=p+1}^{p+q+1}x_j^2=-1\,\right\}
\end{equation*}
\noindent and if $q=0$, we should take a connected component of
$\BH^{p,0}.$ \vskip 0.2cm Another important example is the group
case mentioned before. Let $G$ be a connected semisimple Lie group
with finite center. Then $G_*:=G\times G$ acts on $G$ transitively
via (2.1). Since the stabilizer $d(G_*)$ at the identity $e$ is
given by
\begin{equation*}
d(G_*)=\left\{\,(g,g)\in G_*\,|\ g\in G\ \right\},
\end{equation*}

\noindent $G=G_*/d(G_*).$ We let $\s$ be the involution of $G$
defined by
\begin{equation*}
\s(g,h)=(h,g),\quad g,h\in G.
\end{equation*}

\noindent Then the fixed group $H=G^{\s}$ for $\s$ is equal to
$d(G_*)$. Hence $G=G_*/H$ is a semisimple symmetric space.

\vskip 0.3cm Let $\pi_H$ be the regular representation of $G$ for
a semisimple symmetric space $G/H$ defined by (1.1). An
irreducible representation $\pi$ is said to be a $\textsf{discrete
series representation}$ of $G/H$ if it can be realized as a
minimal closed invariant subspace of $L^2(G/H)$. The closed span
of all such subspaces is denoted by $L^2_d(G/H)$. We denote by $
\textrm{Disc} (G/H)$ the set of all discrete series
representations for $G/H$.

\vskip 0.3cm
\begin{theorem}
([24],\,[36,\,37,\,38]). Let $G/H$ be a semisimple symmetric
space. Then $L^2_d(G/H)$ is nonzero (equivalently $ \textrm{Disc}
(G/H)$ is nonempty) if and only if $ \textrm{rank}\,(G/H)=
\textrm{rank}\,(K/K\cap H).$
\end{theorem}

\vskip 0.2cm  Flensted-Jensen\,[24] proved the fact that rank
condition implies existence of discrete series representations.
This theorem generalizes the analogous result of Harish-Chandra
for the case of the group. \vskip 0.3cm \noindent {\bf Corollary
3.2}\, (Harish-Chandra\,[25]). Let $G$ be a connected semisimple
Lie group. Then the following are equivalent\,: \vskip 0.21cm (1)
$L^2_d(G)\neq 0$, equivalently $ \textrm{Disc} (G)\neq \emptyset$.

\vskip 0.15cm (2) $ \textrm{rank}\,(G)= \textrm{rank}\,(K),$ where
$K$ is a maximal compact subgroup of $G$.

\vskip 0.15cm (3) $G$ has a {\it compact} Cartan subgroup.

\vskip 0.25cm \noindent {\bf Examples 3.1.} (a) If
$G=Sp(n,\BR)\,(n\geq 1),\ SU(p,q)$ with $pq$ even, then $
\textrm{Disc} (G)\neq \emptyset$. \vskip 0.15cm\noindent (b) If
$G=SO(p,q)$ with $pq$ even, then $L^2_d(G)\neq 0$.

\vskip 0.15cm\noindent (c) If $G=SL(n,\BR)\,(n>2), \ SO(p,q)$ with
$pq$ odd, and if $G$ is any complex Lie group, \\ then $
\textrm{Disc} (G)= \emptyset$.

\vskip 0.25cm \noindent {\bf Remark 3.1.} H. Schlichtkrull [41]
and T. Kobayashi [31] generalized the above theorem by
constructing discrete series representations for vector bundles
over symmetric spaces.

\vskip 0.3cm The Plancherel formula for $G/H$ is understood as a
formula decomposing the regular representation $\pi_H$ into
irreducible unitary representations of $G$. In the case of
semisimple symmetric spaces, the Plancherel formula for $G/H$ was
obtained by Peter-Weyl (the compact group case), Gelfand-Naimark
(classical complex semisimple Lie groups) and Harish-Chandra (real
semisimple Lie groups, and Riemannian symmetric space\,; $
\textrm{Disc} (G/H)= \emptyset$\,). We may say that the
determination of $ \textrm{Disc} (G/H)$ plays an important role in
obtaining the Plancherel formula.

\vskip 0.25cm \noindent {\bf Example 3.2.} We shall describe
briefly the harmonic analysis on the unit disk $\BD$. The unit
disk
\begin{equation*}
\BD=\left\{ w\in\BC\,| \ |w|<1\,\right\}
\end{equation*}

\noindent is realized as the homogeneous space $SU(1,1)/K$, where
\begin{equation*}
SU(1,1)=\left\{ \begin{pmatrix} p & q \\ {\overline q} &
{\overline p} \end{pmatrix}\,\Big|\ |p|^2-|q|^2=1,\ p,q\in\BC\
\right\}
\end{equation*}
\noindent and
\begin{equation*}
K=\left\{ \begin{pmatrix} p & 0 \\ 0 & {\overline p}
\end{pmatrix}\,\Big|\ |p|=1,\ p\in\BC\ \right\}.
\end{equation*}

\noindent We note that $SU(1,1)$ acts on $\BD$ transitively by
\begin{equation}
 \begin{pmatrix} p & q \\ {\overline q} &
{\overline p} \end{pmatrix}\cdot w=(pw+q)({\overline
q}w+{\overline p})^{-1},\quad w\in\BD
\end{equation}

\noindent and that $K$ is the stabilizer of the above action at
the origin $o$. Thus $\BD$ is a Riemannian symmetric space of
noncompact type. The metric
\begin{equation}
ds^2={{dw\,d{\overline w}}\over {(1-|w|^2)^2} }={ {du^2+dv^2}\over
{(1-u^2-v^2)^2} },\quad w=u+iv\in\BD
\end{equation}

\noindent is a Riemannian metric invariant under the action (3.3)
of $SU(1,1)$ and its Laplacian $\Delta$ is given by
\begin{equation}
\Delta=4\,(1-|w|^2)^2\,{ {\partial^2}\over{\partial
w\,\partial{\overline w}} }=(1-u^2-v^2)^2\,\left( {
{\partial^2}\over {\partial u^2}}+{ {\partial^2}\over {\partial
v^2}}\right).
\end{equation}

\noindent The differential form
\begin{equation}
dV=(1-u^2-v^2)^{-2}du\wedge dv
\end{equation}

\noindent is a $SU(1,1)$-invariant volume element on $\BD$. Let
$B:=\left\{ w\in\BC\,|\ |w|=1\,\right\}$ be the boundary of $\BD$.
Then the geodesics in $\BD$ are the circular arcs perpendicular to
$B$. The distance between two points $w_1$ and $w_2$ in $\BD$ with
the respect to the metric (3.4) is given by
\begin{equation}
d(w_1,w_2)={\frac 12}\,\log { {|1-{\overline
w_1}w_2|+|w_2-w_1|}\over {|1-{\overline w_1}w_2|-|w_2-w_1|} }.
\end{equation}

\noindent A circle in $\BD$ tangential to $B$ is called a
$\textit{horocycle}$. We see that a horocycle $\xi$ in $\BD$ is
perpendicular to all geodesics in $\BD$ tending to the point of
contact of $\xi$. For a point $w$ on a horocycle $\xi$, with the
point $b$ of contact, we define
\begin{equation}
<w,b>:= \textmd{distance from}\ o\ \textmd{to}\ \xi,
\end{equation}

\noindent where the negative sign is taken if $o$ lies inside
$\xi.$ According to (3.6) and the cosine relation on the triangle,
we see that
\begin{equation}
e^{2<w,b>}= { {1-|w|^2}\over {|w-b|^2} },\quad w\in \xi.
\end{equation}

\noindent Let ${\mathcal S}(B)$ be the space of analytic functions
on $B$ and let ${\mathcal S}'(B)$ its dual space. Since the
elements of ${\mathcal S}'(B)$ generalize measures, it is
convenient to write
\begin{equation*}
T(f)=\int_B f(b)\,dT(b),\quad f\in{\mathcal S}(B),\ T\in {\mathcal
S}'(B).
\end{equation*}

For $\la\in\BC,$ we let
\begin{equation*}
E_{\la}(\BD):=\left\{ f\in C^{\infty}(\BD)\,|\ \Delta
f=-(\la^2+1)f\ \right\}
\end{equation*}

\noindent be the eigenspace of $\Delta$ with the topology induced
from that of $C^{\infty}(\BD)$. It is known (cf.\,[29]) that the
eigenspace representation $T_{\la}$ of $SU(1,1)$ on $E_{\la}(\BD)$
is irreducible if and only if $i\la+1$ is not an even integer. It
can be shown (cf.\,[29]) that the functions
\begin{equation*}
f_{\la,T}(w):=\int_B e^{(i\la+1)<w,b>}dT(b),\quad \la\in\BC,\ T\in
{\mathcal S}'(B)
\end{equation*}

\noindent exhaust all the eigenfunctions of the Lapalcian
$\Delta$, and moreover that if $i\la\neq -1,-3,-5,\cdots,$ then
the mapping $T\mapsto f_{T,\la}$ from ${\mathcal S}'(B)$ to
$E_{\la}(\BD)$ is bijective.

\vskip 0.1cm A $\textit{spherical function}$ on $\BD$ is by
definition a radial eigenfunction of $\Delta$. If $w\in\BD$ and
$r=d(0,w)$, we see from (3.7) that
\begin{equation*}
w=|w|e^{i\theta}=\tanh r\,e^{i\theta}.
\end{equation*}

\noindent In the geodesic polar coordinates $(r,\theta),\ \Delta$
is given by
\begin{equation*}
\Delta={ {\partial^2}\over {\partial r^2}}+ 2\,\coth 2r\,{
{\partial}\over {\partial r}}+4\,\sinh^{-2}(2r)\,{
{\partial^2}\over {\partial \theta^2}}.
\end{equation*}

\noindent Thus a spherical function $\phi$ satisfies the
differential equation
\begin{equation}
{ {\partial^2 \phi}\over {\partial r^2}}+ 2\,\coth 2r\,{
{\partial\phi}\over {\partial r}}=-(\la^2+1)\phi.
\end{equation}

\noindent The following function $\phi_{\la}$ defined by
\begin{equation}
\phi_{\la}(w):=\int_B e^{(i\la+1)<w,b>}db,\quad \la\in\BC
\end{equation}

\noindent satisfies the equation (3.10). By expanding a power
series in $\sinh 2r$, we see that all smooth solutions of (3.10)
are the constant multiples of $\phi_{\la}(w).$

\vskip 0.2cm We put
\begin{equation*}
a_t:= \begin{pmatrix} \cosh t & \sinh t \\ \sinh t & \cosh t
\end{pmatrix}\in SU(1,1),\quad t\in \BR.
\end{equation*}

\noindent It is a well known result of Harish-Chandra
(cf.\,[29],\,p.\,39) that if $ \textmd{Re}(i\la)>0$, the limit
\begin{equation*}
\textbf{c}(\la):=
\lim_{t\lrt\infty}e^{(-i\la+1)t}\,\phi_{\la}(a_t\cdot o)
\end{equation*}

\noindent exists and
\begin{equation}
\textbf{c}(\la)=\pi^{-1/2}\,{ {\G ({\frac 12}i\la)}\over {
\G({\frac 12}(i\la+1))} },\end{equation}

\noindent where $\G (x)$ is the Gamma function. The formula (3.12)
defines a meromorphic function in $\BC$, which is called the
Harish-Chandra's $\textbf{c}$-function.

\vskip 0.2cm $ \textsf{Definition.}$ If $f$ is a radial function
on $\BD$, its $ \textit{spherical transform}\ {\tilde f}$ is
defined by
\begin{equation}
{\tilde f}(\la):=\int_{\BD}f(w)\,\phi_{-\la}(w)\,dV
\end{equation}

\noindent whenever this integral converges. Here $dV$ denotes a
$SU(1,1)$-invariant volume element given by (3.6).

\vskip 0.2cm It can be proved (cf.\,[29],\,p.\,40) that if $f$ is
a radial function in $C_c^{\infty}(\BD)$,
\begin{equation}
f(w)={ 1\over {2\pi^2}}\int_{\BR}{\tilde f}(\la)\,\phi_{\la}(w)\,|
\textbf{c}(\la)|^{-2} \,d\la
\end{equation}

\noindent and
\begin{equation}
\int_{\BD} |f(w)|^2 \,dV={ 1\over {2\pi^2}}\int_{\BR}|{\tilde
f}(\la)|^2\,| \textbf{c}(\la)|^{-2} \,d\la.
\end{equation}

\vskip 0.2cm $ \textsf{Definition.}$ If $f$ is a complex-valued
function on $\BD$, its $ \textit{Fourier transform}\ {\tilde f}$
is defined by
\begin{equation}
{\tilde f}(\la,b):=\int_{\BD}f(w)\,e^{(-i\la+1)<z,b>}\,dV
\end{equation}

\noindent for all $\la\in\BC,\ b\in B$ for which this integral
exists.

\vskip 0.2cm Then we have the following results\,:

\vskip 0.2cm $\bigstar$ {\bf Inversion Formula\,:} If $f\in
C_c^{\infty}(\BD)$, then
\begin{equation}
f(w)= {1\over {4\pi}}\, \int_{\BR}\int_B{\tilde
f}(\la,b)\,e^{(i\la+1)<w,b>}\,\la\,\tanh \left( { {\pi\la}\over
2}\right)\,d\la\, db ,
\end{equation}

\noindent where $d\la$ is the Lebesque measure on $\BR$ and $db$
is the circular measure on $B$ normalized by $\int_B db=1.$

\vskip 0.2cm $\bigstar$ {\bf Paley-Wiener Theorem\,:} $f\mapsto
{\tilde f}$ is a bijection of $C_c^{\infty}(\BD)$ onto the
Paley-Wiener space $\textrm{PW}(\BD)$, where $\textrm{PW}(\BD)$ is
the space of smooth functions $\psi:\BC\times B\lrt \BC$ which is
holomorphic in $\BC$ such that for each $N\in \BZ^+\cup \left\{
0\right\}$ and $R>0$,
\begin{equation*}
\sup_{\la\in\BC,\,b\in B}e^{-R\,|
\textrm{Im}\la|}\,(1+|\la|)^N\,|\psi(\la,b)|<\infty
\end{equation*}

\noindent and
\begin{equation*}
  \int_B e^{(i\la+1)<w,b>}\,\psi(\la,b)\, db =\int_B e^{(-i\la+1)<w,b>}\,\psi(-\la,b)\,
  db
\end{equation*}

\vskip 0.2cm $\bigstar$ {\bf Plancherel Theorem\,:} $f\mapsto
{\tilde f}$ extends to an isometry of the Hilbert space $L^2(\BD)$
onto
\begin{equation*}
L^2\left(\BR^+\times B,\,(2\pi)^{-1}\la\,\tanh \left({
{\pi\la}\over 2}\right)\,d\la\,db\right), \end{equation*}

\noindent where $\BR^+$ denotes the set of all nonnegative real
numbers.

\vskip 0.2cm We refer to [29] for the proof of the above results.
Since $\BD=SU(1,1)/K$ is a Riemannian symmetric space of
noncompact type, the decomposition of $L^2(\BD)$ contains no
discrete series representation.

\vskip 0.3cm \noindent {\bf Example 3.3.} \ If
$X=SO_0(p,q)/SO_0(p-1,q),\
SU(p,q)/(S(U(1)\times U(p-1,q)),\ \textmd{or}$\\
$Sp(p,q)/(Sp(1)\times Sp(p-1,q))\,\,(p>1,\,q>0),$ then $X$ is a
semisimple symmetric space. The decomposition of $L^2(X)$ contains
both a discrete part (discrete series) and a continuous part (the
principal series). In particular, it is known that the
multiplicities are one, except when $X=SO_0(p,1)/SO_0(p-1,1)$,
where the representations of the principal series have
multiplicity 2 (cf.\,[23]).

\end{section}
%
%
\begin{section}{{\bf Homogeneous Spaces of Reductive Type }}
\setcounter{equation}{0}

\vskip 0.2cm In this section, for simplicity we assume that $G$ is
{\it a real reductive linear group}. By definition, $G$ is a
closed subgroup of the general linear group $GL(N,\BR)$ which is
stable under the standard Cartan involution $\theta_0:G\lrt G$
defined by $\theta_0(g)=\,{}^tg^{-1}\,(g\in G).$ Next we assume
that $H$ is a closed subgroup of $G$ which has at most finitely
many connected components. The coset space $G/H$ is said to be a
$\textsf{homogeneous space of reductive type}$ if $G$ and $H$ are
realized as closed subgroups of $GL(N,\BR)$ such that $GL(N,\BR)
\supset G \supset H$ are stable under the standard Cartan
involution $\theta_0$\,(cf.\,[32]).

\vskip 0.2cm We give some examples of reductive homogeneous spaces
of reductive type.
\begin{eqnarray*}
  GL(n,\BR)/O(n)&=&
\textrm{the cone of} \ n\times n\ \textmd{positive real symmetric matrices} ,\\
GL(2n,\BR)/GL(n,\BC)&=& \textmd{the space of complex structures
on}\ \BR^{2n},\\
U(p,q;{\mathbb F})/U(p-m,q;{\mathbb F}) &=& \textmd{the indefinite
Stiefel manifold,\ where}\ {\mathbb F}=\BR,\,\BC,\,\BH.
\end{eqnarray*}

\noindent Here $\BH$ denotes the quaternionic number field and
$U(p,q;\BH)\cong Sp(p,q).$

\vskip 0.2cm Let $G$ be a real reductive linear group and let
$\theta$ be a Cartan involution of $G$. The differential of
$\theta$ is also denoted by the same symbol $\theta$. We have the
so-called Cartan decomposition of $\fg=\fk+\fp$ of $\fg$, where
$\fk$ and $\fp$ are the $\pm 1$ eigenspaces of $\theta$
respectively. There exists a $G$-invariant nondegenerate symmetric
$\BR$-bilinear form $B$ on $\fg$ such that the bilinear form
$B_{\theta}:\fg\times\fg\lrt \BR$ defined by
\begin{equation*}
B_{\theta}(X,Y):=-B(X,\theta Y),\quad X,Y\in \fg
\end{equation*}

\noindent is symmetric and positive definite. If $G$ is a simple
Lie group, $B$ coincides with the Killing form of $\fg$ up to
positive multiples. For instance, we realize $G$ as a subgroup of
$GL(N,\BR)$ which is stable under $\theta_0$. The differential
$\theta_0$ of $\theta_0$ is given by $\theta_0(X)=-\,{}^tX\ (X\in
\fg{\frak l}(N,\BR))$. We set
\begin{eqnarray*}
  K:&=&O(N)\cap G,\\
  \fk:&=&\fg\cap \left\{
 N\times N\ \textmd{real skew-symmetric matrices}\right\} ,\\
\fp:&=& \fg\cap \left\{\,N\times N \ \textmd{real symmetric matrices}\right\},\\
B(X,Y):&=& \textrm{tr} (XY),\quad X,Y\in\fg,\\
B_{\theta_0}(X,Y):&=& \textrm{tr} (X\,{}^tY),\quad X,Y\in \fg,
\end{eqnarray*}

\noindent where $ \textrm{tr} (X)$ denotes the trace of an
$N\times N$ matrix. It is easily seen that $B_{\theta_0}$ is
positive definite.

\vskip 0.2cm We assume that $\s$ is an involution of $G$ commuting
with a Cartan involution $\theta$. Let $H$ be an open subgroup of
the group $G^{\s}$ of fixed points for $\s$. Namely
$G^{\s}=\left\{ g\in G\,|\ \s(g)=g\right\}$. The coset space $G/H$
is said to be a $ \textsf{reductive symmetric space}$.

\vskip 0.2cm We survey the recent progress in harmonic analysis on
reductive symmetric spaces. First we describe the recent works of
van den Ban and Schlichtkrull (cf.\,[3]-[15]). Let $\fk$ and $\fh$
be the subalgebras of $\fg$ corresponding to $K$ and $H$
respectively. Then we have the decomposition (3.1) and (3.2). We
now start with a maximal abelian subspace $\fa_{\fq}$ of $\fp\cap
\fq$ and extend it to a maximal abelian subspace $\fa$ of $\fp$.
Then it is easily seen that $\fa_{\fq}=\fa\cap \fq.$ Let
$P=M_PA_PN_P$ be a Langlands decomposition of $\s\theta$-stable
parabolic subgroup $P$ of $G$ containing $A_{\fq}:=\exp
\fa_{\fq}$. Let $\fm_P,\,\fa_P$ and $\fn_P$ be the subalgebra of
$\fg$ corresponding to the subgroups $M_P,\,A_P$ and $N_P$
respectively. We denote by ${\widehat {M_P}}$ the unitary dual of
$M_P$ and let $\fa_{P\,\BC}^*:= \textmd{Hom}_\BR (\fa_P,\BC)$ be
the complexified linear dual of $\fa_P$. Let $\rho_P$ be the
element of $\fa_P^*$ defined by
\begin{equation*}
\rho_P (a):={\frac 12}\, \textmd{tr}\left(
\textmd{ad}(a)\big|_{\fn_P}\right),\quad a\in \fa_P.
\end{equation*}

\noindent For $(\xi,\CHX)\in {\widehat {M_P}}$ and $\la\in
\fa_{P\,\BC}^*$, we define the induced representation
\begin{equation*}
\pi_{P,\xi,\la}:= \textmd{Ind}_P^G (\xi\otimes e^{\la}\otimes 1)
\end{equation*}

\newcommand\PXL{\pi_{P,\xi,\la}}
\newcommand\HXL{\CH_{P,\xi,\la}}
\noindent as follows\,: The representation space $\HXL$ of $\PXL$
is the Hilbert space of measurable functions $f:G\lrt \CHX$
satisfying the transformation rule
\begin{equation}
f(manx)=a^{\la+\rho_P}\,\xi(m)\,f(x)
\end{equation}

\noindent and for which the restriction $f|_K$ is square
integrable, where $m\in M_P,\,a\in A_P,\, n\in N_P$ and $x\in G.$
The inner product on $\HXL$ is given by
\begin{equation*}
\langle\,f,g\,\rangle:=\int_K
\langle\,f(k),\,g(k)\,\rangle_{\CHX}\,dk,\quad f,g\in \HXL,
\end{equation*}

\noindent where $dk$ denotes the normalized Haar measure on $K$.
$\PXL$ is realized in $\HXL$ as
\begin{equation*}
\pi_{P,\xi,\la}(g)f(x)= f(xg),\quad g,x\in G,\ f\in \HXL.
\end{equation*}

\newcommand\CWP{{}^P{\mathcal W}}
Let $W_P$ denote the centralizer of $\fa_{P,\fq}:=\fa_P\cap \fq$
in $W=N_K(\fa_{\fq})/Z_K(\fa_{\fq})$. Here $N_K(\fa_{\fq})$ and
$Z_K(\fa_{\fq})$ denote the normalizer of $\fa_{\fq}$ in $K$ and
the centralizer of $\fa_{\fq}$ in $K$ respectively. Let $N_{K\cap
H}(\fa_{\fq})$ be the normalizer of $\fa_{\fq}$ in $K\cap H$. We
denote by $W_{K\cap H}$ the image of $N_{K\cap H}(\fa_{\fq})$
under the projection $N_K(\fa_{\fq})\lrt W.$ Clearly $W_{K\cap H}$
is a subgroup of $W$. We fix a collection ${}^P{\mathcal W}$ of
representatives for $W_P\backslash W/W_{K\cap H}$ contained in
$N_K(\fa_{\fq})$. We denote by $P\backslash G/H$ the collection of
$H$-orbits on $P\ba G$ and by $\left( P\ba G/H\right)_{
\textmd{open}}$ the subset of open orbits. Then we can show that
the set $P\backslash G/H$ is finite and that the union
$\cup_{v\in\CWP}PvH$ is an open dense subset of $P\ba G$. Moreover
we see that the map $v\mapsto PvH$ gives a one-to-one
parametrization of $\left( P\ba G/H\right)_{\textmd{open}}$ by
$\CWP$.

\vskip 0.2cm For $v\in \CWP,$ we set
\begin{equation*}
X_{P,v}:=M_P/(M_P\cap vHv^{-1}).
\end{equation*}

\noindent We see that $X_{P,v}$ is a reductive symmetric space.
Let $\textmd{DS}_{P,v}$ be the set of all discrete series
representations of $M_P$ for $X_{P,v}$. We set
\begin{equation*}
{}^*\fa_{P,\fq}:=\fa_{\fq}\cap \fm_P,\quad {}^*A_{P,\fq}= \exp
({}^*\fa_{P,\fq}).
\end{equation*}

For $(\xi,\CHX)\in {\widehat {M_P}}$ and $v\in \CWP$, we denote by
$\CHX^{-\infty}$ the continuous linear dual of the conjugate
Frech{\'e}t space ${\overline{\CHX^{\infty}} }$, that is, the set
of all generalized vectors of $\CHX$ and by $\left(
\CHX^{-\infty}\right)^{M_P\cap vHv^{-1}}$ the space of $(M_P\cap
vHv^{-1})$-fixed vectors in $\CHX^{-\infty}$. If
\begin{equation*}
\textrm{Hom}_{M_P}(\CHX,L^2(X_{P,v}))\neq 0
\end{equation*}

\noindent for $(\xi,\CHX)\in {\widehat {M_P}}$ and $v\in \CWP$, an
element $T$ in $\textrm{Hom}_{M_P}(\CHX,L^2(X_{P,v}))$ restricts
to a continuous linear $M_P$-equivariant map $T^{\infty}$ from
$\CHX^{\infty}$ to the space $L^2(X_{P,v})^{\infty}$ of smooth
vectors in $L^2(X_{P,v})$. By the local Sobolov inequalities, the
space $L^2(X_{P,v})^{\infty}$ is contained in
$C^{\infty}(X_{P,v}).$ By density of $\CHX^{\infty}$ in $\CHX$, it
follows that restriction to the space of smooth vectors induces an
embedding from $\textrm{Hom}_{M_P}(\CHX,L^2(X_{P,v}))$ onto the
subspace $\textrm{Hom}_{M_P}(\CHX,L^2(X_{P,v}))^{R,\infty}$ of $
\textrm{Hom}_{M_P}(\CHX^{\infty},C^{\infty}(X_{P,v})).$ By the
isomorphism\,(cf.\,[5], Lemma 2.1, p.\,14)
\begin{equation*}
{\overline {\left(\CHX^{-\infty}\right)}}^{M_P\cap
vHv^{-1}}\cong\,\,
\textrm{Hom}_{M_P}(\CHX^{\infty},C^{\infty}(X_{P,v})),
\end{equation*}

\noindent the space
$\textrm{Hom}_{M_P}(\CHX,L^2(X_{P,v}))^{R,\infty}$ corresponds to
the subspace
$$\left(\CHX^{-\infty}\right)_{ \textrm{ds}}^{M_P\cap vHv^{-1}}$$
of ${\overline {\left(\CHX^{-\infty}\right)}}^{M_P\cap vHv^{-1}}$.

For $(\xi,\CHX)\in {\widehat {M_P}}$ and $v\in \CWP$, we define
the finite dimensional Hilbert space $V(P,\xi,v)$ by
\begin{equation*}
V(P,\xi,v):=\begin{cases}  \left( \CHX^{-\infty}\right)_{
\textrm{ds}}^{M_P\cap
vHv^{-1}}\ &  \textmd{if}\ \xi\in \textmd{DS}_{P,v}, \\
\ \ \ \ \ 0 & \textmd{otherwise}.\end{cases}
\end{equation*}

\vskip 0.2cm Let $C^{\infty}(P:\xi:\la)$ be the space of smooth
functions $f:G\lrt \CHX$ satisfying the transformation law (4.1).
Let $C^{-\infty}(P:\xi:-{\overline\la})$ be the space of all
generalized vectors of $C^{\infty}(P:\xi:\la)$. We denote by
$C^{-\infty}(P:\xi:-{\overline\la})^H$ the space of $H$-fixed
generalized vectors for $\PXL.$ Let $C^{\infty}(K:\xi)$ denote the
space of smooth functions $\varphi:K\lrt \CHX$ satisfying the
transformation rule
\begin{equation}
\varphi(mk)=\xi(m)\,\varphi(k),\quad k\in K,\ m\in K\cap P=K\cap
M_P.
\end{equation}

\noindent We denote by $C^{-\infty}(K:\xi)$ the space of all
generalized vectors of $C^{\infty}(K:\xi)$.

\vskip 0.2cm If $(\xi,\CHX)\in {\widehat {M_P}}$, we define
$V(P,\xi)$ to be the formal direct sum
\begin{equation*}
V(P,\xi):=\oplus_{v\in \CWP}V(P,\xi,v).
\end{equation*}

\noindent If $\eta\in V(P,\xi)$, the $\eta_{v}$ denotes its
component in $V(P,\xi,v)$.

\vskip 0.2cm $ \textsf{Definition.}$ Let $\eta\in V(P,\xi)$. For
$\la \in (\fa_{P,\fq})_{\BC}^*$ with $-( \textmd{Re}\,\la+\rho_P)$
strictly $P$-dominant (i.e., $-<\textmd{Re}\,\la+\rho_P,\,\al> >0$
for all $\al\in \Sigma (P)),$ we define the function
\begin{equation*}
j(P:\xi:\la:\eta): G\lrt \CHX^{-\infty}
\end{equation*}

\noindent by
\begin{equation*}
j(P:\xi:\la:\eta)(manvh)=\begin{cases}  a^{\la+\rho_P}\xi(m)\eta_v\ &  \textmd{if}\ manvh\in \cup_{v\in\CWP}PvH, \\
\ \ \ \ \ 0 & \textmd{otherwise}.\end{cases}
\end{equation*}

We see that $j(P:\xi:\la:\eta)$ is an element of
$C^{-\infty}(P:\xi:\la)^H$ and has a meromorphic continuation to
$(\fa_{P,\fq})_{\BC}^*$. If $\la$ is regular, the map
\begin{equation*}
j(P:\xi:\la): V(P,\xi)\lrt C^{-\infty}(P:\xi:\la)^H
\end{equation*}

\noindent defined by
\begin{equation*}
j(P:\xi:\la)(\eta)=j(P:\xi:\la:\eta)
\end{equation*}

\noindent is an injective homomorphism. We define
\begin{equation*}
\textmd{DS}_P:=\left\{ \xi\in {\widehat {M_P}}\,|\ V(P,\xi)\neq
0\,\right\}.
\end{equation*}

Let ${\mathcal P}_{\s}$ be the set of all $\s\theta$-stable
parabolic subgroups of $G$ containing $A_{\fq}$. Two parabolic
subgroups $P,Q\in {\mathcal P}_{\s}$ are said to be $\s$-{\it
associated} if their $\s$-split components $\fa_{P,\fq}$ and
$\fa_{Q,\fq}$ are conjugate under the Weyl group $W$. If $P$ and
$Q$ are $\s$-associated, we denote $P\thicksim Q$. It is easily
seen that $\thicksim$ is an equivalence on ${\mathcal P}_{\s}$.
Let ${\mathbb P}_{\s}$ be a set of representatives in ${\mathcal
P}_{\s}$ for the equivalence classes ${\mathcal P}_{\s}/
\thicksim$.

\newcommand\CPS{ {\mathcal P}_{\sigma}}
\newcommand\APC{ (\fa_{P,\fq})_{\BC}^* }

\vskip 0.2cm $ \textsf{Step I.\ The Fourier transform.}$ \vskip
0.2cm Let $f\in C^{\infty}_c(G/H)$. The $ \textsf{Fourier
transform}\ {\widehat f}$ of $f$ is defined by
\begin{equation}
{\widehat
f}(P:\xi:\la)=\int_{G/H}f(x)\pi_{P,\xi,\la}(x)\,j(P:\xi:\la)\,dx
\end{equation}

\noindent for $P\in \CPS,\ \xi\in \textmd{DS}_P$ and generic
$\la\in i\,\fa_{P,\fq}^*$. We can show that for given fixed $P\in
\CPS$ and $\xi\in \textmd{DS}_P$, the map
\begin{equation*}
f\mapsto {\widehat f}(P:\xi:\la)
\end{equation*}

\noindent intertwines the regular representation $\pi_H$ of $G$
and the representation $1\otimes \PXL$ for generic $\la\in
\fa_{P,\fq}^*.$

\vskip 0.2cm $ \textsf{Step II.\ The Plancherel measure.}$ \vskip
0.2cm Let $P,Q\in \CPS$ such that $A_{P,\fq}=A_{Q,\fq}.$ Then
$M_P=M_Q,\ A_P=A_Q,$ and $ \textmd{DS}_P= \textmd{DS}_Q$. Let
${\overline P}$ and ${\overline Q}$ be the opposite parabolic
subgroups of $P$ and $Q$ respectively. Then ${\overline P}$ and
${\overline Q}$ have the Langlands decompositions
\begin{equation*}
{\overline P}=M_PA_P{{\overline N}_P}\quad \textmd{and}\quad
{\overline Q}=M_QA_Q{{\overline N}_Q}.
\end{equation*}

\noindent We write
\begin{equation*}
\Sigma(Q:P)=\left\{ \al\in \Sigma\,| \ \fg_{\al}\subset {\overline
\fn}_Q\cap \fn_P\ \right\}
\end{equation*}

\noindent for $\la\in \APC$ with $< \textmd{Re} \la,\al>$
sufficiently large for each $\al\in \Sigma(Q:P)$, then the
following integral
\begin{equation}
A(Q:P:\xi:\la)f(x):=\int_{N_Q\cap {\overline N}_P }f(nx)\,dn
\end{equation}

\noindent converges absolutely for $f\in C^{\infty}(P:\xi:\la),$
where $dn$ is a suitably normalized Haar measure of $N_Q\cap
{\overline N}_P$. The above operator
\begin{equation*}
A(Q:P:\xi:\la): C^{\infty}(P:\xi:\la)\lrt C^{\infty}(Q:\xi:\la)
\end{equation*}

\noindent is an intertwining operator between $\PXL$ and
$\pi_{Q,\xi,\la}$, and has an analytic continuation to $\APC$. The
operator $A(Q:P:\xi:\la)$ is called the $\textsf{standard
intertwining operator }.$

\vskip 0.2cm For any $P\in \CPS$, we put
\begin{equation*}
\fa_{P,\fq}^{* \textrm{reg}}:=\left\{ \la\in \fa_{P,\fq}^*\,\,|\
<\la,\al>\neq 0\ \textrm{for all}\ \al\in \Sigma (P)\,\right\}.
\end{equation*}

\noindent According to Harish-Chandra's results, $\PXL$ is an
irreducible unitary representation of $G$ when $\xi\in
\textmd{DS}_P$ and $\la\in i\,\fa_{P,\fq}^{* \textrm{reg}}.$ The
standard intertwining operator $A(Q:P:\xi:\la)$ has an adjoint
\begin{equation*}
A(Q:P:\xi:-{\overline{\la}})^*: C^{-\infty}(Q:\xi:\la)\lrt
C^{-\infty}(P:\xi:\la)
\end{equation*}

\noindent which depends meromorphically on $\la\in
(\fa_{P,\fq})^*_{\BC}$ and is equal to the continuous linear
extension of $A(P:Q:\xi:\la)$ which will be also denoted by the
same notation $A(P:Q:\xi:\la)$. The operator
\begin{equation}
A({\overline P}:P:\xi:-{\overline\la})^*\circ A({\overline
P}:P:\xi:\la)
\end{equation}

\noindent is a $G$-intertwining operator from
$C^{\infty}(P:\xi:\la)$ to $C^{\infty}(P:\xi:\la)$ for generic
$\la\in i\,\fa_{P,\fq}^*$. Since $\pi_{P,\xi,\la}$ is irreducible
for any $\xi\in \textmd{DS}_P$ and $\la\in i \fa_{P,\fq}^{*
\textrm{reg}}$, the operator (4.5) is a scalar multiple of the
identity operator
\begin{equation}
\eta(P:\xi:\la)\,I
\end{equation}

\noindent such that $\eta(P:\eta:\,\cdot\, ):\APC\lrt \BC$ is a
meromorphic function. We now define the measure $d\mu_{P,\xi}$ on
$i\,\fa_{P,\fq}^*$ by
\begin{equation}
d\mu_{P,\xi}(\la):={1 \over {\eta(P:\xi:\la)} }\,d\mu_P(\la),
\end{equation}

\noindent where $d\mu_P$ is a suitable Lebesque measure on
$i\,\fa_{P,\fq}^*$.

\vskip 0.2cm From now on, $W_P^*$\,(resp. $W_P$) denotes the
normalizer (resp. centralizer) of $\fa_{P,\fq}$ in the Weyl group
$W$. We define the group
\begin{equation*}
W(\fa_{P,\fq}):=W_P^*/W_P.
\end{equation*}

van den Ban and Schlichtkrull [13] and Delorme [21] obtained the
beautiful Plancherel formula.

\vskip 0.2cm\noindent $ \textbf{Plancherel Formula A\,[13,\,21]:}$
\ $\textit{Let}\ f\in C_c^{\infty}(G/H).\ \textit{Then}$
\begin{equation}
||f||^2_{L^2(G/H)}=\sum_{P\in {\mathbb
P}_{\s}}[W:W_P^*]\sum_{\xi\in \textmd{DS}_P
}\int_{i\,\fa_{P,\fq}^*} ||{\widehat
f}(P:\xi:\la)||^2\,d\mu_{P,\xi}(\la).
\end{equation}

\newcommand\CL{\mathcal L}
\vskip 0.5cm We now describe the Plancherel formula (4.8) in terms
of representations of $G$. For $P\in {\mathcal P}_{\s}$ and
$\xi\in \textmd{DS}_P$, we define the Hilbert space
\begin{equation*}
\CH (P,\xi):={\overline {V(P,\xi)}}\otimes L^2(K:\xi),
\end{equation*}

\noindent equipped with the tensor product inner product. In
addition, we let
\begin{equation*}
\CL^2(P,\xi):=L^2(i\,\fa_{P,\fq}^*,\CH(P,\xi),\,[W:W_P^*]\,d\mu_{P,\xi})
\end{equation*}

\noindent be the space of square integrable functions from
$i\,\fa_{P,\fq}^*$ to $\CH (P,\xi)$ for the measure
$[W:W_P^*]\,d\mu_{P,\xi}$. We define the representation
$\pi_{P,\xi}$ of $G$ on $\CL^2 (P,\xi)$ by
\begin{equation*}
(\pi_{P,\xi}(g)\varphi)(\la):=[1\otimes \PXL
(g)]\varphi(\nu),\quad g\in G,\ \varphi\in \CL^2(P,\xi),\ \la\in
i\,\fa_{P,\fq}^*.
\end{equation*}

\noindent $\pi_{P,\xi}$ is realized as the following direct
integral
\begin{equation*}
\pi_{P,\xi}\cong \int_{i\,\fa_{P,\fq}^*}1_{{\overline {V(P,\xi)}}
}\otimes \PXL\, [W:W_P^*]\,d\mu_{P,\xi}(\la).
\end{equation*}

\newcommand\BPS{ {\mathcal P}_{\sigma}}
\noindent We define $(\pi_P,\,\CL^2(P))$ as the Hilbert direct sum
of the
unitary representations $(\pi_{P,\xi},\,\CL^2(P,\xi)),$\\
$\xi\in \textmd{DS}_P.$ Finally we define the representation
$(\pi,\CL^2)$ as the unitary direct sum of the representations
$(\pi_P,\,\CL^2(P)),\ P\in {\mathbb P}_{\s}.$ Therefore we get
\begin{equation*}
\pi\cong \oplus_{P\in {\mathbb P}_{\s}}{\hat \oplus}_{\xi\in
\textmd{DS}_P }\pi_{P,\xi}.
\end{equation*}

For $P\in {\mathbb P}_{\s}$, there exists a unitary representation
of $W_P^*$ on $\CL^2$ such that this representation commutes with
$\pi$ and is not the trivial representation. Its representation
space is a closed $G$-invariant subspace of $\CL^2$. The map
$f\mapsto {\widehat f}$ extends to an isometry ${\mathcal F}$ from
$L^2(G/H)$ into $\CL^2$ intertwining $\pi_H$ with $\pi$. The image
of the Fourier transform ${\mathcal F}$ is given by
\begin{equation}
\textmd{Im}({\mathcal F})=(\CL^2)^{W_P^*}.
\end{equation}

\noindent The group $W(\fa_{P,\fq})=W_P^*/W_P$ acts freely, but in
general not transitively on the components of $\fa_{P,\fq}^{*
\textmd{reg}}$. Let $\Omega_P$ be a fundamental domain for the
action $W(\fa_{P,\fq})$ on $i\,\fa_{P,\fq}^{* \textmd{reg}}$
consisting of connected components of $i\,\fa_{P,\fq}^{*
\textmd{reg}}$. Then, for each $P\in {\mathcal P}_{\s}$ and
$\xi\in \textmd{DS}_P$, we denote by $\CL^2_{\Omega_P}(P,\xi)$ the
closed $G$-invariant subspace of functions in $\CL^2(P,\xi)$ that
vanish almost everywhere outside $\Omega_P$. We define the closed
$G$-invariant subspace
$$\CL_0^2:=\oplus_{P\in {\mathbb P}_{\s} }{\hat {\oplus}}_{\xi\in
\textmd{DS}_P} \CL^2_{\Omega_P}(P,\xi).$$

\noindent The orthogonal projection of $\CL^2$ onto $\CL^2_0$ is
denoted by $\varphi\mapsto \varphi_0$.

\indent van den Ban and Schlichtkrull [5,\,14] obtained the
following results.

\vskip 0.2cm \noindent $ \textbf{Plancherel Formula B\,[14]:}$ \ \
$ \textit{The map}\ f\mapsto ({\mathcal F}f)_0\ \textit{is an
isometry from} \ L^2(G/H)\ \textit{onto}$\\ $\CL^2_0\
\textit{intertwining}\ \pi_H \textit{with}\  \pi,\
\textit{establishing the Plancherel formula}$
\begin{equation}
\pi_H\cong \oplus_{P\in {\mathbb P}_{\s} }{\hat {\oplus}}_{\xi\in
\textmd{DS}_P} \int_{\Omega_P} 1_{\overline{V(P,\xi)}}\otimes
\PXL\,[W:W_P]\,d\mu_{P,\xi}(\la).
\end{equation}

\noindent $ \textit{In particular, for all}\ P\in {\mathcal
P}_{\s}\ \textit{and}\ \xi\in \textmd{DS}_P,\ \PXL \
\textit{occurs with multiplicity}\ \dim V(P,\xi) \ \textit{for}$\\
$ \textit{almost all}\ \la\in i\,\fa_{P,\fq}^{*}.$

\vskip 0.2cm Next we provide a normalized version of the
Plancherel formula (4.10). For $P\in {\mathcal P}_{\s}$ and
$\xi\in \textmd{DS}_P$, we define
\begin{equation*}
j^{\star}(P:\xi:\la):=A({\overline P}:P:\xi:\la)^{-1}\circ
j({\overline P}:\xi:\la).
\end{equation*}

\vskip 0.3cm \noindent $ \textbf{Regularity Theorem I\,[14]:}
\textit{The}\ \Hom (V(P,\xi),C^{-\infty}(K:\xi))$-$ \textit{valued
function} \ j^{\star}(P:\xi:\,\cdot\,)$\\ $ \textit{is regular on}
\ i\,\fa_{P,\fq}^{*}.$

\vskip 0.2cm For a minimal $\s$-parabolic subgroup $P$, the above
theorem was proved by van den Ban and Schlichtkrull [9], and for a
general $\s$-parabolic subgroup $P$, the above theorem was proved
by Carmona and Delorme [19]. For $P\in {\mathcal P}_{\s}$ and
$\xi\in \textmd{DS}_P$, we now define the $ \textsf{normalized
Fourier transform}\ {\widehat f}_{\star}$ of a function $f\in
C^{\infty}_c(G/H)$ by
\begin{equation}
{\widehat f}_{\star}(P:\xi:\la):=\int_{G/H}f(x)\PXL
(x)j^{\star}(P:\xi:\la)\,dx.
\end{equation}

\noindent According to Regularity Theorem, the function $\la\lrt
{\widehat f}_{\star}(P:\xi:\la)$ is analytic as a ${\overline
{V(P,\xi)}}\otimes C^{\infty}(K:\xi)$-valued function on $
i\,\fa_{P,\fq}^{*}.$ It is shown that for $P\in {\mathcal
P}_{\s},\ \xi\in \textmd{DS}_P$ and $f\in C^{\infty}_c(G/H)$,
\begin{equation*}
||{\widehat f}(P:\xi:\la)||^2\,d\mu_{P,\xi}(\la)=||{\widehat
f}_{\star}(P:\xi:\la)||^2\,d\mu_P.
\end{equation*}

Let
\begin{equation*}
\CL^2_{\star}(P,\xi):=L^2(i\,\fa_{P,\fq}^{*},\,{\mathcal
H}(P,\xi),[W:W_P^*]\,d\mu_P)
\end{equation*}

\noindent be the space of square integrable functions from
$i\,\fa_{P,\fq}^*$ to $\CH (P,\xi)$ with respect to the measure
$[W:W_P^*]\,d\mu_P$. We define the representation
$(\pi_{\star},\CL^2_{\star})$ as the unitary direct sum of the
representation $\CL_{\star}^2(P,\xi)$ $(P\in {\mathbb P}_{\s},\
\xi\in \textmd{DS}_P)$. In other words,
\begin{equation*}
\pi_{\star}\cong \oplus_{P\in {\mathbb P}_{\s} }{\hat
{\oplus}}_{\xi\in \textmd{DS}_P}\CL_{\star}^2(P,\xi).
\end{equation*}

\noindent We denote the continuous linear extension of $f\mapsto
{\widehat f}_{\star}$ by ${\mathcal F}_{\star}:L^2(G/H)\lrt
\CL^2_{\star}.$ It can be shown that there exists a unitary
representation of $W(\fa_{P,\fq})$ in $\CL_{\star}^2$ such that
this representation commutes with $\pi_{\star}$ and is not the
trivial representation. We fix fundamental domains
$\Omega_P\subset i\,\fa_{P,\fq}^{* \textmd{reg}}$ for the action
of $W(\fa_{P,\fq})$, and define $\CL_{\star,0}^2$ in a similar way
as $\CL_0^2$. Then $\CL_{\star,0}^2$ is a closed $G$-invariant
subspace of $\CL^2_{\star}.$ The orthogonal projection of
$\CL^2_{\star}$ onto $\CL^2_{\star,0}$ is denoted by
$\varphi\mapsto \varphi_0.$ The Plancherel formula (4.10) is
essentially equivalent to the following normalized analogue.

\vskip 0.2cm \noindent $ \textbf{Plancherel Formula C\,(A
normalized version)\,[5,\,14]:}$\ \ $\textit{The map}\ f\mapsto
({\mathcal F}_{\star}f)_0\ \textit{is an}$\\ $ \textit{isometry
from}\ L^2(G/H)\ \textit{onto}\ \CL^2_{\star,0}
\textit{intertwining}\ \pi_H \textit{with}\  \pi_{\star},\
\textit{establishing the Plancherel}$\\ $ \textit{formula}$
\begin{equation}
\pi_H\cong \oplus_{P\in {\mathbb P}_{\s} }{\hat {\oplus}}_{\xi\in
\textmd{DS}_P} \int_{i\,\Omega_P} 1_{\overline{V(P,\xi)}}\otimes
\PXL\,[W:W_P]\,d\mu_P(\la).
\end{equation}

\vskip 0.3cm The Plancherel formula (4.12) is derived from the
Plancherel formula for spherical functions. The Plancherel formula
for spherical functions is derived by an application of a residue
calculus.

\vskip 0.2cm We now describe the $\emph{spherical Plancherel
formula}$. Let $(\delta,V_{\delta})\in {\widehat K}$ and let
$L^2(G/H)_{\de}$ be the subspace of $L^2(G/H)$ consisting of left
$K$-finite functions of type $\delta$. We put
$\tau=\tau_{\de}:=\de^*\otimes 1$ and $V_{\tau}=V_{\de}^*\otimes
V_{\de}$. We denote by $L^2(G/H:\tau)$ the space of
$V_{\tau}$-valued functions $f:G/H\lrt V_{\tau}$ satisfying the
transformation rule
\begin{equation}
f(kx)=\tau(k)f(x),\quad k\in K,\ x\in G/H.
\end{equation}

Elements of $L^2(G/H:\tau)$ are called $\tau$-$ \textsf{spherical
functions}$ or simply $ \textit{spherical functions}$. Similarly
we can define $C^{\infty}_c(G/H:\tau),\ C^{\infty}(G/H:\tau)$ and
$C(G/H:\tau)$. We observe that
\begin{eqnarray*}
L^2(G/H)_{\de}&\cong & \Hom_K(V_\de,L^2(G/H))\otimes V_\de \\
&\cong&
\left(L^2(G/H)\otimes V_{\de}^*\right)^K\otimes V_{\de}\\
&\cong& \left(L^2(G/H)\otimes V_{\tau}\right)^K \\ &\cong&
L^2(G/H:\tau).
\end{eqnarray*}

\noindent The natural isomorphism
\begin{equation}
{\mathbb S}:L^2(G/H)_{\de}\lrt L^2(G/H:\tau_\de)
\end{equation}

\noindent is called the $ \textit{sphericalization}$ of
$L^2(G/H)_{\de}$. Let $P\in {\mathcal P}_\s$ and $\xi\in
\textmd{DS}_P.$ For generic $\la\in (\fa_{P,\fq})^*_{\BC},$ we
define the map
\begin{equation*}
M_{P,\xi,\la}:{\overline {V(P,\xi)}}\otimes C^{\infty}(K:\xi)\lrt
C^{\infty}(G/H)
\end{equation*}

\noindent by
\begin{equation}
M_{P,\xi,\la}(\eta\otimes \varphi)(x):=\langle\,
\varphi,\pi_{P,\xi,-{\overline\la}}(x)\,j^{\star}(P:\xi:-{\overline\la})\eta\,\rangle,
\end{equation}

\noindent where $\eta\otimes \varphi\in {\overline
{V(P,\xi)}}\otimes C^{\infty}(K:\xi)$ and $x\in G/H$. We consider
the matrix coefficient
\begin{equation*}
{\mathbb M}_{\xi,v}:{\overline {V(P,\xi,v)}}\otimes \CHX\lrt
L^2(X_{P,v})
\end{equation*}

\noindent defined by
\begin{equation*}
{\mathbb M}_{\xi,v}(\eta\otimes
w)(m):=\langle\,\xi(m)\eta,w\,\rangle,
\end{equation*}

\noindent where $m\in M_P,\ \eta\otimes w\in {\overline
{V(P,\xi,v)}}\otimes \CHX.$ For $\xi\in \textmd{DS}_P$, we set
\begin{equation*}
L^2(X_{P,v})_\xi :={\mathbb M}_{\xi,v}({\overline
{V(P,\xi,v)}}\otimes \CHX).
\end{equation*}

\noindent We note that $L^2(X_{P,v})_\xi\neq 0$ if and only if
$\xi\in \textmd{DS}_P$.

\vskip 0.2cm We set
\begin{equation*}
K_P:=K\cap M_P\quad \textmd{and}\quad
\tau_P=\tau_{\de,P}:=\tau_\de|_{K_P}.
\end{equation*}

\noindent Then
\begin{equation*}
L^2(X_{P,v}:\tau_P)=\left( L^2(X_{P,v})\otimes
V_{\tau}\right)^{K_P}.
\end{equation*}

\noindent We define
\begin{equation*}
L^2(X_{P,v}:\tau_P)_\xi:=\left( L^2(X_{P,v})_\xi\otimes
V_{\tau}\right)^{K_P}.
\end{equation*}

\noindent We can show that there exists a canonical isometry
\begin{equation}
\al_{\xi,\de}:{\overline {V(P,\xi)}}\otimes C^{\infty}(K:\xi)_\de
\lrt \oplus_{v\in {}^P{\mathcal W}}L^2(X_{P,v}:\tau_P)_\xi,
\end{equation}

\noindent where $\oplus$ denotes the formal direct sum of Hilbert
spaces. Moreover ${\overline {V(P,\xi)}}\otimes
C^{\infty}(K:\xi)_\de$ is finite dimensional.

\vskip 0.2cm $ \textsf{Definition of Eisenstein integrals.}$\ \
Let $\psi\in \oplus_{v\in {{}^P\mathcal
W}}\,L^2(X_{P,v}:\tau_P)_\xi$. For $\la\in (\fa_{P,\fq})^*_\BC$,
we define the $ \textsf{normalized Eisenstein integral}\
E^0(P:\psi:\la)$ by
\begin{equation}
E^0(P:\psi:\la):={\mathbb
S}(M_{P,\xi,-\la}(\al_{\xi,\de}^{-1}(\psi)))\in L^2(G/H:\tau).
\end{equation}

From now on, for any finite dimensional unitary representation
$(\chi,V_{\chi})$ of $K$, we let ${\mathcal A}_2(G/H:\chi)$ be the
space of smooth functions $f\in C^{\infty}(G/H:\chi)$ satisfying
the following conditions \vskip 0.1cm \ $\quad \ \ (a)\quad f\in
L^2(G/H:\chi)$\,;
\vskip 0.1cm \ $\quad\ \ (b)\ \ \BD(G/H)f$ is
finite dimensional.

\vskip 0.1cm \noindent Here $\BD(G/H)$ denotes the commutative
algebra of all $G$-invariant differential operators on $G/H.$
According to [37], ${\mathcal A}_2(G/H:\chi)$ is finite
dimensional and is decomposed into the orthogonal direct sum
\begin{equation}
{\mathcal A}_2(G/H:\chi)=\oplus_{\pi\in \textmd{Disc}(G/H) }
L^2(G/H:\chi)_{\pi}.
\end{equation}

\noindent In fact, the number of $\pi\in \textmd{Disc}(G/H)$ such
that $L^2(G/H:\chi)_{\pi}\neq 0$ is finite. We define the formal
sum of Hilbert spaces
\begin{equation*}
{\mathcal A}_{2,P}:=\oplus_{v\in {}^P{\mathcal W}}{\mathcal A}_2
(X_{P,v}:\tau_P).
\end{equation*}

\noindent By the formula (4.18) for $X_{P,v}\,(v\in {}^P{\mathcal
W}),\ {\mathcal A}_{2,P}$ is decomposed into the orthogonal direct
sum
\begin{equation}
{\mathcal A}_{2,P}:=\oplus_{\xi\in \textmd{DS}_P}  \oplus_{v\in
{}^P{\mathcal W}}{\mathcal A}_2 (X_{P,v}:\tau_P)_\xi.
\end{equation}

For $\psi\in {\mathcal A}_{2,P}$, we define the $
\textsf{normalized Eisenstein integral}\ E^0(P:\psi:\la)$ by
\begin{equation}
E^0(P:\psi:\la):=\sum_{\xi\in \textmd{DS}_P} E^0(P:\psi_\xi:\la),
\end{equation}

\noindent where $\psi_{\xi}$ is the component in ${\mathcal
A}_{2,P,\xi}:=\oplus_{v\in {}^P{\mathcal W}}L^2
(X_{P,v}:\tau_P)_\xi.$

\vskip 0.3cm \noindent $ \textbf{Regularity Theorem II:}$ $
\textit{The meromorphic function}$
\begin{equation}
E^0(P:\psi:\,\cdot\,):(\fa_{P,\fq})^*_\BC\lrt C^{\infty}(G/H:\tau)
\end{equation}

\noindent $ \textit{is regular on}\ i\,\fa_{P,\fq}^*\ \textit{for
any}\ \psi\in {\mathcal A}_{2,P}.$

\vskip 0.3cm This regularity theorem is essentially equivalent to
that for $j^{\star}(P:\xi:\,\cdot\,)$. The proof of the above
regularity theorem is based on an asymptotic behavior of
$E^0(P:\psi:\la)$ together with the Maass-Selberg relations among
Harish-Chandra's $C$-functions.

\vskip 0.2cm
For $\la\in (\fa_{P,\fq})^*_\BC$, we consider the map
\begin{equation*}
E^0(P:\la):G/H\lrt \Hom({\mathcal A}_{2,P},V_\tau)
\end{equation*}

\noindent defined by
\begin{equation*}
E^0(P:\la)(x)\psi=E^0(P:\la:x)\psi:=E^0(P:\psi:\la)(x).
\end{equation*}

We define the $ \textsf{dual Eisenstein integral}\ E^*(P:\la:x)$
of $E^0(P:\la:x)$ by
\begin{equation}
E^*(P:\la:x):=E^0(P:-{\overline \la}:x)^*\in \Hom
(V_\tau,{\mathcal A}_{2,P})
\end{equation}

\noindent for $x\in G/H$ and generic $\la\in (\fa_{P,\fq})^*_\BC$.
For a given spherical function $f\in C^{\infty}_c(G/H:\tau)$, we
define the $ \textsf{spherical Fourier transform}\
 {\mathcal F}_P f: i\fa_{P,\fq}^*\lrt  {\mathcal A}_{2,P}$ of $f$
by
\begin{equation}
{\mathcal F}_P f(\la):=\int_{G/H}E^*(P:\la:x)f(x)\,dx.
\end{equation}

We let ${\mathcal S}(i\,\fa_{P,\fq}^*)$ be the space of Schwartz
functions on $i\,\fa_{P,\fq}^*$ and let ${\mathcal C}(G/H)$ be the
$L^2$-Schwartz space of $G/H$ consisting of smooth functions
$f:G/H \lrt \BC$ such that
\begin{equation*}
(1+\log_H)^n Xf\in L^2(G/H)
\end{equation*}

\noindent for any $X\in U(\fg)$ and any $n\in\BZ^+.$ Here
$\log_H:G\lrt [0,\infty)$ is a function on $G/H$ defined by
$\log_H (kah)=||\log\, a||\ (k\in K,\,a\in A,\, h\in H)$ and
$U(\fg)$ denotes the universal enveloping algebra of $\fg$.
According to [4], we get the following facts\,: \vskip 0.1cm $
\textsf{(SF1)}\ \ \langle\, {\widehat
f}_{\star}(P:\xi:\la),T\,\rangle=\langle\,{\mathcal F}{\mathbb
S}(f)(-\la),\,\al_{\xi,\de}(T)\,\rangle$ for $f\in
C_c^{\infty}(G/H:\tau),\ T\in {\overline {V(P,\xi)}}\otimes
L^2(K:\xi)$ and all $\la\in i\,\fa_{P,\fq}^*$.

\vskip 0.1cm $ \textsf{(SF2)}$ For $P\in {\mathcal P}_\s,\
{\mathcal F}_P$ maps $C_c^{\infty}(G/H:\tau)$ into ${\mathcal
S}(i\,\fa_{P,\fq}^*)\otimes {\mathcal A}_{2,P}.$

\vskip 0.1cm $ \textsf{(SF3)}\ {\mathcal F}_P$ extends to a
continuous operator from ${\mathcal C}(G/H:\tau)$ to ${\mathcal
S}(i\,\fa_{P,\fq}^*)\otimes {\mathcal A}_{2,P}.$

\vskip 0.5cm According to (SF3), we can define the $ \textsf{wave
packet operator}\ {\mathcal T}_P:{\mathcal
S}(i\,\fa_{P,\fq}^*)\otimes {\mathcal A}_{2,P}\lrt
C^{\infty}(G/H:\tau)$ by
\begin{equation}
{\mathcal T}_P
\varphi(x):=\int_{i\,\fa_{P,\fq}^*}E^0(P:\la:x)\varphi(\la)\,d\mu_P(\la).
\end{equation}

\noindent It was proved by van den Ban and Schlichtkrull [8,\,10]
that for $f\in C_c^{\infty}(G/H:\tau)$
\begin{equation}
f=\sum_{P\in {\mathbb P}_\s}[W:W_P^*]\, {\mathcal T}_P ({\mathcal
F}_P f)
\end{equation}

\noindent and that for associated $P,Q\in {\mathcal P}_\s$,
\begin{equation*}
{\mathcal T}_P\circ {\mathcal F}_P={\mathcal T}_Q\circ {\mathcal
F}_Q\quad \textmd{on}\ {\mathcal C}(G/H:\tau).
\end{equation*}

\vskip 0.3cm The following asymptotic behavior of the Eisenstein
integral is obtained.

\vskip 0.2cm
\begin{theorem}
Let $P,Q\in {\mathcal P}_\s$ be associated. Let
\begin{equation*}
W(\fa_{Q,\fq}|\fa_{P,\fq}):=\left\{ s|_{\fa_{P,\fq} }\,|\ s\in W,\
s(\fa_{P,\fq})\subset \fa_{Q,\fq}\ \right\}.
\end{equation*}

\noindent Then there are uniquely determined meromorphic
functions, called the $C$-$ \textbf{functions}$
\begin{equation*}
C^0_{Q|P}(s:\,\cdot\,):(\fa_{P,\fq})^*_\BC \lrt  \Hom ({\mathcal
A}_{2,P},{\mathcal A}_{2,Q}),\quad s\in W(\fa_{Q,\fq}|\fa_{P,\fq})
\end{equation*}

\noindent such that
\begin{equation}
E^0(P:\la:mav)\psi \,\thicksim\, \sum_{s\in
W(\fa_{Q,\fq}|\fa_{P,\fq})}
a^{s\la-\rho_Q}\,[C^0_{Q|P}(s:\la)\psi]_v (m)
\end{equation}

\noindent as $a\lrt \infty$ in $A_{Q,\fq}^+$ for all $\la\in
i\,\fa_{P,\fq}^*,\ v\in {\mathcal W}_Q,\ m\in X_{Q,v}$ and
$\psi\in {\mathcal A}_{2,P}.$ Here $s\la:=\la\circ s^{-1}.$
\end{theorem}

\vskip 0.3cm\noindent $ \textbf{Maass-Selberg Relations\,:}$\ \ $
\textit{Let}\ P,Q\in {\mathcal P}_\s\ \textit{be associated and
let}\ s\in W(\fa_{Q,\fq}|\fa_{P,\fq}).\ \textit{Then}$
\begin{equation}
C^0_{Q|P}(s:-{\overline \la})^*\circ C^0_{Q|P}(s:\la)=I.
\end{equation}

\noindent $ \textit{In particular, if}\ \la\in i\,\fa_{P,\fq}^*,\
\textit{then}\ C^0_{Q|P}(s:\la)\ \textit{is unitary}$.

For the group case, the above result is due to Harish-Chandra
[28], for $P$ minimal, it was obtained by van den Ban [4] and for
general $P$, it was due to Carmona and Delorme [19].

\vskip 0.3cm The $C$-functions $C^0_{Q|P}(s:\la)$ have the
following properties.

\vskip 0.2cm $\textsf{(C1)}$\ \ $C^0_{P|P}(1:\la)=I.$ \vskip
0.15cm $ \textsf{(C2)}$\ \ For each $s\in
W(\fa_{Q,\fq}|\fa_{P,\fq})$,
\begin{equation}
E^0(P:\la:x)=E^0(Q:s\la:x)\circ C^0_{Q|P}(s:\la)
\end{equation}

\ \ \ \ \ for any $x\in G/H$ and $\la\in (\fa_{Q,\fq})^*_\BC$.

\vskip 0.15cm $ \textsf{(C3)}$ \ \ Let $P,Q,R\in {\mathcal P}_\s$
be associated and $s\in W(\fa_{Q,\fq}|\fa_{P,\fq})$. Then
\begin{equation}
C^0_{R|P}(ts:\la)=C^0_{R|Q}(t:s\la)\circ C^0_{Q|P}(s:\la)
\end{equation}

\ \ \ \ \ as a $\Hom ({\mathcal A}_{2,P},{\mathcal
A}_{2,R})$-valued identity of meromorphic functions on
$(\fa_{Q,\fq})^*_\BC$.

\vskip 0.15cm $ \textsf{(C4)}$\ \ Let $P,Q\in {\mathcal P}_\s$ be
associated. Then for every $f\in {\mathcal C}(G/H:\tau)$ and $s\in
W(\fa_{Q,\fq}|\fa_{P,\fq})$,
\begin{equation}
{\mathcal F}_Q f(s\la)=C^0_{Q|P}(s:\la)\,{\mathcal F}_P
f(\la),\quad \la\in i\,\fa_{P,\fq}^*.
\end{equation}

The functional equation (4.28) was originally obtained by
Harish-Chandra [28] for the group case. For $P,Q$ minimal, (4.28)
is due to van den Ban [4] and for general $P,Q$, it is due to
Carmona and Delorme [19]. Recently van den Ban and Schlichtkrull
gave another proof of (4.28) based on the principle of induction
[12] involving the idea that both sides of (4.28) are
eigenfunctions depending meromorphically on $\la\in
(\fa_{Q,\fq})^*_\BC$. The formula (4.29) follows from the formula
(4.28) by comparing the coefficients of $a^{ts\la-\rho_R}$ in the
asymptotic expansions along $A_{R,\fq}^+$ for $v\in {\mathcal
W}_R.$ The formula (4.30) follows from the following functional
equation
\begin{equation}
C^0_{Q|P}(s:\la)\circ E^*(P:\la:x)=E^*(Q:s\la:x)
\end{equation}

\noindent obtained from the equation (4.28) and Maass-Selberg
relations.

\vskip 0.2cm Motivated by the formula (4.30) with $P=Q$, we define
\begin{equation*}
[{\mathcal S}(i\,\fa_{P,\fq}^*)\otimes {\mathcal
A}_{2,P}]^{W(\fa_{P,\fq})}
\end{equation*}

\noindent to be the subspace of ${\mathcal
S}(i\,\fa_{P,\fq}^*)\otimes {\mathcal A}_{2,P}$ consisting of the
functions $f$ satisfying the following transformation law
\begin{equation*}
f(s\la)=C^0_{P|P}(s:\la)f(\la),\quad \la\in i\,\fa_{P,\fq}^*,\
s\in W(\fa_{P,\fq}).
\end{equation*}
Then we see that ${\mathcal F}_P$ maps ${\mathcal C}(G/H:\tau)$
into $[{\mathcal S}(i\,\fa_{P,\fq}^*)\otimes {\mathcal
A}_{2,P}]^{W(\fa_{P,\fq})}$.

\vskip 0.3cm van den Ban and Schlichtkrull [13] and Delorme [21]
obtained the following. \vskip 0.2cm \noindent $ \textbf{The
Plancherel formula for Spherical functions\,([13,\,21]: }$ $
\textit{ The map} \ {\mathcal F}:=\sum_{P\in {\mathbb P}_\s}
{\mathcal F}_P$ $ \textit{is a topological linear isomorphism
from}$ ${\mathcal C}(G/H:\tau)$ $ \textit{onto}\ \oplus_{P\in
{\mathbb P}_\s}[{\mathcal S}(i\,\fa_{P,\fq}^*)\otimes {\mathcal
A}_{2,P}]^{W(\fa_{P,\fq})}.\ \textit{The}$\\ $ \textit{inverse
of}\ {\mathcal F}\ \textit{is given by}$
\begin{equation*}
{\mathcal T}:=\oplus_{P\in {\mathbb P}_\s}[W:W_P^*]\,{\mathcal
T}_P.
\end{equation*}

\noindent $ \textit{Moreover, for every}\ f\in C(G/H:\tau),$
\begin{equation*}
||f||^2_{L^2(G/H:\tau)}=\oplus_{P\in {\mathbb
P}_\s}[W:W_P^*]\,||{\mathcal F}_Pf||^2_{L^2}.
\end{equation*}

\end{section}

\vskip 1cm
%
%
\begin{section}{{\bf The Siegel-Jacobi Space}}
\setcounter{equation}{0}

\vskip 0.2cm In this section, we present some results and survey
some known results on the homogeneous space $\BH_n\times
\BC^{(m,n)}$ of non-reductive type which is important
arithmetically and geometrically. The theory of harmonic analysis
on the Siegel-Jacobi space is a vast work which is needed to be
studied in the future because it is important arithmetically and
geometrically. In fact, this work generalizes the theory of
harmonic analysis on the Siegel's fundamental domain which is
still under development and gives many important arithmetic and
geometric results. In this section, we deal with the theory of
harmonic analysis on the Siegel-Jacobi space briefly. In the near
future, I hope that I present a more detailed theory of harmonic
analysis on the Siegel-Jacobi space.

This section is organized as follows. In the subsection 5.1, we
discuss invariant differential operators on the Siegel upper half
plane $\BH_n$. We also present invariant metrics and their
Laplacians on $\BH_n$ obtained explicitly by Siegel and Maass
respectively. In the subsection 5.2, we discuss invariant
differential operators on the Siegel-Jacobi space $\BH_n\times
\BC^{(m,n)}$. Explicit invariant metrics and their Laplacians on
$\BH_n\times \BC^{(m,n)}$ obtained recently by the author are
presented. The author also present a partially Cayley transform of
$\BH_n\times \BC^{(m,n)}$ which gives a partially bounded
realization of $\BH_n\times \BC^{(m,n)}$ by $\BD_n\times
\BC^{(m,n)}$, where
$$\BD_n=\left\{ W\in \BC^{(n,n)}\,|\ W=\,^tW,\ I_n-{\overline
W}W>0\ \right\}$$ is the generalized unit disk of degree $n$.
Using this partial Cayley transform, the author also present
invariant metrics and their Laplacians on $\BD_n\times
\BC^{(m,n)}$ explicitly. In the subsection 5.3, we describe a
fundamental domain for the Siegel-Jacobi space explicitly using
the work of Siegel. In the subsection 5.4, we define Maass-Jacobi
forms on $\BH_n\times \BC^{(m,n)}$ which play an important role in
a spectral theory of a Laplacian on the Siegel-Jacobi space. In
the subsection 5.5, we mention a formal Eisenstein series briefly.
In the subsection 5.6, we give a brief remark on Fourier
expansions of Maass-Jacobi forms and present the problem to be
investigated in the future. In the subsection 5.7, we discuss
singular Jacobi forms which play an important role in the study of
the universal abelian variety, and also present the duality
theorem. In the subsection 5.8, we review a real analytic
Eisenstein series introduced by T. Arakawa [1]. In the subsection
5.9, we give a brief remark on spectral analysis of a Lapalcian on
the Siegel-Jacobi space. In the final subsection, we discuss
roughly the harmonic analysis on the Siegel-Jacobi space and the
Jacobi group modulo an arithmetic subgroup. We present the problem
to be investigated in the future.

\vskip 0.2cm\noindent $ \textbf{5.1.\ Invariant differential
operators on the Siegel upper half plane.}$ \vskip 0.3cm

For a given fixed positive integer $n$, we let
$${\mathbb H}_n=\,\{\,\Om\in \BC^{(n,n)}\,|\ \Om=\,^t\Om,\ \ \ \text{Im}\,\Om>0\,\}$$
be the Siegel upper half plane of degree $n$ and let
$$Sp(n,\BR)=\{ M\in \BR^{(2n,2n)}\ \vert \ ^t\!MJ_nM= J_n\ \}$$
be the symplectic group of degree $n$, where
$$J_n=\begin{pmatrix} 0&I_n \\
                   -I_n&0 \\ \end{pmatrix}.$$
$Sp(n,\BR)$ acts on $\BH_n$ transitively by
\begin{equation}M\cdot \Om=(A\Om+B)(C\Om+D)^{-1}, \end{equation}
where $M=\begin{pmatrix} A&B\\ C&D\end{pmatrix}\in Sp(n,\BR)$ and
$\Om\in \BH_n.$ \vskip 0.2cm For brevity, we write $G=Sp(n,\BR).$
We see that the stabilizer $K$ at $iI_n$ is given by
\begin{equation*}
K=\left\{ \begin{pmatrix}  \ A & B \\ -B & A \end{pmatrix}\in G
\Big|\ A+iB\in U(n)\,\right\}
\end{equation*}

\noindent and hence the map
\begin{equation*}
G/K\lrt \BH_n,\quad gK\mapsto g\cdot (iI_n)
\end{equation*}

\noindent is diffeomorphic (in fact, biholomorphic). Let $\fg$ be
the Lie algebra of $G$ and let $\fg_{\BC}$ be its
complexification. The Killing form $B$ for $\fg$ is given by
\begin{equation*}
B(X,Y)=2(n+1) \,\textrm{tr} (XY),\quad X,Y\in \fg.
\end{equation*}

Let $\theta$ be a Cartan involution of $G$ defined by
\begin{equation}
\theta (g)=\,{}^tg^{-1},\quad g\in G.
\end{equation}

\noindent $\fg$ has a Cartan decomposition
\begin{equation*}
\fg=\fk\oplus \fp,
\end{equation*}

\noindent where
\begin{equation*}
\fk=\left\{ \begin{pmatrix}  \ A & B \\ -B & A \end{pmatrix}\in
\BR^{(2n,2n)}\, \Big|\ A+\,^tA=0,\ B=\,{}^tB\,\right\}
\end{equation*}

\noindent and
\begin{equation*}
\fp=\left\{ \begin{pmatrix}   A & \ B \\ B & -A \end{pmatrix}\in
\BR^{(2n,2n)}\, \Big|\ A=\,^tA,\ B=\,{}^tB\,\right\}.
\end{equation*}

\noindent We note that $\fk$ is the Lie algebra of $K$ and that
$\fk$ and $\fp$ are orthogonal with respect to the Killing form
$B$. The vector space $\fp$ is identified with the tangent space
of $\BH_n$ at $iI_n$. The correspondence
\begin{equation*}
{\frac 12}\begin{pmatrix}  \ A & B \\ B & -A \end{pmatrix}\mapsto
A+iB
\end{equation*}

\noindent yields an isomorphism of $\fp$ onto the space $T_n$ of
symmetric $n\times n$ complex matrices. The differential map
$d\theta$ of $\theta$ extends complex linearly to $\fg_\BC$. The
$(\pm 1)$-eigenspaces of $d\theta$ are
\begin{equation*}
\fk_\BC=\left\{ \begin{pmatrix}  \ A & B \\ -B & A
\end{pmatrix}\in \BC^{(2n,2n)}\, \Big|\ A+\,^tA=0,\
B=\,{}^tB\,\right\}
\end{equation*}

\noindent and
\begin{equation*}
\fp_\BC=\left\{ \begin{pmatrix}   A & \ B \\ B & -A
\end{pmatrix}\in \BC^{(2n,2n)}\, \Big|\ A=\,^tA,\
B=\,{}^tB\,\right\}
\end{equation*}

\noindent respectively. We observe that $\fk_\BC$ and $\fp_\BC$
are the complexifications of $\fk$ and $\fp$ respectively.
$\fp_\BC$ has the following decomposition
\begin{equation*}
\fp_\BC=\fp_+\oplus \fp_-,
\end{equation*}

\noindent where
\begin{equation*}
\fp_+=\left\{ \begin{pmatrix}   X & \ iX \\ iX & -X
\end{pmatrix}\in \BC^{(2n,2n)}\, \Big|\ X=\,^tX\,\right\}
\end{equation*}

\noindent and
\begin{equation*}
\fp_-=\left\{ \begin{pmatrix}   Y & \ -iY \\ -iY & -Y
\end{pmatrix}\in \BC^{(2n,2n)}\, \Big|\ Y=\,^tY\,\right\}.
\end{equation*}

\noindent We observe that $\fp_+$ and $\fp_-$ are abelian
subalgebras of $\fg_\BC$. It is seen that
\begin{equation}
[\fk_\BC,\fk_\BC]\subset \fk_\BC,\quad \ [\fk_\BC,\fp_\BC]\subset
\fp_\BC,\quad
 \ [\fp_\BC,\fp_\BC]\subset \fk_\BC.
\end{equation}

\noindent Since $ \textrm{Ad}(k)X=kXk^{-1} \ (k\in K,\ X\in
\fg_\BC)$, we get the relation
\begin{equation*}
\textrm{Ad}(k)\fp_+\subset \fp_+,\quad \textrm{Ad}(k)\fp_-\subset
\fp_-.
\end{equation*}

\noindent For instance, if $k=\begin{pmatrix}  \ A & B \\ -B & A
\end{pmatrix}\in K,$ then
\begin{equation}
\textrm{Ad}(k)\begin{pmatrix}   X & \ \pm iX \\ \pm iX & -X
\end{pmatrix}=
\begin{pmatrix}   X' & \ \pm iX' \\ \pm iX' & -X' \end{pmatrix},\quad X=\,{}^tX,
\end{equation}

\noindent where
\begin{equation*}
X'=(A+iB)X\,{}^t(A+iB).
\end{equation*}

\noindent If we identify $\fp_-$ with $T_n$, then the action of
$K$ on $\fp_-$ is compatible with the natural representation
$\rho^{[1]}$ of $GL(n,\BC)$ on $T_n$ given by
\begin{equation}
\rho^{[1]}(g)(X):=gX\,{}^tg,\quad g\in GL(n,\BC),\ X\in T_n.
\end{equation}

\newcommand\PO{ {{\partial}\over {\partial \Omega}} }
\newcommand\PE{ {{\partial}\over {\partial \eta}} }
\newcommand\POB{ {{\partial}\over {\partial{\overline \Omega}}} }
\newcommand\PEB{ {{\partial}\over {\partial{\overline \eta}}} }
For a coordinate $\Om=(\om_{\mu\nu})\in {\mathbb H}_n$, we put
\begin{eqnarray*}
\Om\,&=&\,X\,+\,iY,\quad\ \ X\,=\,(x_{\mu\nu}),\quad\ \
Y\,=\,(y_{\mu\nu})
\ \ \text{real},\\
d\Om\,&=&\,(d\om_{\mu\nu}),\quad\ \ d{\overline
\Omega}=(d{\overline\omega}_{\mu\nu}),
\end{eqnarray*}
\begin{eqnarray*}
\PO\,=\,\left(\, { {1+\delta_{\mu\nu}} \over 2}\, {
{\partial}\over {\partial \om_{\mu\nu}} } \,\right),\quad
\POB\,=\,\left(\, { {1+\delta_{\mu\nu}}\over 2} \, {
{\partial}\over {\partial {\overline \om}_{\mu\nu} }  } \,\right),
\end{eqnarray*}

\noindent where $\delta_{ij}$ denotes the Kronecker delta symbol.
C. L. Siegel [43] introduced the symplectic metric
\begin{eqnarray}
ds_n^2= \, \textrm{tr} (Y^{-1}d\Om\,Y^{-1}d{\overline \Om})
\end{eqnarray}

\noindent which is invariant under the action (5.1) of $G$. H.
Maass [34] proved that the differential operator
\begin{eqnarray}
\Delta_n=\,4\, \textrm{tr}  \left( Y \,\,^t\!\left( Y\POB\right)
\PO\right)
\end{eqnarray}
is the Laplacian of ${\mathbb H}_n$ for the symplectic metric
$ds_n^2.$

\vskip 0.2cm According to Harish-Chandra, we see that the algebra
$\BD(\BH_n)$ of all $G$-invariant differential operators on
$\BH_n$ is generated by $n$ algebraically independent commuting
differential operators on $\BH_n$. Here we note that the real rank
of $G$ is $n$. Therefore $\BD(\BH_n)$ is isomorphic to the
polynomial ring $\BC [x_1,\cdots,x_n].$ We see from (5.4) that $K$
acts on $T_n$ by
\begin{equation}
k\cdot X=kX\,{}^tk,\quad k\in K,\ X\in T_n.
\end{equation}

\noindent This action induces naturally the action $\tau$ of $K$
on the polynomial algebra $ \textrm{Pol}(T_n)$ of $T_n$. We denote
by $ \textrm{Pol}(T_n)^K$ the subalgebra of $ \textrm{Pol}(T_n)$
consisting of all $K$-invariants for the action $\tau$ of $K$. The
following $K$-invariant inner product $(\ ,\ )$ on the complex
vector space $T_n$ defined by
\begin{equation*}
(X,Y):= \textrm{tr}(X{\overline Y}),\quad X,Y\in T_n
\end{equation*}

\noindent gives an isomorphism
\begin{equation*}
T_n\cong T_n^*,\quad X\mapsto f_X,\quad X\in T_n,
\end{equation*}

\noindent where $T_n^*$ denotes the dual space of $T_n$ and $f_X$
is the linear functional on $T_n$ defined by
\begin{equation*}
f_X(Y):=(Y,X)= \textrm{tr}(Y{\overline X}),\quad Y\in T_n.
\end{equation*}

Let $E_{ij}$ be the $n\times n$ matrix with entry $1$ where the
$i$-th row and the $j$-th column meet, and all other entries $0$.
We put
\begin{equation*}
e_i=E_{ii}\,(1\leq i\leq n),\quad e_{ij}={1\over {\sqrt
2}}(E_{ij}+E_{ji})\,(1\leq i< j\leq n).
\end{equation*}

\noindent It is easy to see that the matrices $e_i\,(1\leq i\leq
n),\, e_{ij}\,(1\leq i< j\leq n)$ form an orthonormal basis for
$T_n$ with respect to the inner product $(\,,\,)$. We choose a
basis $e_i\,(1\leq i\leq n),\, e_{ij}\,(1\leq i< j\leq n)$ for
$T_n$. Then we get a canonical linear bijection of $S(T_n)^K$ onto
$\BD(\BH_n)$, where $S(T_n)^K$ is the subalgebra of the symmetric
algebra $S(T_n)$ of $T_n$ consisting of all $K$-invariants for the
natural action of $K$ on $S(T_n)$. Identifying $T_n$ with $T_n^*$
by the above isomorphism, we get a canonical linear bijection
\begin{equation*}
\Psi:\textrm{Pol}(T_n)^K \lrt \BD(\BH_n)
\end{equation*}

\noindent of $\textrm{Pol}(T_n)^K$ onto $\BD(\BH_n)$. For more
detail, we refer to [29], p.\,287. In fact, $\textrm{Pol}(T_n)^K$
is generated by
\begin{equation}
p_k(X):=\s ((X{\overline X})^k),\quad 1\leq k\leq n
\end{equation}

\noindent and $\Psi(p_1(X))=\Delta_n.$ As far as I know, it seems
that for general $n$, the invariant differential operators
$\Psi(p_k(X))\, (2\leq k\leq n)$ were not written explicitly so
far.

\vskip 0.2cm Let
\begin{equation*}
\BD_n=\left\{ W\in \BC^{(n,n)}\,|\ W=\,{}^tW,\ I_n-\OW W >
0\,\right\}
\end{equation*}
be the generalized unit disk of degree $n$. The Cayley transform
$\Phi:\BD_n\lrt \BH_n$ defined by
\begin{equation}
\Phi(W):=i(I_n+W)(I_n-W)^{-1},\quad W\in \BD_n
\end{equation}

\noindent is a biholomorphic mapping of $\BD_n$ onto $\BH_n$ which
gives the bounded realization of $\BH_n$ by $\BD_n$
(cf.\,[43],\,pp. 281-283). Let
\begin{equation*}
T={1\over {\sqrt{2}} }  \begin{pmatrix}   I_n & \ iI_n \\ iI_n &
-iI_n \end{pmatrix}
\end{equation*}

\noindent be the $2n\times 2n$ matrix represented by $\Phi$. Then
we see that
\begin{equation*}
T^{-1}G T=\left\{  \begin{pmatrix}   P & Q \\ \OQ & \OP
\end{pmatrix}\in SU(n,n)\,\Big|\ {}^tP\OP-{}^t\OQ Q=I_n,\
{}^tP\OQ=\,{}^t\OQ P\,\right\}.
\end{equation*}

\noindent For brevity, we set
\begin{equation*}
G_*=T^{-1}G T.
\end{equation*}

\vskip 0.2cm \noindent If the case $n=1$, we note that
$G_*=SU(1,1)$. If $n>1,$ then $G_*$ is a {\it proper} subgroup of
$SU(n,n).$ In fact, since ${}^tTJ_nT=-\,i\,J_n$, we get
$$G_*=\left\{\,h\in SU(n,n)\,|\ {}^thJ_nh=J_n\,\right\}.$$

Let
\begin{equation*}
P^+=\left\{\begin{pmatrix} I_n & Z\\ 0 & I_n
\end{pmatrix}\,\Big|\ Z=\,{}^tZ\in\BC^{(n,n)}\,\right\}
\end{equation*}
be the $P^+$-part of the complexification of $G_*\subset SU(n,n).$

Since the Harish-Chandra decomposition of an element
$\begin{pmatrix} P & Q\\ {\overline Q} & {\overline P}
\end{pmatrix}$ in $G_*^J$ is
\begin{equation*}
\begin{pmatrix} P & Q\\ \OQ & \OP
\end{pmatrix}=\begin{pmatrix} I_n & Q\OP^{-1}\\ 0 & I_n
\end{pmatrix} \begin{pmatrix} P-Q\OP^{-1}\OQ & 0\\ 0 & \OP
\end{pmatrix} \begin{pmatrix} I_n & 0\\ \OP^{-1}\OQ & I_n
\end{pmatrix},
\end{equation*}
the $P^+$-component of the following element
$$\begin{pmatrix} P & Q\\ \OQ & \OP
\end{pmatrix}   \cdot\begin{pmatrix} I_n & W\\ 0 & I_n
\end{pmatrix},\quad W\in \BD_n$$ of the complexification of $G_*^J$ is
given by
\begin{equation*}
 \begin{pmatrix} I_n & (PW+Q)(\OQ W+\OP)^{-1}
\\ 0 & I_n
\end{pmatrix}.
\end{equation*}

We note that $Q\OP^{-1}\in \BD_n.$ We get the Harish-Chandra
embedding of $\BD_n$ into $P^+$\,(cf.\,[30],\,p.\,155). Therefore
we see that $G_*$ acts on $\BD_n$ transitively by
\begin{equation}
 \begin{pmatrix}   P & Q \\ \OQ & \OP \end{pmatrix}\cdot W=(PW+Q)(\OQ W+\OP)^{-1},
\end{equation}

\noindent where $\begin{pmatrix}   P & Q \\ \OQ & \OP
\end{pmatrix}\in G_*$ and $W\in \BD_n$. We can show that the
action (5.11) of $G_*$ on $\BD_n$ is compatible with the action
(5.1) of $G$ through the Cayley transform $\Phi.$

\vskip 0.2cm For a coordinate $W\in \BD_n$ with $W=(w_{\mu\nu})$,
we put
\begin{eqnarray*}
dW=(dw_{\mu\nu}),\quad\ \ d{\overline W}\,=\,(d{\overline
w}_{\mu\nu})
\end{eqnarray*}
and
\begin{eqnarray*}
\PW\,=\,\left(\, { {1+\delta_{\mu\nu}} \over 2}\, {
{\partial}\over {\partial w_{\mu\nu}} } \,\right),\quad
\PWB\,=\,\left(\, { {1+\delta_{\mu\nu}}\over 2} \, {
{\partial}\over {\partial {\overline w}_{\mu\nu} }  } \,\right).
\end{eqnarray*}

\noindent Using the Cayley transform $\Phi$, we can see
(cf.\,[43]) that
\begin{equation}
ds_*^2=4 \, \textrm{tr} ((I_n-W{\overline W})^{-1}dW\,(I_n-\OW
W)^{-1}d\OW\,)\end{equation} is a $G_*$-invariant Riemannian
metric on $\BD_n$ and H. Maass [34] showed that its Laplacian is
given by
\begin{equation}
\Delta_*=\, \textrm{tr} \left( (I_n-W\OW)\,{ }^t\!\left(
(I_n-W\OW) \PWB\right)\PW\right).\end{equation}

\vskip 0.3cm \noindent $ \textbf{5.2.\ Invariant differential
operators on the Siegel-Jacobi space.}$ \vskip 0.3cm

 For two positive integers $n$ and $m$, we consider
the Heisenberg group
$$H_{\BR}^{(n,m)}=\{\,(\l,\mu;\k)\,|\ \l,\mu\in \BR^{(m,n)},\ \k\in \BR^{(m,m)},\
\k+\mu\,^t\l\ \text{symmetric}\ \}$$ endowed with the following
multiplication law
$$(\l,\mu;\k)\circ (\l',\mu';\k')=(\l+\l',\mu+\mu';\k+\k'+\l\,^t\mu'-
\mu\,^t\l').$$ It is easy to see that $Sp(n,\BR)$ acts on
$H_{\BR}^{(n,m)}$ on the right as group automorphisms\,:
$$(\lambda,\mu;\kappa)\cdot M:=((\lambda,\mu)M;\kappa),\quad M\in
Sp(n,\BR),\ (\lambda,\mu;\kappa)\in H_{\BR}^{(n,m)}.$$ We observe
that the pairing on $\BR^{(m,n)}$ defined by
\begin{equation*}
((\lambda,\mu),(\lambda',\mu'))\mapsto \,
(\lambda,\mu)J_n\,^t(\lambda',\mu')=\lambda\,^t\mu'-\mu\,^t\lambda'
\end{equation*}

\noindent is invariant under the above right action of $Sp(n,\BR)$
on $H_{\BR}^{(n,m)}$. \vskip 0.3cm The $\textsf{Jacobi group}$
\begin{equation}
G^J:=Sp(n,\BR)\ltimes H_{\BR}^{(n,m)}\quad ( \textrm{the\
semidirect\ product}) $$ is defined to be the semidirect product
endowed with the following multiplication law
$$(M,(\lambda,\mu;\kappa))\cdot(M',(\lambda',\mu';\kappa')) =\,
(MM',(\tilde{\lambda}+\lambda',\tilde{\mu}+ \mu';
\kappa+\kappa'+\tilde{\lambda}\,^t\!\mu'
-\tilde{\mu}\,^t\!\lambda'))
\end{equation}

\noindent with $M,M'\in Sp(n,\BR),\
(\lambda,\mu;\kappa),\,(\lambda',\mu';\kappa') \in
H_{\BR}^{(n,m)}$ and
$(\tilde{\lambda},\tilde{\mu})=(\lambda,\mu)M'$. Then we get the
{\it canonical action} of $G^J$ on $\BH_n\times
\BC^{(m,n)}$\,(cf.\,[46,\,47,\,48]) defined by
\begin{equation}(M,(\lambda,\mu;\kappa))\cdot (\Om,Z)=(M\cdot\Om,(Z+\lambda \Om+\mu)
(C\Om+D)^{-1}), \end{equation}
where $M=\begin{pmatrix} A&B\\
C&D\end{pmatrix} \in Sp(n,\BR),\ (\lambda,\mu; \kappa)\in
H_{\BR}^{(n,m)}$ and $(\Om,Z)\in \BH_n\times \BC^{(m,n)}.$ We note
that the action (5.15) is transitive. The stabilizer $K^J$ of
$G^J$ at $(iI_n,0)$ is given by
\begin{equation*}
K^J=\left\{ (k,(0,0;\ka))\,|\ k\in K,\ \ka=\,^t\ka\in
\BR^{(m,m)}\,\right\}.
\end{equation*}

\noindent Therefore $\Hnm\cong G^J/K^J$ is a homogeneous space of
$ \textit{non-reductive type}$. The Lie algebra $\fg^J$ of $G^J$
has a decomposition
\begin{equation*}
\fg^J=\fk^J+\fp^J,
\end{equation*}

\noindent where
\begin{equation*}
\fk^J=\left\{ (X,(0,0,\ka))\,|\ X\in \fk,\ \ka=\,^t\ka\in
\BR^{(m,m)}\,\right\}
\end{equation*}

\noindent
\begin{equation*}
\fp^J=\left\{ (Y,(P,Q,0))\,|\ Y\in \fp,\ P,Q\in
\BR^{(m,n)}\,\right\}.
\end{equation*}

\noindent Thus the tangent space of the homogeneous space $\Hnm$
at $(iI_n,0)$ is identified with $\fp^J$. We define a complex
structure $I^J$ on the tangent space $\fp^J$ of $\Hnm$ at
$(iI_n,0)$ by
\begin{equation*}
I^J\left( \begin{pmatrix}  \ Y & \ X \\ X & -Y
\end{pmatrix},(P,Q)\right)= \left( \begin{pmatrix}   X &  -Y \\ -Y
& -X \end{pmatrix},(Q,-P)\right),
\end{equation*}

\noindent where $X=\,{}^tX\in \BR^{(m,n)}$ and $Y=\,{}^tY\in
\BR^{(m,n)}$. Identifying $\BR^{(m,n)}\times \BR^{(m,n)}$ with
$\BC^{(m,n)}$ via
\begin{equation*}
(P,Q)\mapsto iP+Q,\quad P,Q\in \BR^{(m,n)},
\end{equation*}

\noindent we may regard the complex structure $I^J$ as a real
linear map
\begin{equation*}
I^J(X+iY,Q+iP)=(-Y+iX,-P+iQ),
\end{equation*}

\noindent where $X+iY\in T_n$ and $Q+iP\in \BC^{(m,n)}.$ Obviously
$I^J$ extends complex linearly on the complexification $\fp_\BC^J$
of $\fp^J$. Then $\fp_\BC^J$ has a decomposition
\begin{equation*}
\fp_\BC^J=\fp^J_+\oplus \fp^J_-,
\end{equation*}

\noindent where $\fp^J_+$\,(resp.\,$ \fp^J_-$) denotes the
$(+i)$-eigenspace (resp. $(-i)$-eigenspace) of $I^J$. Precisely,
both $\fp^J_+$ and $ \fp^J_-$ are given by
\begin{equation*}
\fp_+^J=\left\{ \left( \begin{pmatrix}   X &  iX \\ iX & -X
\end{pmatrix}, (P,iP)\right)\, \Big|\ X=\,^tX\in \BC^{(n,n)},\
P\in \BC^{(m,n)}\,\right\}
\end{equation*}

\noindent and
\begin{equation*}
\fp_-^J=\left\{ \left( \begin{pmatrix}   X &  -iX \\ -iX & -X
\end{pmatrix}, (P,-iP)\right)\, \Big|\ X=\,^tX\in \BC^{(n,n)},\
P\in \BC^{(m,n)}\,\right\}.
\end{equation*}

\noindent Let $\fk^J$ be the Lie algebra of $K^J$ and let
$\fk^J_\BC$ be its complexification. We can see that
\begin{equation}
[\fk^J,\fk^J]\subset \fk^J,\quad [\fk^J,\fp^J]\subset \fp^J
\end{equation}

\noindent and
\begin{equation}
[\fk^J_\BC,\fp^J_+]\subset \fp_+^J,\quad
[\fk^J_\BC,\fp^J_-]\subset \fp^J_-.
\end{equation}

\noindent The complexification $\fg^J_\BC$ of $\fg^J$ has a
decomposition $\fg_\BC^J=\fk_\BC^J+\fp_\BC^J.$

\vskip 0.2cm For brevity, we write $\BH_{n,m}:=\Hnm.$ For a
coordinate $(\Om,Z)\in \Hnm$ with $\Om=(\omega_{\mu\nu})$ and
$Z=(z_{kl})$, we put $d\Om,\,d{\overline \Om},\,\PO,\,\POB$ as
before and set
\begin{eqnarray*}
Z\,&=&U\,+\,iV,\quad\ \ U\,=\,(u_{kl}),\quad\ \ V\,=\,(v_{kl})\ \
\text{real},\\
dZ\,&=&\,(dz_{kl}),\quad\ \ d{\overline Z}=(d{\overline z}_{kl}),
\end{eqnarray*}

\begin{eqnarray*}
\PO\,=\,\left(\, { {1+\delta_{\mu\nu}} \over 2}\, {
{\partial}\over {\partial \om_{\mu\nu}} } \,\right),\quad
\POB\,=\,\left(\, { {1+\delta_{\mu\nu}}\over 2} \, {
{\partial}\over {\partial {\overline \om}_{\mu\nu} }  } \,\right),
\end{eqnarray*}

$$\PZ=\begin{pmatrix} {\partial}\over{\partial z_{11}} & \hdots &
 {\partial}\over{\partial z_{m1}} \\
\vdots&\ddots&\vdots\\
 {\partial}\over{\partial z_{1n}} &\hdots & {\partial}\over
{\partial z_{mn}} \end{pmatrix},\quad \PZB=\begin{pmatrix}
{\partial}\over{\partial {\overline z}_{11} }   &
\hdots&{ {\partial}\over{\partial {\overline z}_{m1} }  }\\
\vdots&\ddots&\vdots\\
{ {\partial}\over{\partial{\overline z}_{1n} }  }&\hdots &
 {\partial}\over{\partial{\overline z}_{mn} }  \end{pmatrix}.$$

\newcommand\bw{d{\overline W}}
\newcommand\bz{d{\overline Z}}
\newcommand\be{d{\overline \eta}}
\newcommand\bo{d{\overline \Omega}}

\vskip 0.3cm Recently the author [51] proved the following
theorem.

\vskip 0.5cm
\begin{theorem}
The following metric
\begin{eqnarray}
ds_{n,m}^2=&\,\textrm{tr}\left(Y^{-1}d\Om\,Y^{-1}d{\overline
\Om}\right)\,+
\,\textrm{tr}\left(Y^{-1}\,^tV\,V\,Y^{-1}d\Om\,Y^{-1}
\bo\right) \nonumber\\
&\ \ \ \ +\,\textrm{tr}\left(Y^{-1}\,^t(dZ)\,\bz\right)\\
&\ \ -\textrm{tr}\left(V\,Y^{-1}d\Om\,Y^{-1}\,^t(\bz)\,
+\,V\,Y^{-1}\bo\, Y^{-1}\,^t(dZ)\,\right)\nonumber
\end{eqnarray}
is a Riemannian metric on $\BH_{n,m}$ which is invariant under the
action (5.15) of the Jacobi group $G^J$. The Laplacian
$\Delta_{n,m}$ of $(\BH_{n,m},ds^2_{n,m})$ is given by
\begin{eqnarray}
\Delta_{n,m}\,=\,& 4\, \textrm{tr} \left(\,Y\,\,
^t\!\left(Y\POB\right)\PO\,\right)\,+\,
4\,\textrm{tr}\left(\, Y\,\PZ\,\,{}^t\!\left( \PZB\right)\,\right) \nonumber \\
&\ \ \ \ +\,4\,\textrm{tr}\left(\,VY^{-1}\,^tV\,\,^t\!\left(Y\PZB\right)\,\PZ\,\right)\\
&\ \
+\,4\,\textrm{tr}\left(V\,\,^t\!\left(Y\POB\right)\PZ\,\right)+\,4\,\textrm{tr}\left(\,^tV\,\,^t\!\left(Y\PZB\right)\PO\,\right).\nonumber
\end{eqnarray}

\noindent The following differential form
$$dv=\,\left(\,\text{det}\,Y\,\right)^{-(n+m+1)}[dX]\w [dY]\w
[dU]\w [dV]$$ is a $G^J$-invariant volume element on $\BH_{n,m}$,
where
$$[dX]=\w_{\mu\leq\nu}dx_{\mu\nu},\quad [dY]=\w_{\mu\leq\nu}
dy_{\mu\nu},\quad [dU]=\w_{k,l}du_{kl}\quad \text{and} \quad
[dV]=\w_{k,l}dv_{kl}.$$
\end{theorem}

\newcommand\SJ{{\mathbb H}_n\times {\mathbb C}^{(m,n)}}
\newcommand\Hn{{\mathbb H}_n}
\noindent $\textit{Sketch of Proof.}$ Let $g=(M,(\la,\mu;\k))$ be
an element of $G^J$ with $M=\begin{pmatrix} A&B\\ C&D
\end{pmatrix}\in Sp(n,\BR)$ and $(\Omega,Z)\in \BH_{n,m}$ with $\Omega\in\BD_n$ and $Z\in\Cmn.$
If we put $(\Omega_*,Z_*):=g\cdot(\Omega,Z),$ then we have
\begin{eqnarray*}
&\Omega_*=M\cdot \Omega=(A\Omega+B)(C\Omega+D)^{-1},\\
&Z_*=(Z+\la \Omega+\mu)(C\Omega+D)^{-1}.
\end{eqnarray*}
Thus we obtain
\begin{eqnarray}
d\Omega_*=d\Omega[(C\Omega+D)^{-1}]={}^t\!(C\Omega+D)^{-1}d\Omega(C\Omega+D)^{-1}
\end{eqnarray}
and
\begin{eqnarray}
\quad\quad\quad dZ_*=dZ(C\Omega+D)^{-1}+ \{ \l-(Z+\l
\Omega+\mu)(C\Omega+D)^{-1}C \} d\Omega(C\Omega+D)^{-1}.
\end{eqnarray}
Here we used the following facts that
$$d(C\Omega+D)^{-1}=-(C\Omega+D)^{-1}C\,d\Omega(C\Omega+D)^{-1}$$
and that $(C\Omega+D)^{-1}C$ is symmetric.

\renewcommand\o{\overline}
We put
$$\Omega_*=X_*+iY_*,\quad Z_*=U_*+iV_*,\quad X_*,Y_*,U_*,V_*\;\text{real}.$$
From [35], p.33 or [44], p.128, we know that
\begin{eqnarray}Y_*=Y\{(C\Omega+D)^{-1}\}={}^t\!(C\overline{\Omega}+D)^{-1}\,Y(C\Omega+D)^{-1}.\end{eqnarray}

For brevity, we write
\begin{eqnarray*}
&&(a)=\s(Y^{-1}d\Omega\,Y^{-1}d\overline{\Omega}),\\
&&(b)=\s(Y^{-1}{}^tVVY^{-1}d\Omega\,Y^{-1}d\overline{\Omega}),\\
&&(c)=\s(Y^{-1}{}^t(dZ)d\overline{Z}),\\
&&(d)=-\s(V\,Y^{-1}d\Omega\,Y^{-1}\,{}^t(d\overline{Z})\,+\,V\,Y^{-1}d\overline{\Omega}\,Y^{-1}\,{}^t(dZ))
\end{eqnarray*}
and
\begin{eqnarray*}
&&(a)_*=\s(Y^{-1}_*d\Omega_*\,Y_*^{-1}d\overline{\Omega}_*),\\
&&(b)_*=\s(Y_*^{-1}\,{}^tV_*V_*Y_*^{-1}\,d\Omega_*\,Y_*^{-1}d\overline{\Omega}_*),\\
&&(c)_*=\s(Y_*^{-1}\,{}^t(dZ_*)\,d\overline{Z}_*),\\
&&(d)_*=-\s(V_*\,Y_*^{-1}d\Omega_*\,Y_*^{-1}\,{}^t(d\overline{Z}_*)+\,V_*\,Y_*^{-1}d\overline{\Omega}_*\,Y_*^{-1}\,{}^t(dZ_*)).
\end{eqnarray*}

Using the formulas (5.20)-(5.22), we can prove that
$$ (a)=(a)_*\quad \textrm{and}\quad (b)+(c)+(d)=(b)_* +(c)_*
+(d)_*.$$ For a complete proof, we refer to [51]. Therefore
$ds_{n,m}^2$ is invariant under the action of $G^J$. In
particular, for $(\Omega,Z)=(iI_n,0),$ we have
\begin{eqnarray*}
ds_{n,m}^2&=&  \textrm{tr}(d\Omega\,d\overline{\Omega})+ \textrm{tr}\left({}^t(dZ)\,d\overline{Z} \right)\\
&=&\sum_{\mu=1}^n(dx_{\mu\mu}^2+dy_{\mu\mu}^2)+2\sum_{1\leq\mu<\nu\leq
n}
(dx_{\mu\nu}^2+dy_{\mu\nu}^2)\\
&&\quad+\sum_{ 1\leq k\leq m,\, 1\leq l\leq n}
(du_{kl}^2+dv_{kl}^2),
\end{eqnarray*}
which is clearly positive definite. Since $G^J$ acts on $\mathbb
H_{n,m}$ transitively, $ds_{n,m}^2$ is positive definite
everywhere in $\mathbb H_{n,m}.$

\newcommand\POZ{ {{\partial}\over {\partial \Omega}} }
\newcommand\PZW{ {{\partial}\over {\partial Z}} }
\newcommand\POOB{ { {\partial}\over {\partial{\overline {\Omega}} }
}}
\newcommand\PZZB{ {{\partial}\over {\partial{\overline Z}}} }

\vskip 0.2cm We can show that
\begin{eqnarray}
\frac{\partial}{\partial \Omega_*}&=&(C\Omega+D)\,\, {}^t\!\left\{
(C\Omega+D)\,\POZ\right\} \\
&&\quad +\,(C\Omega+D)\,\,{}^t\!\left\{
(C\,{}^tZ+C\,\,{}^t\mu-D\,\,{}^t\la)\,\,{}^t\!\left(\PZW\right)\right\}\nonumber
\end{eqnarray}
and
\begin{equation}\frac{\partial}{\partial Z_*}=(C\Omega+D)\PZW.\end{equation}
For brevity, we put
\begin{eqnarray*}
&&(\alpha):=4\, \textrm{tr} \left(\,Y\,\,
^t\!\left(Y\POOB\right)\POZ\,\right),\\
&&(\beta):=\;4\,\textrm{tr}\left(\, Y\,\PZW\,\,{}^t\!\left( \PZZB\right)\,\right),\\
&&(\g):=\;4\,\textrm{tr}\left(\,VY^{-1}\,^tV\,\,^t\!\left(Y\PZZB\right)\,\PZW\,\right),\\
&&(\delta):=\;4\,\textrm{tr}\left(V\,\,^t\!\left(Y\POOB\right)\PZW\,\right)
\end{eqnarray*}
and
\begin{eqnarray*}
(\epsilon):=\;4\,\textrm{tr}\left(\,^tV\,\,^t\!\left(Y\PZZB\right)\POZ\,\right).
\end{eqnarray*}
We also set
\begin{eqnarray*}
&&(\alpha)_*:=4\,\s\left(\,Y_*\,\,
^t\!\left(Y_*\POOB_*\right)\POZ_*\,\right),\\
&&(\beta)_*:=\;4\,\textrm{tr}\left(\, Y_*\,\PZW_*\,\,{}^t\!\left( \PZZB_*\right)\,\right),\\
&&(\g)_*:=\;4\,\textrm{tr}\left(\,V_*\,Y_*^{-1}\,^tV_*\,\,^t\!\left(Y_*\PZZB_*\right)\,\PZW_*\,\right),\\
&&(\delta)_*:=\;4\,\textrm{tr}\left(V_*\,\,^t\!\left(Y_*\POOB_*\right)\PZW_*\,\right)
\end{eqnarray*}
and
\begin{eqnarray*}
(\epsilon)_*:=\;4\,\textrm{tr}\left(\,^tV_*\,\,^t\!\left(Y_*\PZZB_*\right)\POZ_*\,\right).
\end{eqnarray*}

By a complicated computation using the formulas (5.23)-(5.24), we
can prove that
$$(\beta)=(\beta)_*\quad \textrm{and}\quad
(\alpha)+(\g)+(\delta)+(\epsilon)=(\alpha)_*+(\g)_*+(\delta)_*+(\epsilon)_*.$$
For a complete proof, we refer to [51]. Hence $\Delta_{n,m}$ is
invariant under the action (5.15) of $G^J$. In particular, for
$(\Omega,Z)=(iI_n,0),$ the differential operator $4\,\Delta_{n,m}$
coincides with the Laplacian for the metric $ds^2_{n,m}$. It
follows from the invariance of $\Delta_{n,m}$ under the action
(5.15) and the transitivity of the action of $G^J$ on $\Hnm$ that
$4\,\Delta_{n,m}$ is the Laplacian of $(\Hnm,\,ds_{n,m}^2).$ The
invariance of the differential form $dv$ follows from the fact
that the following differential form
$$(\,\det Y\,)^{-(n+1)}[dX]\wedge [dY]$$
is invariant under the action (1.1) of
$Sp(n,\BR)$\,(cf.\,[44],\,p.\,130). \hfill$\square$

\vskip 0.3cm \noindent {\bf Remark 5.1.} In the special case
$n=m=1,$ by a direct computation, we see that the scalar curvature
of $(\BH_{1,1},\, ds^2_{1,1})$ is -3.

\vskip 0.4cm \noindent {\bf Remark 5.2.} We let ${\mathbb M}_1$
and ${\mathbb M}_2$ be the differential operators on $\BH_{n,m}$
defined by
\begin{equation*}
{\mathbb M}_1=\,4\,\textrm{tr}\left(\, Y\,\PZ\,\,{}^t\!\left(
\PZB\right)\,\right)
\end{equation*}

\noindent and
\begin{equation*}
{\mathbb M}_2=\Delta_{n,m}-{\mathbb M}_1.
\end{equation*}

\noindent Then ${\mathbb M}_1$ and ${\mathbb M}_2$ are invariant
under the action (5.15) of $G^J.$

\vskip 0.2cm  We can identify an element $g=(M,(\la,\mu;\kappa))$
of $G^J,\ M=\begin{pmatrix} A&B\\
C&D\end{pmatrix}\in Sp(n,\BR)$ with the element
\begin{equation*}
\begin{pmatrix} A & 0 & B & A\,^t\mu-B\,^t\la  \\ \la & I_m & \mu
& \kappa \\ C & 0 & D & C\,^t\mu-D\,^t\la \\ 0 & 0 & 0 & I_m
\end{pmatrix}
\end{equation*}
of $Sp(m+n,\BR).$ \vskip 0.3cm We set
\begin{equation*}
T_*={1\over {\sqrt 2}}\,
\begin{pmatrix} I_{m+n} & I_{m+n}\\ iI_{m+n} & -iI_{m+n}
\end{pmatrix}.
\end{equation*}
We now consider the group $G_*^J$ defined by
\begin{equation*}
G_*^J:=T_*^{-1}G^JT_*.
\end{equation*}
If $g=(M,(\la,\mu;\kappa))\in G^J$ with $M=\begin{pmatrix} A&B\\
C&D\end{pmatrix}\in Sp(n,\BR)$, then $T_*^{-1}gT_*$ is given by
\begin{equation*}
T_*^{-1}gT_*=
\begin{pmatrix} P_* & Q_*\\ {\overline Q}_* & {\overline P}_*
\end{pmatrix},
\end{equation*}
where
\begin{equation*}
P_*=
\begin{pmatrix} P & {\frac 12} \left\{ Q\,\,{}^t(\la+i\mu)-P\,\,{}^t(\la-i\mu)\right\}\\
{\frac 12} (\la+i\mu) & I_m+i{\frac \kappa 2}
\end{pmatrix},
\end{equation*}

\begin{equation*}
Q_*=
\begin{pmatrix} Q & {\frac 12} \left\{ P\,\,{}^t(\la-i\mu)-Q\,\,{}^t(\la+i\mu)\right\}\\
{\frac 12} (\la-i\mu) & -i{\frac \kappa 2}
\end{pmatrix},
\end{equation*}
and $P,\,Q$ are given by the formulas
\begin{equation}
P= {\frac 12}\,\left\{ (A+D)+\,i\,(B-C)\right\}
\end{equation}
and
\begin{equation}
 Q={\frac
12}\,\left\{ (A-D)-\,i\,(B+C)\right\}.
\end{equation}

\noindent Thus we can see that $G^J_*$ is of the form
\begin{equation*}
G^J_*=\left\{ (M_*,(\xi,{\overline \xi};\,i\ka)\,|\ M_*\in G_*,\
\xi\in\Cmn,\ \ka\in\BR^{(m,m)}\,\right\}
\end{equation*}

\noindent with the multiplication (5.14). We get the canonical
action of $G^J_*$ on $\Dnm$ defined by
\begin{equation}
\left(\begin{pmatrix} P & Q\\
{\overline Q} & {\overline P}
\end{pmatrix},\left( \xi, {\overline \xi};\,i\kappa\right)\right)\cdot
(W,\eta)=((PW+Q)(\OQ W+\OP)^{-1},(\eta+\xi W+{\overline\xi})(\OQ
W+\OP)^{-1}),\end{equation}

\noindent where $\begin{pmatrix} P & Q\\
{\overline Q} & {\overline P}
\end{pmatrix}\in G_*,\ \xi\in\Cmn,\ \ka\in\BR^{(m,m)}, \
W\in\BD_n$ and $\eta\in\Cmn.$ For brevity, we write
$\BD_{n,m}:=\Dnm$. We define the $ \textsf{partial Cayley
transform}\ \Phi_*:\BD_{n,m}\lrt \BH_{n,m}$ by
\begin{equation}
\Phi_*(W,\eta):=\left((i(I_n+W)(I_n-W)^{-1},\,2i\eta
(I_n-W)^{-1}\right),\quad (W,\eta)\in \BD_{n,m}.
\end{equation}

\noindent Then $\Phi_*$ is a biholomorphic mapping of $\BD_{n,m}$
onto $\BH_{n,m}$ which gives a partially bounded realization of
$\BH_{n,m}$ by $\BD_{n,m}$. Recently the author [52] proved the
following theorem.

\vskip 0.3cm
\begin{theorem}
The action (5.15) of
$G^J$ on $\BH_{n,m}$ is compatible with the action (5.27) of
$G^J_*$ on $\BD_{n,m}$ through the partial Cayley transform
$\Phi_*$ defined by (5.28). In other words, if $g_0\in G^J$ and
$(W,\eta)\in \BD_{n,m}$, then
\begin{equation*}
g_0\cdot \Phi_*(W,\eta)=\Phi_*(g_*\cdot (W,\eta)),
\end{equation*}

\noindent where $g_*=T_*^{-1}g_0 T_*.$ The inverse of $\Phi_*$ is
given by
\begin{equation*}
\Phi_*^{-1}(\Omega,Z)=\left((\Omega-iI_n)(\Omega+iI_n)^{-1},\,Z(\Omega+iI_n)^{-1}\right).
\end{equation*}
\end{theorem}
\vskip 0.1cm \noindent $\textit{Sketch of Proof.}$
 Let
 \begin{equation*}
 g_0=\left(\begin{pmatrix} A&B\\
C&D\end{pmatrix},(\la,\mu;\kappa) \right)
\end{equation*}
be an element of $G^J$ and let $g_*=T_*^{-1}g_0T_*.$ Then
 \begin{equation*}
g_*=\left(\begin{pmatrix} P & Q\\ {\overline Q} & {\overline P}
\end{pmatrix},\left( {\frac 12}(\la+i\mu),\,{\frac 12}(\la-i\mu);\,-i{\kappa\over
2}\right)\right),
\end{equation*}
where $P$ and $Q$ are given by (5.25) and (5.26).

\vskip 0.2cm For brevity, we write
\begin{equation*}
(\Omega,Z):=\Phi_*(W,\eta)\quad \emph{and}\quad
(\Omega_*,Z_*):=g_0\cdot (\Omega,Z).
\end{equation*}
Then
\begin{equation*}
\Omega=i(I_n+W)(I_n-W)^{-1}\quad \emph{and}\quad Z=2\,i\,\eta
(I_n-W)^{-1}.
\end{equation*}
On the other hand, we set
\begin{equation*}
(W_*,\eta_*):=g_*\cdot(W,\eta)\quad \emph{and}\quad
({\widehat\Omega},{\widehat Z}):=\Phi_*(W_*,\eta_*).
\end{equation*}
Then
\begin{equation*}
W_*=(PW+Q)(\OQ W+\OP)^{-1}  \quad \emph{and}\quad
\eta_*=(\eta+\la_* W+\mu_*)(\OQ W+\OP)^{-1},
\end{equation*}
where $\la_*={\frac 12}(\la+\,i\,\mu)$ and $\mu_*={\frac
12}(\la-\,i\,\mu).$ We can prove that
$$\Omega_*={\widehat\Omega}=\left\{ (i\,A-B)W+(i\,A+B)\right\}\left\{
(i\,C-D)W+(i\,C+D)\right\}^{-1}$$ and
$$Z_*={\widehat Z}=\left\{ 2\,i\,\eta+(\la\, i-\mu)W+\la\, i+\mu\right\}\left\{
(i\,C-D)W+(i\,C+D)\right\}^{-1}.$$ The complete proof can be found
in [52].\hfill $\square$

\vskip 0.3cm For a coordinate $(W,\eta)\in\Dnm$ with
$W=(w_{\mu\nu})\in {\mathbb D}_n$ and $\eta=(\eta_{kl})\in \Cmn,$
we put
\begin{eqnarray*}
dW\,&=&\,(dw_{\mu\nu}),\quad\ \ d{\overline W}\,=\,(d{\overline w}_{\mu\nu}),\\
d\eta\,&=&\,(d\eta_{kl}),\quad\ \
d{\overline\eta}\,=\,(d{\overline\eta}_{kl})
\end{eqnarray*}
and
\begin{eqnarray*}
\PW\,=\,\left(\, { {1+\delta_{\mu\nu}} \over 2}\, {
{\partial}\over {\partial w_{\mu\nu}} } \,\right),\quad
\PWB\,=\,\left(\, { {1+\delta_{\mu\nu}}\over 2} \, {
{\partial}\over {\partial {\overline w}_{\mu\nu} }  } \,\right),
\end{eqnarray*}

$$\PE=\begin{pmatrix} {\partial}\over{\partial \eta_{11}} & \hdots &
 {\partial}\over{\partial \eta_{m1}} \\
\vdots&\ddots&\vdots\\
 {\partial}\over{\partial \eta_{1n}} &\hdots & {\partial}\over
{\partial \eta_{mn}} \end{pmatrix},\quad \PEB=\begin{pmatrix}
{\partial}\over{\partial {\overline \eta}_{11} }   &
\hdots&{ {\partial}\over{\partial {\overline \eta}_{m1} }  }\\
\vdots&\ddots&\vdots\\
{ {\partial}\over{\partial{\overline \eta}_{1n} }  }&\hdots &
 {\partial}\over{\partial{\overline \eta}_{mn} }  \end{pmatrix}.$$

Recently the author [53] obtained the following results.
\newcommand\ot{\overline\eta}
\vskip 0.5cm
\begin{theorem}
The following metric $ds^2$ defined by

\vskip 0.2cm
\begin{eqnarray*}
{\frac 14}\,ds^2&=&\,\textrm{tr}\left( (I_n-W\OW)^{-1}dW(I_n-\OW W)^{-1}\bw\right)\,+\,\textrm{tr}\left( (I_n-W\OW)^{-1}\,{}^t(d\eta)\,\be\,\right)\\
& &\ +\,\textrm{tr}\left(  (\eta\OW-{\overline\eta})(I_n-W\OW)^{-1}dW(I_n-\OW W)^{-1}\,{}^t(d\ot)\right)\\
& &\  +\,\textrm{tr}\left( (\ot W-\eta)(I_n-\OW
W)^{-1}d\OW(I_n-W\OW)^{-1}\,{}^t(d\eta)\,\right)    \\
& &\ -\,\textrm{tr}\left( (I_n-W\OW)^{-1}\,{}^t\eta\,\eta\,
(I_n-\OW W)^{-1}\OW
dW (I_n-\OW W)^{-1}d\OW \, \right)\\
& &\ -\,\textrm{tr}\left( W(I_n-\OW W)^{-1}\,{}^t\ot\,\ot\,
(I_n-W\OW )^{-1}
dW (I_n-\OW W)^{-1}d\OW \,\right)\\
& &\ +\,\textrm{tr}\left( (I_n-W\OW)^{-1}{}^t\eta\,\ot \,(I_n-W\OW)^{-1} dW (I_n-\OW W)^{-1} d\OW\,\right)\\
& &\ +\,\textrm{tr}\left( (I_n-\OW)^{-1}\,{}^t\ot\,\eta\,\OW\,(I_n-W\OW)^{-1} dW (I_n-\OW W)^{-1} d\OW\,\right)\\
& &\ +\,\textrm{tr} ( (I_n-\OW)^{-1}(I_n-W)(I_n-\OW
W)^{-1}\,{}^t\ot\,\eta\,(I_n-\OW W)^{-1}\\
& &\ \ \ \ \ \times\, (I_n-\OW)(I_n-W)^{-1}dW
(I_n-\OW W)^{-1}d\OW\,)\\
& &\ -\,\textrm{tr}
((I_n-W\OW)^{-1}(I_n-W)(I_n-\OW)^{-1}\,{}^t\ot\,\eta\,(I_n-W)^{-1}\\
& & \ \ \ \ \ \times\,dW (I_n-\OW W)^{-1}d\OW\,)
\end{eqnarray*}
\vskip 0.21cm\noindent is a Riemannian metric on $\Dnm$ which is
invariant under the action (5.27) of the Jacobi group $G^J_*$.
\end{theorem}

\newcommand\OZ{\overline{Z}} \vskip 0.2cm \noindent $
\textit{Sketch of Proof.}$ The main ingredients for the proof are
Theorem 5.1 and the partial Cayley transform $\Phi_*.$ \vskip
0.1cm For $(W,\eta)\in \Dnm,$ we write
\begin{equation*}
(\Omega,Z):=\Phi_*(W,\eta).
\end{equation*}
Thus
\begin{equation}
\Om=i(I_n+W)(I_n-W)^{-1},\qquad Z=2\,i\,\eta\,(I_n-W)^{-1}.
\end{equation}

Since
\begin{equation*}
d(I_n-W)^{-1}=(I_n-W)^{-1}dW\,(I_n-W)^{-1}
\end{equation*}
and
\begin{equation*}
I_n+(I_n+W)(I_n-W)^{-1}=2\,(I_n-W)^{-1},
\end{equation*}
we get the following formulas from (5.29)

\begin{eqnarray}
Y&=&{1 \over {2\,i}}\,(\Om-{\overline
\Om}\,)=(I_n-W)^{-1}(I_n-W\OW\,)(I_n-\OW\,)^{-1},\\
V&=&{1 \over {2\,i}}\,(Z-\OZ\,)=\eta\, (I_n-W)^{-1}+\ot\,
(I_n-\OW\,)^{-1},\\
d\Om &=& 2\,i\,(I_n-W)^{-1}dW\,(I_n-W)^{-1},\\
dZ&=&2\,i\,\left\{
d\eta+\eta\,(I_n-W)^{-1}dW\,\right\}(I_n-W)^{-1}.
\end{eqnarray}

By Theorem 5.1 and the formulas (5.30)-(5.33), we obtain the
desired result. We refer to [53] for a complete proof. \hfill
$\square$

\vskip 0.2cm We note that if $n=m=1,$ we get
\begin{eqnarray*}
{\frac 14}\,ds^2&=& { {dW\,d\OW}\over
{(1-|W|^2)^2}}\,+\,{ 1 \over {(1-|W|^2)} }\,d\eta\,d\ot\\
& &+{ {(1+|W|^2)|\eta|^2-\OW \eta^2-W\ot^2}\over {(1-|W|^2)^3} }\,dW\,d\OW\\
& & + { {\eta\OW -\ot}\over {(1-|W|^2)^2} }\,dWd\ot\,+\,{ {\ot W
-\eta}\over {(1-|W|^2)^2} }\,d\OW d\eta.
\end{eqnarray*}

\vskip 0.3cm
\begin{theorem}
The Laplacian $\Delta$ of $(\Dnm,ds^2)$ is given by
\begin{eqnarray*}
\Delta\,&=&\ \  \textrm{tr} \left( (I_n-W\OW)\,{}^t\!\left(
(I_n-W\OW)\PWB\right)\PW\right)\,\\
& &
+\,\textrm{tr}\left( (I_n-\OW W) \PE\,{}^t\!\left( \PEB\right)\right)\\
& & +\,\textrm{tr} \left(\,{}^t(\eta-\ot\,W)\,{}^t\!\left(
\PEB\right)
(I_n-\OW W)\PW  \right)\,\\
& & +\,\textrm{tr}\left( (\ot-\eta\,\OW)\,{}^t\!\left(
(I_n-W\OW)\PWB\right)\PE\right)\\
& &-\,\textrm{tr}\left( \eta \OW
(I_n-W\OW)^{-1}\,{}^t\eta\,\,{}^t\!\left(\PEB\right)(I_n-\OW
W)\PE\right)\\
& &-\,\textrm{tr}\left( \ot W (I_n-\OW W)^{-1}
\,{}^t\ot\,\,{}^t\!\left(\PEB\right)(I_n-\OW
W)\PE\right)\\
& &+\,\textrm{tr}\left( \ot (I_n-W\OW)^{-1}{}^t\eta\,{}^t\!\left(
\PEB\right)
(I_n-\OW W)\PE \right)\\
& & +\,\textrm{tr}\left( \eta\,\OW W (I_n-\OW
W)^{-1}\,{}^t\ot\,\,{}^t\!\left( \PEB\right) (I_n-\OW W)\PE
\right).
\end{eqnarray*}
\end{theorem}

\vskip 0.2cm \noindent $ \textit{Sketch of Proof.}$ For
$(W,\eta)\in \Dnm,$ we write
\begin{equation*}
(\Omega,Z):=\Phi_*(W,\eta).
\end{equation*}
Thus
\begin{equation*}
\Om=i(I_n+W)(I_n-W)^{-1},\qquad Z=2\,i\,\eta\,(I_n-W)^{-1}.
\end{equation*}
From the formulas (5.29),\,(5.32) and (5.33), we get

\begin{equation}
\PO=\,{1\over {2\,i}}\,(I_n-W)\left[ \,{}^t\left\{
(I_n-W)\,\PW\right\} -\,{}^t\left\{ \,{}^t\eta\,\,{}^t\left(
\PE\right)\right\}\,\right]
\end{equation}
and
\begin{equation}
\PZ=\,{1\over {2\,i}}\,(I_n-W)\PE.
\end{equation}
By Theorem 5.1 and the formulas (5.34)-(5.35), we obtain the
desired result. We refer to [53] for a complete proof. \hfill
$\square$

\newcommand\ddww{{{\partial^2}\over{\partial W\partial \OW}}}
\newcommand\ddtt{{{\partial^2}\over{\partial \eta\partial \ot}}}
\newcommand\ddwe{{{\partial^2}\over{\partial W\partial \ot}}}
\newcommand\ddew{{{\partial^2}\over{\partial \OW\partial \eta}}}

\vskip 0.32cm We note that if $n=m=1,$ we get
\begin{eqnarray*}
\Delta&=&\ \ (1-|W|^2)^2 \ddww\,+\,(1-|W|^2)\,\ddtt\\
& &+\,(1-|W|^2)(\eta-\ot\,W)\,\ddwe\,+\,(1-|W|^2)(\ot-\eta\,\OW)\,\ddew\\
& &-(\OW\,\eta^2+W\ot^2)\,\ddtt\,+\,(1+|W|^2)|\eta|^2\,\ddtt.
\end{eqnarray*}

\vskip 0.2cm \noindent {\bf Remark 5.3.} We see from Remark 5.2
that the following differential operators
\begin{equation*}
{\mathbb S}_1:=\,\textrm{tr}\left( (I_n-\OW
W)\PE\,^t\!\left(\PEB\right)\right)\quad \textrm{and}\quad
{\mathbb S}_2:=\Delta-{\mathbb S}_1
\end{equation*}
are invariant under the action (5.27) of $G_*^J.$

\vskip 0.5cm We now describe the algebra $\BD(\BH_{n,m})$ of all
differential operators on $\BH_{n,m}$ invariant under the action
(5.15) of $G^J$. The adjoint action of $K$ on the tangent space
$\fp^J$ of $\Hnm$ at $(iI_n,0)$ induces the action of $K$ on the
complex vector space $T_n\times \BC^{(m,n)}$ defined by
\begin{equation}
h\cdot (\omega,z)=(h\,\om \,^th,\,z\,^th),
\end{equation}

\noindent where $h\in K,\ w\in T_n$ and $z\in \Cmn$. Here we
regard the complex vector space $T_n\times \BC^{(m,n)}$ as a real
vector space. The action (5.36) induces the action $\rho$ of $K$
on the polynomial algebra
$\text{Pol}_{m,n}=\,\text{Pol}\,(T_n\times\Cmn).$ We denote by
$\text{Pol}_{m,n}^K$ the subalgebra of $\text{Pol}_{m,n}$
consisting of all $K$-invariants of the action $\rho$ of $K.$ For
brevity, we put
$$T_{n,m}:=T_n \times \BC^{(m,n)}.$$
The following $K$-invariant inner product $(\,\,,\,)_*$ of the
complex vector space $T_{n,m}$ defined by
\begin{equation*}
((\om,z),(\om',z'))_*:= \textrm{tr}(\om{\overline {\om'}}\,)+
\textrm{tr}(z\,^t{\overline {z'}}\,),\quad (\om,z),\,(\om',z')\in
T_{n,m}
\end{equation*}

\noindent gives a canonical isomorphism
\begin{equation*}
T_{n,m}\cong\,T_{n,m}^*,\quad (\om,z)\mapsto f_{\om,z},\quad
(\om,z)\in T_{n,m},
\end{equation*}

\noindent where $f_{\om,z}$ is the linear functional on $T_{n,m}$
defined by
\begin{equation*}
f_{\om,z}((\om',z'\,)):=((\om,z),(\om',z'))_*,\quad (\om',z'\,)\in
T_{n,m}.
\end{equation*}

\noindent Let $e_i\,(1\leq i\leq n),\ e_{ij}\ (1\leq i < j\leq n)$
as before. We let $f_{kl}\ (1\leq k\leq m,\ 1\leq l\leq n)$ be the
$m\times n$ matrix where the $k$-th row and the $l$-th column
meet, and all other entries $0$. Then we see that $e_i\,(1\leq
i\leq n),\ e_{ij}\ (1\leq i < j\leq n),\ f_{kl}\ (1\leq k\leq m,\
1\leq l\leq n)$ form an orthonormal basis for $T_{n,m}$ with
respect to the inner product $(\,\,,\,)_*.$ Once and for all we
choose the above orthonormal basis for $T_{n,m}$. Then one gets a
canonical linear bijection of $S(T_{n,m})^K$ onto
$\BD(\BH_{n,m})$. Identifying $T_{n,m}$ with $T_{n,m}^*$ by the
above isomorphism, one gets a natural linear bijection
$$\Psi_*:\,\text{Pol}^K_{m,n}\lrt \BD(\BH_{n,m})$$
of $\text{Pol}^{K}_{m,n}$ onto $\BD(\BH_{n,m}).$ We refer to [29],
p. 287. It may be shown that the algebra $\BD(\BH_{n,m})$ is
generated by the images under the mapping $\Psi_*$ of the
following invariants in $\text{Pol}_{m,n}^K$\,: \par\vskip 0.2cm \
\ \ \ \ \ \ $\text{(I1)} \ p_j(\om,z)=\,\text{tr}((\om{\overline
\om})^j),$\par\ \ \ \ \ \ \ $\text{(I2)}\
\psi_k^{(1)}(\om,z)=\,(z\,^t{\overline z})_{kk},$\par\ \ \ \ \ \ \
$ \text{(I3)}\ \psi_{kp}^{(2)}(\om,z)=
\,\text{Re}\,(z\,^t{\overline z})_{kp}, $\par \ \ \ \ \ \ \
$\text{(I4)}\ \psi_{kp}^{(3)}(\om,z)
=\,\text{Im}\,(z\,^t{\overline z})_{kp}, $\par \ \ \ \ \ \ \
$\text{(I5)}\ f_{kp}^{(1)}(\om,z)= \,\text{Re}\,(z{\overline
\om}\,^tz)_{kp}, $\par \ \ \ \ \ \ \ $\text{(I6)}\
f_{kp}^{(2)}(z,w)= \,\text{Im}\,(z{\overline \om}\,^tz\,)_{kp},
$\par\vskip 0.2cm\noindent where $1\leq j\leq n,\ 1\leq k\leq m$
and $1\leq k\leq p\leq m$. In the case $n=m=1,$ the algebra
$\BD(\BH_{1,1})$ can be described explicitly as follows.

\def\ddx{{{\partial^2}\over{\partial x^2}}}
\def\ddy{{{\partial^2}\over{\partial y^2}}}
\def\ddu{{{\partial^2}\over{\partial u^2}}}
\def\ddv{{{\partial^2}\over{\partial v^2}}}
\def\px{{{\partial}\over{\partial x}}}
\def\py{{{\partial}\over{\partial y}}}
\def\pu{{{\partial}\over{\partial u}}}
\def\pv{{{\partial}\over{\partial v}}}
\def\pxu{{{\partial^2}\over{\partial x\partial u}}}
\def\pyv{{{\partial^2}\over{\partial y\partial v}}}
\def\DSPR{{\Bbb D}(\SPR)}
\def\dx{{{\partial}\over{\partial x}}}
\def\dy{{{\partial}\over{\partial y}}}
\def\du{{{\partial}\over{\partial u}}}
\def\dv{{{\partial}\over{\partial v}}}
\vskip 0.2cm
\begin{theorem}
The algebra $\BD({\mathbf H}_{1,1})$ is generated by the following
differential operators
\begin{align*} D=&y^2\,\left(\,{{\partial^2}\over {\partial x^2}}+
{{\partial^2}\over {\partial y^2}}\,\right)
+v^2\,\left(\,\ddu\,+\,\ddv\,\right) \\
&\ \ +2\,y\,v\,\left(\,\pxu\,+\,\pyv\,\right),
\end{align*}
$$\Psi=y\left(\,{{\partial^2}\over {\partial u^2}}+
{{\partial^2}\over {\partial v^2}}\,\right),\hskip 3.74cm$$
\begin{align*}D_1=&\,2y^2\,{{\partial^3}\over {\partial x\partial u
\partial v}}-y^2\,{{\partial}\over{\partial y}}
\left(\,{{\partial^2}\over{\partial u^2}}-
{{\partial^2}\over{\partial v^2}}\,\right)\hskip 1cm\\ &\ \ \
+\left(\, v\,{{\partial}\over{\partial v}}\,+\,1\,\right)\Psi
\end{align*} and
\begin{align*} D_2=&\,y^2\,{{\partial}\over{\partial x}}\left(\,
{{\partial^2}\over{\partial v^2}}\,-\,{{\partial^2}\over {\partial
u^2}}\,\right)\,-\,2\,y^2\,{{\partial^3}\over{\partial y\partial u
\partial v}}\\ &\ \ \ \ -\,v\,{{\partial}\over{\partial u}}\Psi,
\end{align*} where $\tau=x+iy$ and $z=u+iv$ with real variables
$x,y,u,v.$ Moreover, we have \begin{align*} D\Psi-&\Psi D\,=\,
2\,y^2\,\dy\left(\,\ddu\,-\,\ddv\,\right)\\ & -
4\,y^2\,{{\partial^3}\over{\partial x\partial u\partial
v}}-2\,\left(\,v\,\dv\Psi+\Psi\,\right).
\end{align*}\par
\noindent In particular, the algebra $\BD({\mathbf H}_1\times
\BC)$ is not commutative.
\end{theorem}
\vskip 0.2cm \noindent $ \textit{Proof.}$ The proof can be found
in [18,\,50].\hfill $\square$

\def\bz{d{\overline Z}}
\def\bo{d{\overline \O}}

\vskip 0.53cm \noindent $ \textbf{5.3.\ A Fundamental Domain for
the Siegel-Jacobi Space}$

\vskip 0.53cm Before we describe a fundamental domain for the
Siegel-Jacobi space, we review the Siegel's fundamental domain for
the Siegel upper half plane. \vskip 0.2cm We let
\begin{equation*}
{\mathcal P}_n=\left\{ Y\in\BR^{(n,n)}\,|\ Y=\,^tY>0\,\right\}
\end{equation*}

\noindent be an open cone in $\BR^{n(n+1)/2}$. The general linear
group $GL(n,\BR)$ acts on ${\mathcal P}_n$ transitively by
\begin{equation*}
g\circ Y:=gY\,^tg,\quad g\in GL(n,\BR),\ Y\in {\mathcal P}_n.
\end{equation*}

\noindent Thus ${\mathcal P}_n$ is a symmetric space diffeomorphic
to $GL(n,\BR)/O(n).$

\newcommand\Mg{{\mathcal M}_n}
\newcommand\Rg{{\mathcal R}_n}

The fundamental domain $\Rg$ for $GL(n,\BZ)\ba {\mathcal P}_n$
which was found by H. Minkowski\,[39] is defined as a subset of
${\mathcal P}_n$ consisting of $Y=(y_{ij})\in {\mathcal P}_n$
satisfying the following conditions (M.1)-(M.2)\ (cf.\,[35]
p.\,123): \vskip 0.1cm (M.1)\ \ \ $aY\,^ta\geq y_{kk}$\ \ for
every $a=(a_i)\in\BZ^n$ in which $a_k,\cdots,a_n$ are relatively
prime for $k=1,2,\cdots,n$. \vskip 0.1cm (M.2)\ \ \ \
$y_{k,k+1}\geq 0$ \ for $k=1,\cdots,n-1.$ \vskip 0.1cm We say that
a point of $\Rg$ is {\it Minkowski reduced} or simply {\it M}-{\it
reduced}.

\vskip 0.1cm Siegel\,[43] determined a fundamental domain
${\mathcal F}_n$ for $\G_n\ba \BH_n,$ where $\G_n=Sp(n,\BZ)$ is
the Siegel modular group of degree $n$. We say that $\Om=X+iY\in
\BH_n$ with $X,\,Y$ real is {\it Siegel reduced} or {\it S}-{\it
reduced} if it has the following three properties: \vskip 0.1cm
(S.1)\ \ \ $\det (\text{Im}\,(\g\cdot\Om))\leq \det
(\text{Im}\,(\Om))\qquad\text{for\ all}\ \g\in\G_n$; \vskip 0.1cm
(S.2)\ \ $Y=\text{Im}\,\Om$ is M-reduced, that is, $Y\in \Rg\,;$
\vskip 0.1cm (S.3) \ \ $|x_{ij}|\leq {\frac 12}\quad \text{for}\
1\leq i,j\leq n,\ \text{where}\ X=(x_{ij}).$ \vskip 0.1cm
${\mathcal F}_n$ is defined as the set of all Siegel reduced
points in $\BH_n.$ Using the highest point method, Siegel proved
the following (F1)-(F3)\,(cf. [35] p.\,169): \vskip 0.1cm (F1)\ \
\ $\G_n\cdot {\mathcal F}_n=\BH_n,$ i.e.,
$\BH_n=\cup_{\g\in\G_n}\g\cdot {\mathcal F}_n.$ \vskip 0.1cm (F2)\
\ ${\mathcal F}_n$ is closed in $\BH_n.$ \vskip 0.1cm (F3)\ \
${\mathcal F}_n$ is connected and the boundary of ${\mathcal F}_n$
consists of a finite number of hyperplanes.

The metric $ds_*^2$ given by (5.12) induces a metric $ds_{\mathcal
F}^2$ on ${\mathcal F}_n$. Siegel\,[43] computed the volume of
${\mathcal F}_n$
\begin{equation}
\text{vol}\,(\CCF)=2\prod_{k=1}^n\pi^{-k}\G
(k)\zeta(2k),\end{equation} where $\G (s)$ denotes the Gamma
function and $\zeta (s)$ denotes the Riemann zeta function. For
instance,
$$\text{vol}\,({\mathcal F}_1)={{\pi}\over 3},\quad \text{vol}\,({\mathcal F}_2)={{\pi^3}\over {270}},
\quad \text{vol}\,({\mathcal F}_3)={{\pi^6}\over {127575}},\quad
\text{vol}\,({\Cal F}_4)={{\pi^{10}}\over {200930625}}.$$

Let $f_{kl}\,(1\leq k\leq m,\ 1\leq l\leq n)$ be the $m\times n$
matrix with entry $1$ where the $k$-th row and the $l$-th column
meet, and all other entries $0$. For an element $\Om\in \BH_n$, we
set for brevity
\begin{equation*}
h_{kl}(\Om):=f_{kl}\Om,\qquad 1\leq k\leq m,\ 1\leq l\leq
n.\end{equation*}
 \indent For each $\Om\in {\mathcal F}_n,$ we define a
subset $P_{\Om}$ of $\BC^{(m,n)}$ by
\begin{equation*}
P_{\Om}=\left\{ \,\sum_{k=1}^m\sum_{j=1}^n \la_{kl}f_{kl}+
\sum_{k=1}^m\sum_{j=1}^n \mu_{kl}F_{kl}(\Om)\,\Big|\ 0\leq
\la_{kl},\mu_{kl}\leq 1\,\right\}. \end{equation*} \indent For
each $\Om\in {\mathcal F}_n,$ we define the subset $D_{\Om}$ of
$\BH_n\times \BC^{(m,n)}$ by
\begin{equation*} D_{\Om}:=\left\{\,(\Om,Z)\in\BH_n\times \BC^{(m,n)}\,\vert\ Z\in
P_{\Om}\,\right\}.\end{equation*}
\newcommand\Fgh{{\mathcal F}_{n,m}}
We define
\begin{equation*} \Fgh:=\cup_{\Om\in {\mathcal F}_n}D_{\Omega}.\end{equation*}

\vskip 0.3cm
\begin{theorem}
Let
$$\Gamma_{n,m}:=Sp(n,{\mathbb Z})\ltimes H_{\mathbb Z}^{(n,m)}$$
be the discrete subgroup of $G^J$, where
$$H_{\BZ}^{(n,m)}=\left\{ (\lambda,\mu;\kappa)\in
H_{\BR}^{(n,m)}\,|\ \lambda,\mu,\kappa \ \textrm{are integral}\
\right\}.$$ Then $\Fgh$ is a fundamental domain for $\G_{n,m}\ba
\BH_{n,m}.$
\end{theorem}
\vskip 0.2cm \noindent $\textit{Proof.}$ The proof can be found in
[54]. \hfill $\square$

\vskip 0.53cm \noindent $ \textbf{5.4.\ Maass-Jacobi Forms}$

\vskip 0.3cm In the case $n=m=1$, R. Berndt [18] introduced the
notion of Maass-Jacobi forms. Now we generalize this notion to the
general case.

\vskip 0.2cm
\begin{definition}
Let
$$\Gamma_{n,m}:=Sp(n,{\mathbb Z})\ltimes H_{\mathbb Z}^{(n,m)}$$
be the discrete subgroup of $G^J$, where
$$H_{\BZ}^{(n,m)}=\left\{ (\lambda,\mu;\kappa)\in
H_{\BR}^{(n,m)}\,|\ \lambda,\mu,\kappa \ \textrm{are integral}\
\right\}.$$ A smooth function $f:\Hnm\lrt \BC$ is called a
$\textsf{Maass}$-$\textsf{Jacobi form}$ on $\Hnm$ if $f$ satisfies
the following conditions (MJ1)-(MJ3)\,:\vskip 0.1cm (MJ1)\ \ \ $f$
is invariant under $\G_{n,m}.$\par (MJ2)\ \ \ $f$ is an
eigenfunction of the Laplace-Beltrami operator $\Delta_{n,m}$.\par
(MJ3)\ \ \ $f$ has a polynomial growth, that is, there exist a
constant $C>0$ and a positive \par \ \ \ \ \ \ \ \ \ \ \ integer
$N$ such that
\begin{equation*}
|f(X+iY,Z)|\leq C\,|p(Y)|^N\quad \textrm{as}\ \det
Y\longrightarrow \infty,
\end{equation*}

\ \ \ \ \ \ \ \ \ \ \ where $p(Y)$ is a polynomial in
$Y=(y_{ij}).$
\end{definition}

\vskip 0.3cm It is natural to propose the following problems.

\vskip 0.3cm\noindent {\bf {Problem\ A}\,:} Construct Maass-Jacobi
forms.

\vskip 0.3cm\noindent {\bf {Problem\ B}\,:} Find all the
eigenfunctions of $\Delta_{n,m}.$

\vskip 0.3cm  We consider the simple case $n=m=1.$ A metric
$ds_{1,1}^2$ on ${\mathbf H}_1\times \BC$ given by
\begin{align*} ds^2_{1,1}\,=\,&{{y\,+\,v^2}\over
{y^3}}\,(\,dx^2\,+\,dy^2\,)\,+\, {\frac 1y}\,(\,du^2\,+\,dv^2\,)\\
&\ \ -\,{{2v}\over {y^2}}\, (\,dx\,du\,+\,dy\,dv\,)\end{align*} is
a $G^J$-invariant K{\"a}hler metric on ${\mathbf H}_1\times \BC$.
Its Laplacian $\Delta_{1,1}$ is given by
\begin{align*} \Delta_{1,1}\,=\,& y^2\,\left(\,\ddx\,+\,\ddy\,\right)\,\\ &+\,
(\,y\,+\,v^2\,)\,\left(\,\ddu\,+\,\ddv\,\right)\\ &\ \
+\,2\,y\,v\,\left(\,\pxu\,+\,\pyv\,\right).\end{align*}

\vskip 0.2cm We provide some examples of eigenfunctions of
$\Delta_{1,1}$.
\par
(1) $h(x,y)=y^{1\over 2}K_{s-{\frac12}}(2\pi |a|y)\,e^{2\pi iax} \
(s\in \BC,$ $a\not=0\,)$ with eigenvalue $s(s-1).$ Here
$$K_s(z):={\frac12}\int^{\infty}_0 \exp\left\{-{z\over
2}(t+t^{-1})\right\}\,t^{s-1}\,dt,$$ \indent \ \ \ where
$\mathrm{Re}\,z
> 0.$ \par (2) $y^s,\ y^s x,\ y^s u\ (s\in\BC)$ with eigenvalue
$s(s-1).$
\par
 (3) $y^s v,\ y^s uv,\ y^s xv$ with eigenvalue $s(s+1).$
\par
(4) $x,\,y,\,u,\,v,\,xv,\,uv$ with eigenvalue $0$.
\par
(5) All Maass wave forms.

\vskip 0.5cm
\begin{center}
{\bf 5.5. Formal Eisenstein Series}
\end{center}
\vskip 0.3cm Let
$$\Gamma_{1,1}^{\infty}=\left\{ \left( \begin{pmatrix} \pm 1 & m\\
0 & \pm 1\end{pmatrix},(0,n,\kappa)\right)\,\Big|\ m,n,\kappa\in
{\mathbb Z}\,\right\}$$ be the subgroup of
$\Gamma_{1,1}=SL_2(\mathbb Z)\ltimes H_{\mathbb Z}^{(1,1)}.$
For $\gamma=\left( \begin{pmatrix} a & b\\
c & d\end{pmatrix},(\lambda,\mu,\kappa)\right)\in \Gamma_{1,1},$
we put \\ $(\tau_{\gamma},z_{\gamma})=\gamma\cdot (\tau,z).$ That
is,
\begin{align*}
\tau_{\gamma}&=(a\tau+b)(c\tau+d)^{-1},\\ z_{\gamma}&=(z+\lambda
\tau+\nu)(c\tau+d)^{-1}.\end{align*} We note that if
$\gamma\in\Gamma_{1,1},$
$$\mathrm{Im}\,\tau_{\gamma}=\mathrm{Im}\,\tau,\ \
\mathrm{Im}\,z_{\gamma}=\mathrm{Im}\,z$$ if and only if
$\gamma\in\Gamma_{1,1}^{\infty}.$ For $s\in\BC,$ we define an
Eisenstein series formally by
$$E_s(\tau,z)=\sum_{\gamma\in \Gamma_{1,1}^{\infty}\backslash
\Gamma_{1,1}}(\mathrm{Im}\,\tau_{\gamma})^s\cdot
\mathrm{Im}\,z_{\gamma}.$$ Then $E_s(\tau,z)$ satisfies formally
$$E_s(\gamma\cdot(\tau,z))=E_s(\tau,z),\ \ \gamma\in
\Gamma_{1,1}$$ and $$\Delta E_s(\tau,z)=s(s+1)E_s(\tau,z).$$ But
the series does not converge.

\vskip 0.5cm
\begin{center}
{\bf 5.6. Remarks on Fourier Expansions of Maass-Jacobi Forms}
\end{center}
\vskip 0.3cm
\newcommand{\bHC}{{\mathbf H}_1\times \BC}
We let $f:{\mathbf H}_1\times \BC\lrt \BC$ be a Maass-Jacobi form
with $\Delta f=\lambda f.$ Then $f$ satisfies the following
invariance relations
$$f(\tau+n,\,z)\,=\,f(\tau,z)\ \ \ \mathrm{for\ all}\ n\in
{\mathbb Z}$$ and
$$f(\tau,\,z\,+\,n_1\tau\,+\,n_2)\,=\,f(\tau,z)$$
$\mathrm{for\ all}\ n_1,\,n_2\in {\mathbb Z}.$ Therefore $f$ is a
smooth function on $\bHC$ which is periodic in $x$ and $u$ with
period $1.$ So $f$ has the following Fourier series
$$f(\tau,z)\,=\,\sum_{n\in \mathbb Z}\sum_{r\in \mathbb Z}\,c_{n,r}(y,v)\,
e^{2\pi i(nx+ru)}.$$ For two fixed integers $n$ and $r$, we have
to calculate the function $c_{n,r}(y,v).$ For brevity, we put
$F(y,v)=\,c_{n,r}(y,v).$ Then $F$ satisfies the following
differential equation
\begin{equation}
\left[y^2\,\ddy\,+\,(y+v^2)\,\ddv\,+\,2yv\, \pyv\,\right]\,F =
\left\{\,(ay+bv)^2\,+\,b^2y\,+\,\lambda\,\right\}
\,F.\end{equation}
Here $a=2\pi n$ and $b=2\pi r$ are constant. We
note that the function $u(y)=y^{\frac 12} K_{s-{\frac 12}}(2\pi
\vert n\vert y)$ satisfies the above differential equation with
$\lambda=s(s-1).$ Here $K_s(z)$ is the $K$-Bessel function before.

\vskip 0.2cm\noindent {\bf {Problem\ C}\,:} Find the solutions of
the above differential equations (5.38) explicitly.

\vskip 0.5cm
\begin{center}
{\bf 5.7. Jacobi Forms}
\end{center}

\vskip 0.3cm
Let $\rho$ be an irreducible representation of $K$ on a finite
dimensional vector space $V_{\rho}$ with highest weight
$\ell=(\ell_1,\cdots,\ell_n).$ Then $\rho$ is extended to a
rational representation of $GL(n,\BC)$ denoted also by $\rho.$ Let
${\mathcal M}$ be a symmetric half-integral semi-positive definite
matrix of degree $m$.
The $ \textsf{canonical automorphic factor}$ \\
$$J_{\rho,\CM}:G^J\times \BH_{n,m}\lrt GL(V_{\rho})$$
is given by
\begin{align*}
\hskip 1cm J_{\rho,\CM}(g,(\Omega,Z)) &=& e^{-2\pi i\,\textrm{tr}
(\CM [Z +\la \Omega+\mu](C\Omega+D)^{-1}C)}\cdot e^{2\pi
i\,\textrm{tr} (\CM(\la \Omega\,{}^t\!\la+2\la\,
{}^tZ+\kappa+\mu\,{}^t\!\la))}\hskip 4cm\\ & &\ \times
\rho(C\Omega+D)^{-1},\hskip 11cm
\end{align*}

\noindent where $g=(M,(\la,\mu,\kappa))\in G^J$ with
$M=\begin{pmatrix} A & B\\ C & D\end{pmatrix}\in Sp(n,\BR).$

\vskip 0.1cm For a function $f\in C^{\infty}(\BH_{n,m},V_{\rho}),$
we define
\begin{equation*}
(f|_{\rho,\CM}g)(\Omega,Z):=J_{\rho,\CM}(g,(\Omega,Z))f(g\cdot
(\Omega,Z)).
\end{equation*}

\indent Let $\G_*$ be an arithmetic subgroup of the Siegel modular
group $\G_n=Sp(n,\BZ)$ and let
$$\G_*^J:=\G_*\ltimes H_{\BZ}^{(n,m)}$$
be the discrete subgroup of $\G_{n,m}.$ A $ \textsf{Jacobi form}$
of index $\CM$ with respect to $\G_*$ is defined to be a
$V_{\rho}$-holomorphic function on $\BH_{n,m}$ satisfying the
conditions (A) and (B): \vskip 0.1cm (A) $f|_{\rho,\CM}\g=f$ for
all $\g\in \G_*^J$. \vskip 0.2cm (B) $f$ has a Fourier expansion
of the following form
\begin{equation*}
f(\Omega,Z)=\sum_{T\geq 0}\sum_{R\in\BZ^{(n,m)}}c(T,R)\cdot
e^{2\pi i\,\textrm{tr}(T\Omega+RZ)}
\end{equation*}

\noindent where $c(T,R)\neq 0$ only if $\begin{pmatrix} T & {\frac 12}R\\
{\frac 12}\,{}^t\!R & \CM\end{pmatrix}\geq 0$ and $T$ runs over
the set of all half-integral matrices of degree $n$. We denote by
$J_{\rho,\CM}(\G_*)$ the space of all Jacobi forms of index $\CM$
with respect to $\G_*$.

\vskip 0.1cm
\begin{definition}
A Jacobi form $f\in
J_{\rho,\CM}(\G_*)$ is said to be a $\textsf{cuspidal}$ form if
$\begin{pmatrix} T & {\frac 12}R\\ {\frac 12}\,{}^t\!R &
\CM\end{pmatrix}> 0$ for any $T,R$ with $c(T,R)\neq 0$. A Jacobi
form $f\in J_{\rho,\CM}(\G_*)$ is said to be a $\textsf{singular}$
form if it admits a Fourier expansion such that a Fourier
coefficient $c(T,R)$ vanishes unless
$$\det\begin{pmatrix} T & {\frac 12}R\\ {\frac 12}\,{}^t\!R & \CM\end{pmatrix}= 0.$$
\end{definition}
\indent Under the assumption that $\CM$ is positive definite, the
author [47] introduced the differential operator $M_{n,m,\CM}$
characterizing singular Jacobi forms. We let
\begin{equation*}
{\mathcal P}_n=\left\{ Y\in \BR^{(n,n)}\,|\ Y=\,^tY>0\ \right\}
\end{equation*}
be the open convex cone of positive definite real matrices of
degree $n$. We define the differential operator $M_{n,m,\CM}$ on
${\mathcal P}_n\times \BR^{(m,n)}$ by
\begin{equation}
M_{n,m,\CM}:=\det (Y)\cdot \det \left( {{\partial}\over {\partial
Y}}+{ 1\over {8\pi}}\,\,{}^t\!\left( {{\partial}\over {\partial
V}}\right) \CM^{-1}\left( {{\partial}\over {\partial
V}}\right)\right),
\end{equation}
where
\begin{equation*}
Y=(y_{\mu\nu})\in {\mathcal P}_n,\quad V=(v_{kl})\in
\BR^{(m,n)},\quad {{\partial}\over {\partial Y}}=\left( {
{1+\delta_{\mu\nu}}\over 2} {{\partial}\over {\partial y_{\mu\nu}}
} \right)\end{equation*} and

\begin{equation*}
{{\partial}\over {\partial V}}= \left({{\partial}\over {\partial
v_{kl}} }\right).
\end{equation*}
\vskip 0.3cm The author [47] proved that the following conditions
$(S1)$ and $(S2)$ are equivalent\,: \vskip 0.1cm \ \ \ $(S1)\ f$
is a singular Jacobi form in $J_{\rho,\CM}(\G_n)$. \vskip 0.1cm \
\ \ $(S2)\ f$ satisfies the differential equation
$M_{n,m,\CM}f=0.$ \vskip 0.2cm Let $GL_{n,m}$ be the semidirect
product of the general linear group $GL(n,\BR)$ and the Heisenberg
group $H_{\BR}^{(n,m)}$ equipped with the multiplication
\begin{eqnarray}
(A,(\la,\mu;\ka))\cdot (B,(\la',\mu';\ka'))\hskip 5cm\\
=(AB,(\la B+\la',\mu\,^tB^{-1}+\mu';\,\k+\k'+\la
B\,^t\mu'-\mu\,^tB^{-1}\,^t\la')),\nonumber
\end{eqnarray}
where $A,B\in GL(n,\BR)$ and $(\la,\mu;\k),\,(\la',\mu';\k')\in
H_{\BR}^{(n,m)}.$ The action (5.15) of $G^J$ on $\Hnm$ yields the
natural action of $GL_{n,m}$ on the homogeneous space ${\mathcal
P}_n\times \BR^{(m,n)}$ of non-reductive type defined by
\begin{equation}
(A,(\la,\mu;\ka))\cdot (Y,V):=(AY\,^tA,(V+\la Y+\mu)\,^tA),
\end{equation}
where $A\in GL(n,\BR),\ (\la,\mu;\ka)\in H_\BR^{(n,m)}$ and
$(Y,V)\in {\mathcal P}_n\times \BR^{(m,n)}$. We can show that the
differential operator $M_{n,m,\CM}$ is invariant under the action
of the subgroup
$$GL(n,\BR)\ltimes \left\{ (0,\mu;0)\,|\ \mu\in
\BR^{(m,n)}\,\right\}$$ on ${\mathcal P}_n\times \BR^{(m,n)}$.
Selberg [42] showed that the algebra $\BD({\mathcal P}_n)$ of all
differential operators on ${\mathcal P}_n$ invariant under the
action of $GL(n,\BR)$ on ${\mathcal P}_n$ defined by
\begin{equation*}
A\cdot Y:=AY\,^tA,\quad A\in GL(n,\BR),\ Y\in {\mathcal P}_n
\end{equation*}

\noindent is generated by the following $n$ algebraically
independent commuting differential operators
\begin{equation*}
{\mathbb B}_j=\text{tr}\left( \left(Y {{\partial}\over {\partial
Y}}\right)^j\right),\quad j=1,2,\cdots,n.
\end{equation*}

\noindent We refer to [35,\,45] for more details on ${\mathbb
B}_j$. I present the following problem.

\vskip 0.3cm \noindent {\bf {Problem\ D}\,:} Describe all the
generators of the algebra $\BD({\mathcal P}_n\times \BR^{(m,n)})$
of all differential operators on ${\mathcal P}_n\times
\BR^{(m,n)}$ invariant under the action (5.41) of $GL_{n,m}$ on
${\mathcal P}_n\times \BR^{(m,n)}$.

\vskip 0.2cm I want to mention that the theory of singular Jacobi
forms may be applied to the study of the geometry of the universal
abelian variety as the theory of singular modular forms to the
geometry of the Siegel modular variety\,(cf.\,[49]).

\vskip 0.3cm
\newcommand\CF{{\mathcal F}}
We define the {\it lifting}
\begin{equation*}
\Phi_{\rho,\CM}:\CF(\BH_{n,m},V_{\rho})\lrt \CF(G^J,V_{\rho})
\end{equation*}
by
\begin{equation*}
(\Phi_{\rho,\CM}f)(g):=(f|_{\rho,\CM}g)(iI_n,0)
=J_{\rho,\CM}(g,(iI_n,0))f(g\cdot (iI_n,0)),
\end{equation*}

\noindent where $\CF(\BH_{n,m},V_{\rho})$\,(resp.\,$
\CF(G^J,V_{\rho})$) denotes the space consisting of
$V_{\rho}$-valued functions on $\BH_{n,m}$\,(resp.\,$G^J$). We let
\begin{equation*}
\CF_{\rho,\CM}^{\G}(\BH_{n,m})
\end{equation*}
be the space of all $V_{\rho}$-valued functions $f$ on $\BH_{n,m}$
satisfying the transformation formula
\begin{equation*}
f|_{\rho,\CM}\g=f\quad \textmd{for all}\ \g\in \G^J.
\end{equation*}

We let
\begin{equation*}
\CF_{\rho,\CM}^{\G}(G^J)
\end{equation*}
be the space of all $V_{\rho}$-valued functions $\Phi$ on $G^J$
satisfying the transformation formulas
\begin{equation*}
\Phi(\g g)=\Phi(g)\quad \textmd{for all}\ \g\in \G^J,\ g\in G^J
\end{equation*}
and
\begin{equation*}
\Phi(g\cdot (k,(0,0,\kappa)))=e^{2\pi i\,
\textrm{tr}(\CM\kappa)}\rho(k)^{-1}\Phi(g)
\end{equation*}

\noindent for all $(k,(0,0,\kappa))\in K^J.$ Then we have the
isomorphism
\begin{equation*}
\CF_{\rho,\CM}^{\G}(\BH_{n,m})\cong \CF_{\rho,\CM}^{\G}(G^J)
\end{equation*}

We denote by $E(\rho,\CM)$ the Hilbert space consisting of
$V_{\rho}$-valued measurable functions $\varphi$ on $\BH_{n,m}$
such that
\begin{equation*}
\int_{\BH_{n,m}}(\rho(Y)\varphi(\Omega,Z),\varphi(\Omega,Z))\kappa_{\CM}(\Omega,Z)dv<\infty,
\end{equation*}

\noindent where
$$\kappa_{\CM}(\Omega,Z)=e^{-4\pi\,\textrm{tr}(\,^tV\CM VY^{-1})}.$$
Let $\chi$ be the unitary character of the center ${\mathcal Z}=
\textrm{Symm}(\BR^n)$ defined by
\begin{equation*}
\chi_\CM (\kappa):=e^{2\pi i\,\textrm{tr}(\CM\kappa)},\ \
\kappa\in {\mathcal Z}.
\end{equation*}

\noindent The induced representation $
\textrm{Ind}_{K^J}^{G^J}(\rho\otimes {\overline {\chi}}_\CM)$ is
realized on $E(\rho,\CM)$ by
\begin{equation*}
\left(\textrm{Ind}_{K^J}^{G^J}(\rho\otimes {\overline
{\chi}}_\CM)(g)\varphi\right)(\Omega,Z)
=J_{\rho,\CM}(g^{-1},(\Omega,Z))^{-1}\varphi(g^{-1}\cdot(\Omega,Z)),
\end{equation*}

\noindent where $g\in G^J, \ \varphi\in E(\rho,\CM)$ and $K^J\cong
U(n)\times {\mathcal Z}.$ The subspace
\begin{equation*}
H(\rho,\CM):=\left\{ f\in E(\rho,\CM)\,|\ f\ \textrm{is
holomorphic}\,\right\}
\end{equation*}
is a closed $G^J$-invariant subspace of $E(\rho,\CM)$. Let
$\pi_{\rho,\CM}$ be the restriction of $
\textrm{Ind}_{K^J}^{G^J}(\rho\otimes {\overline {\chi}}_\CM)$ to
$H(\rho,\CM)$. We can show that if $\ell_n>n+{\frac 12}$, then
$H(\rho,\CM)\neq 0$ and $\pi_{\rho,\CM}$ is an irreducible unitary
representation of $G^J$ which is square integrable modulo
${\mathcal Z}.$ Let
$$L^2\left(\G_*^J\backslash G^J; V_{\rho}\right)$$
be the Hilbert space of all $\G^J_*$-invariant $V_{\rho}$-valued
measurable functions on $G^J$ with finite norm. We denote by
$L^2_0\left(\G_*^J\backslash G^J; V_{\rho}\right)$ the subspace
consisting of functions $\varphi\in L^2\left(\G_*^J\backslash G^J;
V_{\rho}\right)$ satisfying the condition
\begin{equation*}
\int_{(N^g\cap \G_*^J)\ba N^g}\varphi(ng_0)dn=0
\end{equation*}
for any cuspidal subgroup $N^g$ of $G^J$ and almost all $g_0\in
G^J.$
\vskip 0.1cm\noindent
\begin{theorem}
\textbf{(Duality Theorem)} Let $\rho$ be an irreducible
representation of $K$ on $V_{\rho}$ with highest weight
$\ell=(\ell_1,\cdots,\ell_n)\in \BZ^n$ with $\ell\geq
\cdots\geq\ell_n.$ Suppose $\ell_n>n+{\frac 12}$. Let $\CM$ be a
half-integral positive definite symmetric matrix of degree $m$.
Then the multiplicity $m_{\rho,\CM}$ of $\pi_{\rho,\CM}$ in the
regular representation $\pi_{\G_*^J,\rho}$ of $G^J$ in
$L^2_0\left(\G_*^J\backslash G^J; V_{\rho}\right)$ is equal to the
dimension of $J_{\rho,\CM}^{ \textrm{cusp}}(\G_*),$ that is,
\begin{equation*}
m_{\rho,\CM}= \dim_\BC J_{\rho,\CM}^{ \textrm{cusp}}(\G_*).
\end{equation*}
\end{theorem}
\vskip 0.2cm \noindent $ \textit{Idea of Proof.}$ Let $R^J$ be the
right regular representation of $G^J$ on
$L^2_0\left(\G_*^J\backslash G^J; V_{\rho}\right)$ given by
\begin{equation*}
R^J(g_0)\varphi(g)=\varphi(gg_0),\quad g_0,g\in G^J,\ \varphi\in
L^2_0\left(\G_*^J\backslash G^J; V_{\rho}\right).
\end{equation*}

Using the cuspidality condtion on $R^J$, we can prove the
following Step I. \vskip 0.2cm \noindent $ \textbf{Step I:}$ $R^J\
\textit{is completely reducible, and each irreducible component
occurs only a finite }$\par\ \ \ \ \ \ \ \  $ \textit{number of
times in it. }$ \vskip 0.2cm We let $T_{\rho,\CM}^+$ be the
discrete series representation of $G^J$ of lowest weight
$(\rho,\CM)$ such that
\begin{equation*}
T_{\rho,\CM}^+(k,\k)\,v_{\rho,\CM}=\,e^{2\pi i\, \textrm{tr}(\CM
\k)}\rho(k)\,v_{\rho,\CM},\ \quad k\in K,\ \k\in \BR^{(m,m)} \
\textrm{with}\ \k=\,^t\k
\end{equation*}

\noindent for a lowest weight vector $v_{\rho,\CM}.$  Let
\begin{equation*}
L^2_0\left(\G_*^J\backslash G^J; V_{\rho}\right)=\,\oplus\,
{\mathcal H}_{\alpha}
\end{equation*}

\noindent be a decomposition of $R^J$ into irreducible
representation spaces. In each ${\mathcal H}_{\alpha}$ isomorphic
to the representation space of $T_{\rho,\CM}^+$ there is an
analytic lowest weight vector $\Phi_{\rho,\CM,\alpha}$ satisfying
\begin{equation*}
\Phi_{\rho,\CM,\alpha}(g\cdot (k,\k))=\,e^{2\pi i\,
\textrm{tr}(\CM \k)}\rho(k)\Phi_{\rho,\CM,\alpha}(g),
\end{equation*}

\noindent where $g\in G^J,\ k\in K,\ \k\in\BR^{(m,m)}$ with
$\k=\,^t\k.$ We can show that $\Phi_{\rho,\CM,\alpha}$ is the
lifting of a certain cusp form in $J_{\rho,\CM}^{
\textrm{cusp}}(\G_*)$. Thus we get the following Step II. \vskip
0.2cm \noindent $ \textbf{Step II:}$
\begin{equation*}
m_{\rho,\CM}\leq \ \dim_{\BC}J_{\rho,\CM}^{ \textrm{cusp}}(\G_*).
\end{equation*}

Let $f$ be a cusp form in $J_{\rho,\CM}^{ \textrm{cusp}}(\G_*)$.
We can show that the lifting $\Phi_f$ of $f$ to $G^J$ is a lowest
weight vector of a discrete series representation
$T_{\rho,\CM}^+$. Using some properties of $T_{\rho,\CM}^+$ and a
decomposition of $R^J$, we can prove the following Step III.
\vskip 0.2cm \noindent $ \textbf{Step III:}$
\begin{equation*}
m_{\rho,\CM}\geq \ \dim_{\BC}J_{\rho,\CM}^{ \textrm{cusp}}(\G_*).
\end{equation*}

Combining Step I,\,II,\,III, we obtain the duality theorem. The
complete proof will appear in the future. \hfill $\square$

\vskip 0.21cm \noindent {\bf Remark 5.4.} Berndt and B{\"o}cherer
[17] proved the duality theorem for the Jacobi group in the case
$n=m=1.$

\vskip 0.5cm
\begin{center}
{\bf 5.8. Real Analytic Eisenstein Series}
\end{center}

\vskip 0.3cm In the case $n=m=1,$ T. Arakawa [1] introduced the
following real analytic Eisenstein series
\begin{equation*}
E_{k,m,r}((\tau,z),s)\ (\tau\in\BH_1,\ z\in\BC)
\end{equation*}

\noindent for each integer $r$ with $r^2\equiv 0\,( \textmd{mod}\
4m)$ and $s\in\BC$ $(m\in \BZ^+,\ k\in \BZ)$\,:

\begin{align*}
\hskip 2cm E_{k,m,r}((\tau,z),s)
&=& \sum_{\g\in\G^J_{\infty}\ba \G^J} e^m\left( \la^2\tau+2\la z-{ {c(z+\la \tau+\mu)^2}\over {c\tau+d}}\right)\hskip 5cm\\
& &\ \ \ \times (c\tau+d)^{-k}\phi_{r,s}\left( M\circ \tau,{
{z+\la \tau+\mu}\over {c\tau+d}}\right),\hskip 4.5cm
\end{align*}

\noindent where $\G^J=SL(2,\BZ)\ltimes H_\BZ$,\ $e^m(a):=e^{2\pi
ima}$,
\begin{equation*}
\g=\left( \begin{pmatrix} a & b \\ c & d
\end{pmatrix},(\la,\mu,\kappa)\right)\in \G^J,
\end{equation*}
\begin{equation*}
\G^J_{\infty}=\left\{ \left( \begin{pmatrix} 1 & n \\ 0 & 1
\end{pmatrix},(0,\mu,\kappa)\right)\, \Big|\
n,\mu,\kappa\in\BZ\right\},
\end{equation*}
\begin{equation*}
\phi_{r,s}(\tau,z):=e^m\left( {{r^2\tau}\over {4m^2}}+{ {rz}\over
m}\right)( \textmd{Im}\,\tau)^{s-{\frac 12}(k-{\frac 12})}.
\end{equation*}

\noindent It can be seen that $E_{k,m,r}((\tau,z),s)$ converges
absolutely for $ \textmd{Re}\,s>{\frac 54}$ and that if $k>3$ and
it is evaluated at $s={\frac 12}(k-{\frac 12})$, it coincides with
the holomorphic Jacobi form mentioned above. Arakawa proved that
under the assumption that $m$ is a square-free positive integer
and $k$ is even, the above Eienstein series satisfies the {\it
functional equation}
\begin{equation*}
E_{k,m,r}((\tau,z),1-s)=\Phi(1-s)E_{k,m,r}((\tau,z),s)
\end{equation*}

\noindent for a finite number of special integers $r$ and a
certain function $\Phi(s)$ involving the Riemann zeta function and
the Gamma function.

\vskip 0.5cm
\begin{center}
{\bf 5.9. Remark on Spectral Analysis on $\G_{n,m}\ba\BH_{n,m}$}
\end{center}
We fix two positive integers $m$ and $n$ throughout this section.
\vskip 0.1cm For an element $\Om\in \BH_n,$ we set
\begin{equation*}
L_{\Om}:=\BZ^{(m,n)}+\BZ^{(m,n)}\Om\end{equation*} We use the
notations in the subsection 5.3. It follows from the positivity of
$\text{Im}\,\Om$ that the elements $f_{kl},\,h_{kl}(\Om)\,(1\leq
k\leq m,\ 1\leq l\leq n)$ of $L_{\Om}$ are linearly independent
over $\BR$. Therefore $L_{\Om}$ is a lattice in $\BC^{(m,n)}$ and
the set $\left\{\,f_{kl},\,h_{kl}(\Om)\,|\ 1\leq k\leq m,\ 1\leq
l\leq n\, \right\}$ forms an integral basis of $L_{\Om}$. We see
easily that if $\Om$ is an element of $\BH_n$, the period matrix
$\Om_*:=(I_n,\Om)$ satisfies the Riemann conditions (RC.1) and
(RC.2)\,: \vskip 0.1cm (RC.1) \ \ \ $\Om_*J_n\,^t\Om_*=0\,$;
\vskip 0.1cm (RC.2) \ \ \ $-{1 \over
{i}}\Om_*J_n\,^t{\overline{\Om}}_*
>0$.

\vskip 0.2cm \noindent Thus the complex torus
$A_{\Om}:=\BC^{(m,n)}/L_{\Omega}$ is an abelian variety.  \vskip
0.2cm It might be interesting to investigate the spectral theory
of the Laplacian $\Delta_{n,m}$ on a fundamental domain ${\mathcal
F}_{n,m}$. But this work is very complicated and difficult at this
moment. It may be that the first step is to develop the spectral
theory of the Laplacian $\Delta_{\Omega}$ on the abelian variety
$A_{\Omega}.$ The second step will be to study the spectral theory
of the Laplacian $\Delta_n$\,(see (5.7)) on the moduli space
$\Gamma_n\backslash \BH_n$ of principally polarized abelian
varieties of dimension $n$. The final step would be to combine the
above steps and more works to develop the spectral theory of the
Lapalcian $\Delta_{n,m}$ on ${\mathcal F}_{n,m}.$ Maass-Jacobi
forms play an important role in the spectral theory of
$\Delta_{n,m}$ on ${\mathcal F}_{n,m}.$ In this subsection, we
deal only with the spectral theory $\Delta_{\Omega}$ on
$L^2(A_{\Omega}).$

 \vskip 0.1cm We fix
an element $\Om=X+iY$ of $\BH_n$ with $X=\text{Re}\,\Om$ and
$Y=\text{Im}\, \Om.$ For a pair $(A,B)$ with $A,B\in\BZ^{(m,n)},$
we define the function $E_{\Om;A,B}:\BC^{(m,n)}\lrt \BC$ by
\begin{equation*}
E_{\Om;A,B}(Z)=e^{2\pi i\left( \textrm{tr}\,(\,^tAU\,)+\,
\textrm{tr}\, ((B-AX)Y^{-1}\,^tV)\right)},\end{equation*} where
$Z=U+iV$ is a variable in $\BC^{(m,n)}$ with real $U,V$. \vskip
0.1cm\noindent
\begin{lemma} For any $A,B\in \BZ^{(m,n)},$ the function
$E_{\Om;A,B}$ satisfies the following functional equation
\begin{equation*}
E_{\Om;A,B}(Z+\la \Om+\mu)=E_{\Om;A,B}(Z),\quad
Z\in\BC^{(m,n)}\end{equation*} for all $\la,\mu\in\BZ^{(m,n)}.$
Thus $E_{\Om;A,B}$ can be regarded as a function on $A_{\Om}.$
\vskip 0.1cm \end{lemma}
\begin{proof}
We write $\Om=X+iY$ with real $X,Y.$ For any
$\la,\mu\in\BZ^{(m,n)},$ we have
\begin{align*}
E_{\Om;A,B}(Z+\la\Om+\mu)&=E_{\Om;A,B}((U+\la X+\mu)+i(V+\la Y))\\
&=e^{ 2\pi i \left\{\,\textrm{tr}\,(\,^t\!A(U+\la
X+\mu))+\,\textrm{tr}\,((B-AX)Y^{-1}\,^t\!(V+\la Y))\,\right\} }\\
&=e^{ 2\pi i \left\{\,\textrm{tr}\,(\,^t\!AU+\,^t\!A\la
X+\,^t\!A\mu)
+\,\textrm{tr}\,((B-AX)Y^{-1}\,^tV+B\,^t\la-AX\,^t\la) \right\} }\\
&=e^{2\pi i \left\{\,\textrm{tr}\,(\,^t\!AU)\,+\,\textrm{tr}\,((B-AX)Y^{-1}\,^tV)\right\} }\\
&=E_{\Om;A,B}(Z).\end{align*} Here we used the fact that
$^t\!A\mu$ and $B\,^t\la$ are integral. \end{proof}
\newcommand\AO{A_{\Omega}}
\newcommand\Imm{\text{Im}}
\begin{lemma}
The metric
$$ds_{\Om}^2=\textrm{tr}\left((\textrm{Im}\,\Om)^{-1}\,\,^t(dZ)\,d{\overline Z})\,\right)$$
is a K{\"a}hler metric on $A_{\Om}$ invariant under the action
(5.15) of $\G^J=Sp(n,\BZ)\ltimes H_{\BZ}^{(m,n)}$ on $(\Om,Z)$
with $\Om$ fixed. Its Laplacian $\Delta_{\Om}$ of $ds_{\Om}^2$ is
given by
\begin{equation*}
\Delta_{\Om}=\,\textrm{tr}\left( (\Imm\,\Omega)\,{ {\partial}\over
{\partial {Z}} }\,^t\!\left(  {{\partial}\over {\partial
{\overline Z}}} \right)\,
 \right). \end{equation*}
\end{lemma}
\begin{proof} Let ${\tilde \gamma}=(\gamma,(\la,\mu;\kappa))\in
\Gamma^J$ with $\gamma=\begin{pmatrix} A & B\\ C & D
\end{pmatrix}\in Sp(n,\BZ)$ and $({\tilde \Omega},{\tilde
Z})={\tilde \gamma}\cdot (\Omega,Z)$ with $\Omega\in {\mathbb
H}_n$ fixed. Then according to [35],\,p.\,33,
$$\Imm\,\g \cdot\Omega=\,^t(C{\overline
\Om}+D)^{-1}\,\Imm\,\Om\,(C\Om+D)^{-1}$$ and by (5.15),
$$d{\tilde Z}=dZ\,(C\Om+D)^{-1}.$$
Therefore
\begin{eqnarray*}
& & (\Imm\,{\tilde\Om})^{-1}\,^t(d{\tilde Z})\,d{\overline{\tilde Z}}  \\
&=&(C{\overline\Om}+D)\, (\Imm\,\Om)^{-1}\,^t(C
{\Om}+D)\,^t(C{ \Om}+D)^{-1}\,^t(d{ Z})\,d{\overline Z}\,(C{\overline \Om}+D)^{-1} \\
&=& (C{\overline\Om}+D)\,(\Imm\,\Om)^{-1}\,^t(dZ)\, d{\overline
Z}\,(C{\overline \Om}+D)^{-1}.
\end{eqnarray*}
The metric $ds_{iI_n}=\textrm{tr} (dZ\,^t(d{\overline Z}))$ at
$Z=0$ is positive definite. Since $G^J$ acts on ${\mathbb
H}_{n,m}$ transitively, $ds^2_{\Om}$ is a Riemannian metric for
any $\Om\in {\mathbb H}_n.$ We note that the differential operator
$\Delta_{\Om}$ is invariant under the action of $\Gamma^J.$ In
fact, according to (5.15),
$${ {\partial}\over {\partial {\tilde Z}} }=(C\Om +D)\,{
{\partial}\over {\partial Z} }.$$ Hence if $f$ is a differentiable
function on $A_{\Om}$, then
\begin{eqnarray*}
& & \Imm\,{\tilde \Omega}\, { {\partial}\over {\partial {\tilde
Z}} }
\,^t\!\left( { {\partial f}\over {\partial \overline{\tilde Z} }  } \right) \\
&=&\,^t(C{\overline\Om}+D)^{-1}\,(\Imm\,\Om)\,(C{\Om}+D)^{-1} (C{
\Om}+D)\, { {\partial}\over {\partial Z} }\,^t\!\left(
(C{\overline\Om}+D){ {\partial f}\over {\partial \overline Z} }
\right)\\
&=&\,^t(C{\overline\Om}+D)^{-1} \, \Imm\,\Omega\,{ {\partial}\over
{\partial Z} }\,^t\!\left( { {\partial f}\over {\partial \overline
Z} } \right)(C{\overline\Om}+D).
\end{eqnarray*}
Therefore
\begin{equation*}
\textrm{tr}\left(\Imm\,{\tilde \Omega}\, { {\partial}\over
{\partial {\tilde Z}} } \,^t\!\left( { {\partial }\over {\partial
\overline{\tilde Z} }  } \right) \right)=\,\textrm{tr}\left(\,
\Imm\,\Omega\,{ {\partial}\over {\partial Z} }\,^t\!\left( {
{\partial f}\over {\partial \overline Z} } \right) \right).
\end{equation*}

 By the induction on $m$, we can compute the Laplacian
$\Delta_{\Om}.$

\end{proof}
  \vskip 0.1cm We
let $L^2(\AO)$ be the space of all functions $f:\AO\lrt\BC$ such
that
$$||f||_{\Om}:=\int_{\AO}|f(Z)|^2dv_{\Om},$$
where $dv_{\Om}$ is the volume element on $\AO$ normalized so that
$\int_{\AO}dv_{\Om}=1.$ The inner product $(\,\,,\,\,)_{\Om}$ on
the Hilbert space $L^2(\AO)$ is given by
\begin{equation}
(f,g)_{\Om}:=\int_{\AO}f(Z)\,{\overline{g(Z)} }\,dv_{\Om},\quad
f,g\in L^2(\AO).\end{equation}
\begin{theorem}
The set $\left\{\,E_{\Om;A,B}\,|\ A,B\in\BZ^{(m,n)}\,\right\}$ is
a complete orthonormal basis for $L^2(\AO)$. Moreover we have the
following spectral decomposition of $\Delta_{\Om}$:
$$L^2(\AO)=\oplus_{A,B\in \BZ^{(m,n)}}\BC\cdot E_{\Om;A,B}.$$
\end{theorem}

\noindent $\textrm{Proof.}$ The complete proof can be found in
[54] \hfill $\square$

\vskip 0.5cm
\begin{center}
{\bf 5.10. Harmonic Analysis on $\BH_{n,m}$ and $\G_*^J\backslash
G^J$}
\end{center}

\vskip 0.3cm \indent Since $\BH_n=G/K$ is a Riemannian symmetric
space of noncompact type, the Plancherel formula for $L^2(\BH_n)$
can be obtained explicitly by the work of Harish-Chandra.
$L^2(\BH_n)=L^2(G/K)$ has no discrete series representation. Since
the Hilbert space $L^2(\BH_{n,m})$ is isomorphic to
$L^2(\BH_n)\otimes L^2(\BR^{(m,n)})$, the Plancherel formula for
$L^2(\BH_{n,m})$ can be obtained.

\vskip 0.3cm \indent Let $\G_*$ be an arithmetic subgroup of the
Siegel modular group $\G_n=Sp(n,\BZ).$ From now on, for brevity,
we set
\begin{equation*}
L^2:=L^2\left(\G_*^J\backslash G^J\right).
\end{equation*}

\noindent Then it is decomposed into
\begin{equation*}
L^2=L^2_0\oplus L^2_{ \textrm{res}}\oplus L^2_c.
\end{equation*}

\noindent The continuous part $L^2_c$ can be understood by the
theory of Eisenstein series initiated by Robert Langlands [33].
For instance, in the case $n=m=1,\ L^2_c$ was described by R.
Berndt [16,\,18] in some detail.

\vskip 0.2cm At this moment it is hard to describe the cuspidal
part $L_0^2$ in an explicit way. I propose the following problem.

\vskip 0.3cm \noindent {\bf Problem E\,:} Find the trace formula
for the Jacobi group $G^J.$

\vskip 0.3cm \noindent {\bf Acknowledgement\,:} I sincerely thank
Erik van den Ban who gave useful comments on the first draft.
These comments were helpful in improving this article. The referee
pointed out mistakes and presented many useful comments which made
this paper more readable. I would like to give my hearty thanks to
Erik van den Ban and the referee for their valuable comments.

\end{section}


\end{document}